\documentclass[12pt,twoside,a4paper]{amsart}
\usepackage{amssymb}
\usepackage{latexsym}
\usepackage{amsmath}
\usepackage{euscript}
\usepackage{graphics}

   \def\sH{{\mathfrak H}}   
   \def\sK{{\mathfrak K}}   \def\sL{{\mathfrak L}}
   \def\sN{{\mathfrak N}}

      \def\dC{{\mathbb C}}

      \def\dR{{\mathbb R}}

   \def\cH{{\mathcal H}}   
\def\cJ{{\mathcal J}}   \def\cK{{\mathcal K}}   
      
\def\cP{{\mathcal P}}

\def\bR{{\boldmath R}}  

\def\h#1{{{\hat #1} }}
\def\wt#1{{{\widetilde #1} }}
\def\wh#1{{{\widehat #1} }}

\def\bm\chi{\mbox{\boldmath$\chi$}}

\def\RE{{\rm Re\,}}
\def\IM{{\rm Im\,}}

\def\ran{{\rm ran\,}}
\def\cran{{\rm \overline{ran}\,}}
\def\dom{{\rm dom\,}}
\def\mul{{\rm mul\,}}
\def\cmul{{\rm \overline{mul}\,}}
\def\cdom{{\rm \overline{dom}\,}}
\def\clos{{\rm clos\,}}
\def\col{{\rm col\,}}
\def\dim{{\rm dim\,}}

\let\xker=\ker \def\ker{{\xker\,}}

\def\cspan{{\rm \overline{span}\, }}

\def\cmr{{\dC \setminus \dR}}
\def\uphar{{\upharpoonright\,}}

\DeclareMathOperator{\hplus}{\, \widehat + \,}

\newtheorem{theorem}{Theorem}[section]
\newtheorem{proposition}[theorem]{Proposition}
\newtheorem{corollary}[theorem]{Corollary}
\newtheorem{lemma}[theorem]{Lemma}

\theoremstyle{definition}
\newtheorem{example}[theorem]{Example}
\newtheorem{remark}[theorem]{Remark}
\newtheorem{definition}[theorem]{Definition}

\numberwithin{equation}{section}

\begin{document}

\title[Boundary relations and generalized resolvents]
{Boundary relations and generalized resolvents\\ of symmetric
operators}
\author{Vladimir Derkach}
\author{Seppo Hassi}
\author{Mark Malamud}
\author{Henk de Snoo}

\address{Department of Mathematics \\
Donetsk State University \\
Universitetskaya str. 24 \\
83055 Donetsk \\
Ukraine} \email{derkach.v@gmail.com}
\address{Department of Statistics \\
University of Helsinki \\
PL 54, 00014 Helsinki \\
Finland}
\email{hassi@cc.helsinki.fi}
\address{Department of Mathematics \\
Donetsk State University \\
Universitetskaya str. 24 \\
83055 Donetsk \\
Ukraine} \email{mmm@telenet.dn.ua}
\address{Department of Mathematics \\
University of Groningen \\
Postbus 800, 9700 AV Groningen \\
Nederland}
\email{desnoo@math.rug.nl}
\date{\today}
\thanks{The present research was supported
by the Academy of Finland (project 116842).}

\begin{abstract} The Kre\u{\i}n-Naimark formula provides a
parametrization of all selfadjoint exit space extensions of a, not
necessarily densely defined, symmetric operator, in terms of maximal
dissipative (in $\dC_+$) holomorphic linear relations on the
parameter space (the so-called Nevanlinna families). The new notion
of a boundary relation makes it possible to interpret these
parameter families as Weyl families of boundary relations and to
establish a simple coupling method to construct the generalized
resolvents from the given parameter family. The general version of
the coupling method is introduced and the role of boundary relations
and their Weyl families for the Kre\u{\i}n-Naimark formula is
investigated and explained.
\end{abstract}

\maketitle

\section{Introduction}

Let $\sH$ be a separable Hilbert space and let $A$ be a not
necessarily densely defined closed symmetric operator or relation in
$\sH$ with equal defect numbers $n_+(A)=n_-(A) \le \infty$. Denote
by $A^*$ the adjoint linear relation of $A$. The Kre\u{\i}n-Naimark
formula
\begin{equation}\label{01}
\bR_\lambda :=P_\sH(\wt A-\lambda)^{-1}\uphar\sH
=(A_0-\lambda)^{-1}-\gamma(\lambda)(M(\lambda)
     +\tau(\lambda))^{-1}\gamma(\bar{\lambda})^*,
     \quad \lambda\in\rho(A_0)\cap \rho(\wt A),
\end{equation}
establishes  a bijective  correspondence between the set of all
selfadjoint  (canonical and exit space) extensions $\wt A$ of $A$
and the set all Nevanlinna families $\tau({\lambda})$.  Here
$A_0=A^*_0$ is a fixed canonical extension of $A,\gamma(\lambda)$ is
the so-called $\gamma$-field, and $M(\lambda)$ is a $Q$-function of
the pair $\{A,A_0\}$. The correspondence in \eqref{01} will also be
indicated by the notation ${\wt A}=A^{(\tau)}$. The
Kre\u{\i}n-Naimark formula  plays an important role in the extension
theory of the operator $A$ (see \cite{AG, AGHH, AKur, DHMS1, DM1,
DM2, GKMT} and references therein) and its numerous applications to
classical interpolation problems (\cite{Ber, KraKre, Kr0, Kr2, KL1,
Naj1, DM1, DM2}), boundary value problems (\cite{DLS, DLS1, GG}) and
different type of physical problems (see \cite{AGHH, AKur, BruGeil,
BruPan, Pav, Pos} and references therein).

During the last two decades a new approach to the extension theory,
based on the concepts of boundary triplets and the corresponding
Weyl functions, has been developed. Recall the basic definitions.

\begin{definition}\label{een}\cite{GG}
A collection  ${\Pi}=\{ \cH, \Gamma_0, \Gamma_1\}$  consisting of a
Hilbert space $\cH$ with $\dim \cH=n_\pm(A)$ and two linear mappings
$\Gamma_0$ and $\Gamma_1$ from $A^*$ to $\cH$, is said to be a
boundary triplet for $A^*$ if
\begin{enumerate}
\def\labelenumi{\rm (\roman{enumi})}
\item[(BT1)] the abstract Green's identity holds
\begin{equation}\label{Greendef}
(f', g)-(f, g') = (\Gamma_1 \wh f, \Gamma_0 \wh g)_{\cH}-
(\Gamma_0f, \Gamma_1g)_{\cH}, \quad \wh f=\{f,f'\},\,\, \wh
g=\{g,g'\} \in A^*;
\end{equation}
\item[(BT2)] the linear mapping
$\Gamma:=\{\Gamma_0,\Gamma_1\}:\,A^* \to {\cH}\oplus{\cH}$ is
surjective.
\end{enumerate}
\end{definition}

The mappings $\Gamma_0$ and $\Gamma_1$ induce two selfadjoint
extensions $A_0=\ker \Gamma_0$ and $A_1=\ker \Gamma_1$ of $A$.  In
\cite{DM0,DM1} the concept of a Weyl function was associated to an
ordinary boundary triplet as an abstract version of the $m$-function
appearing in boundary value problems for differential operators.

\begin{definition} (\cite{DM0,DM1})
\label{W000} Let $\Pi=\{\cH,\Gamma_0,\Gamma_1\}$ be a boundary
triplet for $A^*$. The operator-valued function $M(\lambda)$ defined
by
\begin{equation}\label{0.2}
 \Gamma_1f_{\lambda}=M(\lambda)\Gamma_0f_{\lambda},
 \quad
 f_{\lambda}\in \sN_{\lambda}:=\ker(A^*-\lambda),
\quad \lambda\in\cmr,
\end{equation}
is said to be the \textit{Weyl function}, corresponding to the
triplet $\Pi$.
\end{definition}

A connection between the approach via boundary triplets and the
Kre\u{\i}n-Naimark theory of generalized resolvents has been
established in \cite{DM1,DM2}. It was shown  that all objects in
\eqref{01} can be expressed in terms of boundary triplets.  In fact,
one has
\begin{equation}
A_0=\ker\Gamma_0, \quad   \gamma(\lambda)=(\Gamma_0\lceil\frak
N_{\lambda})^{-1},\quad \{\Gamma_0,\Gamma_1\}R_{\lambda}f\in
-\tau(\lambda), \quad f\in\sH.
\end{equation}
In formula \eqref{01} the Weyl function $M({\lambda})$ is always a
uniformly strict Nevanlinna function, whereas the parameter
$\tau({\lambda})$ is an arbitrary Nevanlinna family. It is known
that any uniformly strict Nevanlinna function is the Weyl function
in the sense of Definitions \ref{een} and \ref{W000}.

If the parameter $\tau({\lambda})$ in formula \eqref{01}
 is a uniformly strict Nevanlinna  function one
can use the inverse problem for Weyl functions in order to construct
the exit space extension ${\widetilde A}={\widetilde A}^{(\tau)}$
connected with $\bR_\lambda$, via \eqref{01}.   This construction
will be briefly recalled; cf. \cite{DHMS1}. Let $S_1:=A$ and let
$S_2$ be a symmetric operator in a Hilbert space ${\frak H}_2$ such
that $\tau(\lambda)$ is the Weyl function of $S_2$ corresponding to
a boundary triplet $\Pi_2 = \{\mathcal H, \chi_0, \chi_1\}$. Then
the linear relation
\begin{equation}\label{coupl01}
 \wt A=\left\{\,\wh f_1\oplus\wh f_2\in A^*\oplus T_2:\, \Gamma_0\wh
f_1-\chi_0\wh f_2=\Gamma_1\wh f_1+\chi_1\wh f_2=0\,\right\}
\end{equation}
is a selfadjoint (exit space) extension of  $S_1\oplus S_2$ and
satisfies  equation \eqref{01}.  Unfortunately this coupling
approach was restricted to uniformly strict Nevanlinna  functions
$\tau(\lambda)$.  In order to extend this  method to arbitrary
Nevanlinna families the new concepts of boundary relations and their
Weyl families were introduced by the authors in \cite{DHMS2},
\cite{Trans}. These concepts generalize the notions of the boundary
triplet and the corresponding  Weyl functions. In \cite{DHMS2} it
was proved that every Nevanlinna family $\tau(\lambda)$  can be
realized as the Weyl family of a boundary relation. The main purpose
of the paper is to show that, due to this new inverse result, the
coupling construction in \eqref{coupl01} can be extended to the case
of any Nevanlinna family.

The paper is organized as follows. In Section \ref{preli} the basic
notions are introduced and various preliminary results are
established. In particular, some new and useful facts on unitary
relations in Kre\u{\i}n spaces are presented, for instance,
concerning the composition of unitary relations; see
Theorem~\ref{product}. In Section 3 the notion of boundary relations
for $S^*$, the corresponding Weyl families, orthogonal couplings,
and $J$-unitary transformations of boundary relations are discussed.
In particular, it is shown that if two boundary relations $\wt
\Gamma$ and $\Gamma$ are connected by means of a standard
$J$-unitary operator $W$ via $\wt \Gamma=W\Gamma$, then the
corresponding Weyl families are connected by means of Shmulyan's
transform. Besides, the following equality is derived
      \begin{equation}
\dim{\cH} - n_{\pm}(A) = \mul\Gamma,
     \end{equation}
showing, in particular, that the equality $\dim{\mathcal
H}=n_{\pm}(A)$ is true if and only if $\Gamma$ is single-valued.

In Section 4 the connection between boundary relation
$({\widetilde\Gamma},{\mathcal H}^2)$ and ordinary boundary triplet
$\Pi=\{{\cH},\Gamma_0,\Gamma_1\}$ for $S^*$ is investigated. In
particular, it is shown (see Proposition~\ref{WGamma}) that formula
        \begin{equation}\label{Wconnect01}
    \wt \Gamma=W\Gamma
      \end{equation}
establishes a bijective correspondence between the set of all
boundary relations for $S^*$ and the set of unitary relations $W$ in
$(\cH^2,J_\cH)$ for which $\ker W=\{0\}$. Observe, that formula
\eqref{Wconnect01} leads to another (equivalent) definition of a
boundary relation at the expense of extending the group of
$J$-unitary operators in \eqref{Wconnect01} to  a (wider) set of
$J$-unitary relations $W$ with $\ker W=\{0\}$. In this section also
generalized boundary triplets as well as boundary triplets whose
Weyl functions take values in $[\cH]$ are investigated.

In Section 5 there are some general transformation results
concerning boundary relations $\Gamma:\sH^2\to\cH^2$ for $S^*$ whose
Weyl family $M(\lambda)$ belongs to the class $R[\cH]$, that is
$M(\cdot)$ is the Weyl function with values in $[\cH].$ In this case
an arbitrary orthogonal decomposition $\cH=\cH_1 \oplus \cH_2\ $ of
$\cH$ induces the corresponding block operator representation
       \begin{equation} \label{decom02}
  M(\lambda)=(M_{ij}(\lambda))_{i,j=1}^2.
       \end{equation}
of $M(\cdot)$. It is shown how one can identify intermediate closed
symmetric extensions $H$ of $A$ and associate boundary relations for
$H^*$, such that the corresponding Weyl function is a given
transform of the blocks of $(M_{ij}(\lambda))$ including, for
instance, linear combinations of $M_{ij}(\lambda)$ and Schur
complements. In particular, there appear induced boundary relations
$\wt \Gamma$ for $H^*$ whose Weyl function $\wt M(\cdot)$ equals
either to $M_{11}+M_{22}$ or to $(M_{11} + M_{22})^{-1}.$

Similar results for ordinary boundary triplets ${\Pi}=\{ \cH,
\Gamma_0, \Gamma_1\}$ for $A^*$ have earlier been published in our
previous paper \cite{DHMS1}. However, the present generalizations
are needed here for applications involving generalized resolvents.

In Section 6 the coupling method from \cite{DHMS1}, as briefly
described above, is extended to the case of arbitrary Nevanlinna
families $\tau(\cdot)$. This approach leads to new results and
further geometric insight into various questions on this area. In
the coupling method the selfadjoint exit space extension $\wt A$ in
$\wt\sH\supset \sH$ is constructed by means of a boundary triplet of
$A^*$, whose Weyl function is $M(\cdot)$, and a boundary relation
that corresponds to the family $\tau(\cdot)\in {\wt R}(\cH)$. The
coupling method makes it possible to treat the families
$\tau(\cdot)$ and $-(\tau(\cdot) + M(\cdot))^{-1}$ appearing in
\eqref{01} as the Weyl families of $S_2 :=\wt A \cap(\wt\sH
\ominus\sH)^2$ (see formula \eqref{ST}) and some intermediate
extension of $A$ (see formula \eqref{couplH1}), respectively.

In Section 7 coupling method is applied to give a complete solution
to the problem of $M$-admissibility (cf. \cite{DHMS1}). Recall, that
if $A$ is nondensely defined, then $A^{(\tau)}$ may be either a
(selfadjoint) linear relation or an operator (i.e, a single-valued
linear relation). Based on a coupling construction, the following
simple criterion for $\tau(\cdot)$ to generate an operator
$A^{(\tau)}$ is established:

The Nevanlinna family $\tau(\cdot)$ in \eqref{01} corresponds to an
operator $A^{(\tau)}$ (that is, it is  $\Pi$-admissible) if and only
if the following two conditions are satisfied:
\begin{equation} \label{Adm01}
  w-\lim_{y \uparrow \infty}
  \frac{(\tau(iy)+M(iy))^{-1}}{y}=0, \quad  \qquad \lim_{y \uparrow \infty}
  \frac{(\tau(iy)^{-1} + M(iy)^{-1})^{-1}}{y}=0.
           \end{equation}
Moreover, results on intermediate extensions given in Section 5 (a
geometric treatment of $(\tau(iy)+M(iy))^{-1}$ as a Weyl function)
allow us to show that if additionally $A_0$ (resp. $A_1$) is an
operator, then $A^{(\tau)}$ is an operator if and only if the first
(resp. the second) of conditions \eqref{Adm01} is satisfied.

In a forthcoming paper the coupling method is applied to the
characterization of the Naimark extensions in terms of
$\tau(\cdot)$.

\section{Preliminaries}\label{preli}

\subsection{Linear relations in Hilbert spaces}

The \textit{Cartesian product} $\sH \times\sH'$ of linear spaces
$\sH$ and $\sH'$ is the set of all ordered pairs (of $1 \times 2$
matrices) $\{f,f'\}$ with $f \in \sH$ and $f' \in \sH'$. Frequently
it will be convenient to denote the Cartesian product $\sH \times
\sH'$ and the elements of it (as $2\times 1$ matrices) by
\[
   \begin{pmatrix} f \\ f' \end{pmatrix}
   \in
   \begin{pmatrix}\sH \\ \sH' \end{pmatrix},
   \quad f \in \sH, \quad f'\in \sH'.
\]
If $\sL \subset \sH$ and $\sL' \subset \sH'$ are linear subspaces
then $\sL \times \sL'$ denotes the Cartesian product of the
subspaces; in agreement with the ordered pairs this product will
also be denoted by $\{\sL,\sL'\}$, or by $\sL^2$ if $\sL=\sL'$.

A linear relation $T$ from $\sH$ to $\sH'$ is a linear subspace of
$\sH \times \sH'$. Systematically a linear operator $T$ will be
identified with its graph. It is convenient to write $T:\sH\to\sH'$
and interpret the linear relation $T$ as a multi-valued linear
mapping from $\sH$ into $\sH'$. If $\sH'=\sH$ one speaks of a linear
relation $T$ in $\sH$.

For a linear relation $T:\sH \to \sH'$ the symbols $\dom T$, $\ker
T$, $\ran T$, and $\mul T$ stand for the domain, kernel, range, and
the multi-valued part, respectively. The inverse $T^{-1}$ is a
relation from $\sH'$ to $\sH$ defined by $\{\,\{f',f\}:\,\{f,f'\}\in
T\,\}$. The adjoint $T^*$ is the closed linear relation from $\sH'$
to $\sH$ defined by (see~\cite{Ben}, \cite{Co})
\begin{equation}
\label{adjo}
 T^*=\{\,\{h,k\} \in \sH' \oplus \sH :\,
     (k,f)_{\sH}=(h,g)_{\sH'}, \, \{f,g\}\in T \,\}.
\end{equation}
The sum $T_1+T_2$ and the componentwise sum $T_1 \wh + T_2$ of two
linear relations $T_1$ and $T_2$ are defined by
\[
\begin{split}
  &T_1+T_2=\{\,\{f, g+h\} :\, \{f,g\} \in T_1, \{f,h\} \in T_2\,\}, \\
  &T_1 \hplus T_2=\{\, \{f+h,g+k\} :\, \{f,g\} \in T_1, \{h,k\}\in
T_2\,\}.
\end{split}
\]
If the componentwise sum is orthogonal it will be denoted by $T_1
\oplus T_2$.
The null spaces of $T-\lambda$, $\lambda \in \dC$, are defined by
\begin{equation}
\label{defect0} \sN_\lambda(T)=\ker (T-\lambda), \quad \wh
\sN_\lambda(T) =\{\, \{f,\lambda f\} \in T :\,
f\in\sN_\lambda(T)\,\}.
\end{equation}
Moreover, $\rho(T)$ ($\hat\rho(T)$) stands for the set of regular
(regular type) points of $T$. The closure of a linear relation $T$
will be denoted by $\clos T$.

The product of linear relations is defined in the standard way.
Some basic facts concerning the product of operators remain valid also for
the product of relations. For instance, the following statement is easy to
check.

\begin{lemma}
Let $\sH_1$, $\sH_2$, and $\sH_3$ be Hilbert spaces and let
$B:\sH_1\to\sH_2$ and $A:\sH_2\to\sH_3$ be linear relations, and let
$C=AB$. Then:
\begin{enumerate}
\def\labelenumi{\rm (\roman{enumi})}
\item $\ker B \subset \ker C$ and $\mul A \subset \mul C$;

\item if $\ker A=\{0\}$, then $\ker B = \ker C$, and
      if $\mul B=\{0\}$, then $\mul A = \mul C$.
\end{enumerate}
\end{lemma}

The next lemma gives some basic facts concerning the inverse, the
product, and the adjoint of linear relations; these facts are well
known for linear operators, the proofs for linear relations are left
to the reader.

\begin{lemma}
\label{basicle} Let $\sH_1$, $\sH_2$, and $\sH_3$ be Hilbert spaces
and let $B:\sH_1\to\sH_2$ and $A:\sH_2\to\sH_3$ be linear relations.
Then:
\begin{enumerate}
\def\labelenumi{\rm (\roman{enumi})}

\item $(A B)^{-1} = B^{-1}A^{-1}$ and $(A^*)^{-1}=(A^{-1})^*$;

\item $(A B)^* \supset B^{*} A^{*}$;

\item if $A \in [\sH_1,\sH_2]$ or $B^{-1}\in[\sH_3,\sH_2]$,
then $(A B)^* = B^{*} A^{*}$.
\end{enumerate}
\end{lemma}

Recall that a linear relation $T$ in $\sH$ is called
\textit{symmetric} (\textit{dissipative}) or \textit{accumulative})
if $\mbox{Im }(h',h)=0$ ($\ge 0$) or $\le 0$, respectively) for all
$\{h,h'\}\in T$. These properties remain invariant under closures.
By polarization it follows that a linear relation $T$ in $\sH$ is
symmetric if and only if $T \subset T^*$. A linear relation $T$ in
$\sH$ is called \textit{selfadjoint} if $T=T^*$, and it is called
\textit{essentially selfadjoint} if $\clos T=T^*$. A dissipative
(accumulative) linear relation $T$ in $\sH$ is called maximal
dissipative (maximal accumulative) if it has no proper dissipative
(accumulative) extensions.

Assume that $T$ is closed. If $T$ is dissipative or accumulative,
then $\mul T \subset \mul T^*$.
In this case the orthogonal decomposition
$\sH=(\mul T)^\perp \oplus \mul T$
induces an orthogonal decomposition of $T$ as
\begin{equation}
\label{eq0}
 T=T_s \oplus T_\infty, \quad T_\infty=\{0\} \times \mul T, \quad
 T_s=\{\, \{f,g\} \in T:\, g \perp \mul T\,\},
\end{equation}
where $T_\infty$ is a selfadjoint relation in $\mul T$
and $T_s$ is an operator in $\sH \ominus \mul T$
with $\cdom T_s=\cdom T=(\mul T^*)^\perp$,
which is dissipative or accumulative.
Moreover, if the relation $T$ is maximal
dissipative or accumulative, then $\mul T=\mul T^*$.
In this case the orthogonal decomposition
$(\dom T)^\perp=\mul T^*$ shows that $T_s$
is a densely defined dissipative or accumulative
operator in $(\mul T)^\perp$, which is maximal
(as an operator).
In particular, if $T$ is a selfadjoint relation,
then there is such a decomposition where $T_s$
is a selfadjoint operator
(densely defined in $(\mul T)^\perp$).

Let $T$ be a linear relation in a Hilbert space $\sH$.  If  $T$ is
closed, then also the eigenspace $\wh\sN_\lambda(T)$ is closed for
every $\lambda \in \dC$.

\begin{lemma}
\label{defectlem} \cite{HSSW} Let $T$ be a linear relation in $\sH$,
let $H$ be a restriction of $T$ with a nonempty resolvent set, and
assume that $\lambda \in \rho(H)$. Then $H$ is closed and
\begin{equation}
\label{vonneu1N} T=H \hplus \wh\sN_\lambda(T).
\end{equation}
\end{lemma}

Observe that if $S$ is a closed symmetric relation and $H$ is a
(maximal symmetric, maximal dissipative) selfadjoint extension, then
\[
 S^*=H \hplus \wh\sN_\lambda(S^*).
\]
In particular, if $S$ is a symmetric relation in $\sH$, then $H:=S
\hplus \wh \sN_\lambda(S^*)$ is a restriction of $S^*$ and,
moreover, if $S$ is closed,
\[
 \bar\lambda \in \rho(H), \quad \lambda\in\cmr.
\]
Hence, the adjoint relation $S^*$ of a closed symmetric linear
relation $S$ in a Hilbert spaces $\sH$ can be decomposed via the von
Neumann formula:
\begin{equation}
\label{JvonN0}
 S^*=S \hplus
     \wh\sN_\lambda(S^*) \hplus \wh\sN_{\bar{\lambda}}(S^*),
 \quad \lambda \in \cmr,
 \quad
 \mbox{direct sums,}
\end{equation}
where $\wh\sN_\lambda(S^*)$ is defined as in \eqref{defect0}.
When $\lambda =\pm i$ the decomposition \eqref{JvonN0}
is orthogonal:
\begin{equation}
\label{JvonN}
 S^*=S\oplus\wh\sN_i(S^*)\oplus\wh\sN_{-i}(S^*),
\end{equation}
where the orthogonality is with respect to the inner product
topology in $S^*$, cf. \cite{Ben}, \cite{Co}. A symmetric linear
relation $S$ is called \textit{simple} if there is no nontrivial
orthogonal decomposition of the Hilbert space $\sH=\sH_1
\oplus\sH_2$ and no corresponding orthogonal decomposition $S=S_1
\oplus S_2$ with $S_1$ a symmetric relation in $\sH_1$ and $S_2$ a
selfadjoint relation in $\sH_2$. The decomposition \eqref{eq0} for
$S=S_s\oplus S_\infty$ shows that a simple closed symmetric relation
is necessarily an operator. Recall that (cf. e.g. \cite{LT}) a
closed symmetric linear relation $S$ in a Hilbert space $\sH$ is
simple if and only if
\[
 \sH=\cspan \{\, \sN_\lambda(S^*) :\, \lambda \in \cmr\,\}.
\]

\subsection{Linear relations in Kre\u{\i}n spaces}

Recall that a signature operator $j$ in a Hilbert space is a bounded
linear operator such that $j=j^*=j^{-1}$. A signature operator
provides the Hilbert space with a Kre\u{\i}n space structure with
the fundamental symmetry $j$. Let $\sH$ and $\cH$ be Hilbert spaces
with signature operators $j_\sH$ and $j_\cH$, respectively, and
denote the corresponding Kre\u{\i}n spaces by $(\sH,j_\sH)$ and
$(\cH,j_\cH)$. Then the adjoint $T^{[*]}$ of a linear relation $T$
from the Kre\u{\i}n space $(\sH,j_\sH)$ to the Kre\u{\i}n space
$(\cH,j_\cH)$ is given by $T^{[*]}=j_{\sH}T^*j_{\cH}$.

The following result was given in \cite[Proposition~2.2]{Trans} for
the Hilbert space case.

\begin{proposition}
\label{closo} Let $T$ be a closed linear relation from the
Kre\u{\i}n space $(\sH,j_\sH)$ to the Kre\u{\i}n space
$(\cH,j_\cH)$. Then:
\begin{enumerate}
\def\labelenumi{\rm (\roman{enumi})}

\item $\dom T$ is closed if and only if $\dom T^{[*]}$ is closed;

\item $\ran T$ is closed if and only if $\ran T^{[*]}$ is closed.
 \end{enumerate}
\end{proposition}

\begin{proof}
The Kre\u{\i}n space adjoint $T^{[*]}$ of $T$ is connected to the
Hilbert space adjoint $T^*$ via $T^{[*]}=j_{\sH}T^*j_{\cH}$. Hence
it is clear that $\dom T^{[*]}$ ($\ran T^{[*]}$) is closed if and
only if $\dom T^*$ (resp. $\ran T^*$) is closed. Therefore, the
statements follow from \cite[Proposition~2.2]{Trans}.
\end{proof}

\begin{definition}
(i)\ A linear relation $T$ from the Kre\u{\i}n space $(\sH,j_\sH)$
to the Kre\u{\i}n space $(\cH,j_\cH)$ is said to be
\textit{isometric} if $T^{-1}\subset T^{[*]}$ and
\textit{co-isometric} if $T^{[*]}\subset T^{-1}$.

(ii)\ A linear relation $T$ is said to be \textit{unitary} if it is
simultaneously isometric and co-isometric, that is, if
$T^{-1}=T^{[*]}$.
\end{definition}

\begin{lemma} \label{UNITlem}
(\cite{Trans}) Let $T$ be a linear relation from the Kre\u{\i}n
space $(\sH,j_\sH)$ to the Kre\u{\i}n space $(\cH,j_\cH)$. Then:
\begin{enumerate}
\def\labelenumi{\rm (\roman{enumi})}

\item if $T:\sK_1\to\sK_2$ isometric then the inverse
$T^{-1}:\sK_2\to\sK_1$ is isometric and the adjoint
$T^{[*]}:\sK_2\to\sK_1$ is co-isometric;

\item if $T:\sK_1\to\sK_2$ is unitary then the inverse
$T^{-1}$ and the adjoint $T^{[*]}$ are also unitary.
\end{enumerate}
\end{lemma}
\begin{proof}
(i) Since $T$ is isometric, one has $T^{-1}\subset T^{[*]}$. Taking
inverses one obtains
\[
(T^{-1})^{-1}\subset(T^{[*]})^{-1}=(T^{-1})^{[*]}
\]
by Lemma~\ref{basicle}, so that $T^{-1}$ is isometric. Taking
adjoints one obtains
\[
(T^{[*]})^{[*]}\subset(T^{-1})^{[*]}=(T^{[*]})^{-1},
\]
again by Lemma~\ref{basicle}, so that $T^{[*]}$ is co-isometric.

(ii) This statement is clear from $(T^{-1})^{[*]}=(T^{[*]})^{-1}$;
cf. Lemma~\ref{basicle}.
\end{proof}

The following two statements are due to Yu.L.~Shmul'jan \cite{Sh0}.
They can be obtained also directly from the equality
$T^{[*]}=T^{-1}$ and Proposition~\ref{closo}; see also~\cite{Trans}.

\begin{proposition}\label{UNIT}
Let $T$ be a unitary relation from the Kre\u{\i}n space
$(\sH,j_\sH)$ to the Kre\u{\i}n space $(\cH,j_\cH)$. Then:
\begin{enumerate}
\def\labelenumi{\rm (\roman{enumi})}
\item $\dom T$ is closed if and only if $\ran T$ is closed; \item
the following equalities hold:
\begin{equation}\label{eee}
\ker T=(\dom T)^{[\perp]}, \quad \mul T=(\ran T)^{[\perp]}.
\end{equation}
\end{enumerate}
\end{proposition}

A unitary relation $T:(\sH,j_\sH)\to(\cH,j_\cH)$ may be
multi-valued, nondensely defined, or unbounded. The following
characterization is useful.

\begin{lemma} \label{UNITcor}
(\cite{Trans}) Let $T$ be a unitary relation from the Kre\u{\i}n
space $(\sH,j_\sH)$ to the Kre\u{\i}n space $(\cH,j_\cH)$. Then:
\begin{enumerate}
\def\labelenumi{\rm (\roman{enumi})}
\item $T$ is single-valued if and only if $\cran T=\cH$;

\item $T$ is single-valued and densely defined if and only if
$\cran T=\cH$ and $\ker T=\{0\}$;

\item $T$ is single-valued and bounded (not
necessarily densely defined) if and only if $\ran T=\cH$;

\item $T\in[\sH,\cH]$ if and only if $\ran T=\cH$ and $\ker T=\{0\}$.
\end{enumerate}
\end{lemma}

A unitary relation $T$ is the graph of an operator if and only if
its range is dense. In this case it need not be densely defined or
bounded; and if it is bounded it need not be densely defined.

\begin{corollary}\label{UNITcor+}
Let $T$ be a unitary relation from the Kre\u{\i}n space
$(\sH,j_\sH)$ to the Kre\u{\i}n space $(\cH,j_\cH)$. Then
$T\in[\sH,\cH]$ if and only if $T^{-1} \in [\cH,\sH]$.
\end{corollary}

\begin{proof}
Assume $T\in[\sH,\cH]$. Then $\dom T=\sH$ and $\mul T=\{0\}$, or
equivalently, $\ran T^{-1}=\sH$ and $\ker T^{-1}=\{0\}$. By (iv) of
Lemma~\ref{UNITcor} this implies that $T^{-1} \in [\cH,\sH]$.
\end{proof}

Observe that for a unitary relation $T$ from $(\sH,j_\sH)$ to
$(\cH,j_\cH)$, both $T$ and $T^{-1}$ are operators if and only if
$\cdom T=\sH$ and $\cran T=\cH$. Moreover, in this case $\dom T
=\sH$ if and only if $\ran T=\cH$, cf. Proposition \ref{UNIT}, which
also leads to Corollary \ref{UNITcor+}.

\begin{remark}\label{uj}
In the present terminology an operator $T$ is unitary if it
satisfies $T^{-1}=T^{[*]}$. However, the terminology in
\cite[Chapter 2, Definition 5.1 and Corollary 5.8]{AI} is different.
An operator $T$ from the Kre\u{\i}n space $(\sH,j_\sH)$ to the
Kre\u{\i}n space $(\cH,j_\cH)$ is unitary in the sense of M.G.
Kre\u{\i}n (see\cite{AI}), if $\dom T=\sH$, $\ran T=\cH$, and
\begin{equation}\label{ai}
 [Tf,Tf]_\cH=[f,f]_\sH, \quad f \in \sH.
\end{equation}
To see the connection with the present setting, observe that
\eqref{ai} implies by polarization that
\begin{equation}\label{ai+}
 [Tf,Tg]_\cH=[f,g]_\sH, \quad f,g \in \sH.
\end{equation}
The identity \eqref{ai+} shows that $T$ is isometric, i.e., the
graph of $T$ satisfies $T^{-1} \subset T^{[*]}$, and that $\ker
T$=\{0\}. Since $\dom T=\sH$ and $\ran T=\cH$ it follows that $T$ is
unitary, i.e., $T^{-1} = T^{[*]}$, cf. \cite[Proposition
2.5]{Trans}. Moreover, $T \in [\sH,\cH]$ by (iv) of
Lemma~\ref{UNITcor} and then also $T^{-1} \in [\cH,\sH]$ by
Corollary \ref{UNITcor+}. Conversely, if $T$ is a unitary relation,
i.e., $T^{-1} = T^{[*]}$ and $T \in [\sH,\cH]$ (or equivalently
$T^{-1} \in [\cH,\sH]$), then $T$ is an operator satisfying
\eqref{ai}, $\dom T=\sH$, and $\ran T=\cH$ by (iv) of
Lemma~\ref{UNITcor}. Therefore, a unitary relation $T$ is \textit{a
standard unitary operator} (in the sense of M.G. Kre\u{\i}n)
precisely when $T$ in addition belongs to $[\sH,\cH]$, i.e., $T$ is
everywhere defined and single-valued, in which case also
$T^{-1}\in[\cH,\sH]$ is a standard unitary operator. In the present
paper a unitary operator need not belong to $[\sH,\cH]$ and it need
not be even densely defined, in which case $\ker T$ is also
nontrivial; cf. Proposition~\ref{UNIT}.
 \end{remark}

On a finite-dimensional space the set of injective unitary operators
coincides with the set of standard unitary operators.

\begin{corollary} \label{cor1}
Let $\sK_1$ and $\sK_2$ be Kre\u{\i}n spaces, let $T:\sK_1\to\sK_2$
be a unitary operator with $\ker T=\{0\}$, and assume that $\dim
\sK_1<\infty$. Then $\dim \sK_2=\dim \sK_1<\infty$ and $T$ is a
standard unitary operator in $[\sK_1,\sK_2]$.
\end{corollary}

 \begin{proof}
Since $\ker T=\{0\}$ and $\dim \sK_1<\infty$, it follows from
Proposition~\ref{UNIT} that $\dom T=\sK_1$ and that $\ran T$ is
closed; furthermore, since $T$ is single-valued, one has $\ran
T=\sK_2$. Thus $T\in[\sK_1,\sK_2]$ by (iv) of Lemma~\ref{UNITcor},
so that $T$ is a standard unitary operator; cf. Remark~\ref{uj}.
Therefore also $\dim \sK_2=\dim \sK_1<\infty$.
 \end{proof}

It is emphasized that the condition $\ker T=\{0\}$ in
Corollary~\ref{cor1} is essential: as will be seen below (see also
\cite{Trans}) boundary triplets of symmetric operators $S$, $\dom S
\neq \{0\}$, (acting on a finite-dimensional or infinite dimensional
space) are typical examples of bounded unitary operators which are
not standard. They are nondensely defined with a nontrivial kernel
that is equal to $S$.

Unitary relations between Kre\u{\i}n spaces admit a couple of useful
properties under composition. First a result which concerns the
adjoint of the product of linear relations in the case that the
domain or the range of one of the relations is closed; observe that
Lemma~\ref{basicle} is still true in the Kre\u{\i}n space situation.

\begin{lemma}\label{prodlemma}
Let $\sK_j$, $j=0,1,2,3$, be Kre\u{\i}n spaces and let
$S:\sK_1\to\sK_2$ be a closed relation. Then:
\begin{enumerate}
\def\labelenumi{\rm (\roman{enumi})}

\item if $\dom S$ is closed then for every linear relation
$X:\sK_0\to\sK_1$ with $\ran X\subset \dom S$ one has
\[
  (SX)^{[*]}=X^{[*]}S^{[*]};
\]

\item if $\ran S$ is closed then for every linear relation $Y:\sK_2\to\sK_3$
with $\dom Y\subset\ran S$ one has
\[
  (YS)^{[*]}=S^{[*]}Y^{[*]}.
\]
\end{enumerate}
\end{lemma}
\begin{proof}
(i) The inclusion $(SX)^{[*]}\supset X^{[*]}S^{[*]}$ is always
satisfied, cf. (ii) in Lemma~\ref{basicle}. To prove the reverse
inclusion let $\{f,g\} \in (SX)^{[*]}$, so that
\begin{equation}
\label{prod00+}
  [g,h]_{\sK_1}=[f,k]_{\sK_3} \quad\mbox{for all } \{h,k\} \in SX.
\end{equation}
Since the linear relation $SX$ contains the set
\[
 \{\{0,f_0\}:f_0\in \mul S\}
\]
it follows from~\eqref{prod00+} that $[f,f_0]=0$ for all $f_0\in
\mul S$, so that $f\in (\mul S)^{[\perp]}=\cdom S^{[*]}$. Since $S$
is closed and $\dom S$ is closed, also $\dom S^{[*]}$ is closed by
Proposition~\ref{closo}. Hence $f\in \dom S^{[*]}$ and $\{f,f'\} \in
S^{[*]}$ for some $f' \in \sK_1$. Now it suffices to show that
$\{f',g\} \in X^{[*]}$, because then $\{f,g\} \in X^{[*]}S^{[*]}$.
Indeed, for each $\{h,u\} \in X$ there is $u'\in \sK_2$ such that
$\{u,u'\} \in S$, due to the condition $\ran X \subset \dom S$. Then
for all $\{f,f'\} \in S^{[*]}$ one has
\begin{equation}
\label{prod00++}
  [g,h]-[f',u]=[g,h]-[f,u'].
\end{equation}
Clearly, $\{h,u'\} \in SX$ and thus \eqref{prod00+} implies that
$[g,h]=[f',u]$ for all $\{h,u\} \in X$. This means that $\{f',g\}
\in X^{[*]}$. Thus $(SX)^{[*]} \subset X^{[*]}S^{[*]} $.

(ii) This statement is obtained by applying part (i) to the inverse
$(YS)^{-1}=S^{-1}Y^{-1}$.
\end{proof}

The following theorem concerns the composition of two unitary
relations; the results therein will be important in the sequel.

\begin{theorem}
\label{product} Let $\sK_1$, $\sK_2$, and $\sK_3$ be Kre\u{\i}n
spaces and let the linear relations $T:\sK_1\to\sK_2$ and
$S:\sK_2\to\sK_3$ be isometric. Then:
\begin{enumerate}
\def\labelenumi{\rm (\roman{enumi})}

\item the linear relation $S T:\sK_1\to\sK_3$ is isometric.
\end{enumerate}
In addition, let the linear relations $T:\sK_1\to\sK_2$ and
$S:\sK_2\to\sK_3$ be unitary. Then:

\begin{enumerate}
\def\labelenumi{\rm (\roman{enumi})}
\setcounter{enumi}{1}

\item
if $\dom S$ and $T \hplus (\{0\}\times\ker S)$ are closed then $S T:
\sK_1\to\sK_3$ is unitary;

\item if
\begin{equation}\label{UnRel}
   \ran T\subset\dom S \quad\text{and }\quad \dom S \,\text{ is
   closed},
\end{equation}
then $S T: \sK_1\to\sK_3$ is unitary and $\dom ST=\dom T$;

\item
if $\dom T$ and $S \hplus (\mul T\times\{0\})$ are closed then $S T:
\sK_1\to\sK_3$ is unitary;

\item if
\begin{equation}\label{UnRel2}
   \ran T\supset\dom S \quad\text{and }\quad \dom T \,\text{ is
   closed},
\end{equation}
then $S T: \sK_1\to\sK_3$ is unitary and $\ran ST=\ran S$;

\item if $\ran T=\dom S$ and $\ran S=\sK_3$, then the unitary
relation $S T:\sK_1\to\sK_3$ is bounded and single-valued (not
necessarily densely defined);

\item if $T\in[\sK_1,\sK_2]$ or $S\in[\sK_2,\sK_3]$, then $S
T:\sK_1\to\sK_3$ is unitary;

\item if $T\in[\sK_1,\sK_2]$ and $S\in[\sK_2,\sK_3]$, then $ST$ is
a unitary operator which belongs to $[\sK_1,\sK_3]$.

\end{enumerate}
\end{theorem}

\begin{proof}

(i) Since $ S$ and $ T$ are isometric, one has $S^{-1}\subset
S^{[*]}$ and $T^{-1}\subset T^{[*]}$. The definition of the product
of relations implies that $T^{-1} S^{-1}\subset T^{[*]} S^{[*]}$.
Lemma~\ref{basicle} yields
\begin{equation}
\label{prod00} (S T)^{-1}=T^{-1} S^{-1}\subset T^{[*]} S^{[*]}
\subset (S T)^{[*]}.
\end{equation}
Hence, the relation $S T$ is isometric.

(ii) Since $ S$ and $ T$ are unitary, $ST$ is isometric by part (i),
i.e., $(S T)^{-1} \subset (S T)^{[*]}$. To see that $ S T$ is
unitary it suffices to prove the inclusion $( S T)^{[*]} \subset (S
T)^{-1}= T^{[*]} S^{[*]}$ (where the last identity is due to $S$ and
$T$ being unitary). The linear relation $T_0$ defined by
\[
 T_0:=T\cap(\sH_1\times\dom S)=\{\,\{h,h'\}\in T:\, h'\in\dom S \,\}
\]
satisfies the inclusion $\ran T_0\subset\dom S$. Hence from
Lemma~\ref{prodlemma} one obtains
\[
 (ST)^{[*]}\subset (ST_0)^{[*]}=T_0^{[*]}S^{[*]}.
\]
Now it is enough to prove that $T_0^{[*]}S^{[*]}\subset
T^{[*]}S^{[*]}$ (then also $T_0^{[*]}S^{[*]}=T^{[*]}S^{[*]}$ holds).
Since $T$ is unitary, it follows from the assumptions in (ii) that
\begin{equation}\label{T0T}
 T_0^{[*]}=T^{[*]} \hplus (\ker S\times \{0\}).
\end{equation}
Now let $\{f,g\} \in T_0^{[*]}S^{[*]}$. Then for some $f' \in \sK_2$
one has $\{f,f'\} \in S^{[*]}$ and $\{f',g\} \in T_0^{[*]}$. Hence
due to~\eqref{T0T} $\{f'-f_0,g\} \in T^{[*]}$ for some $f_0\in\ker
S$. Since $S$ is unitary one has $f_0\in\mul S^{[*]}\,(=\ker S)$.
Thus $\{f,f'-f_0\} \in S^{[*]}$ and therefore $\{f,g\} \in
T^{[*]}S^{[*]}$. This completes the proof of part (ii).

(iii) By the assumptions in \eqref{UnRel} one obtains the statement
directly from Lemma~\ref{prodlemma}:
$(ST)^{[*]}=T^{[*]}S^{[*]}=T^{-1}S^{-1}=(ST)^{-1}$. The equality
$\dom ST = \dom T$ is clear due to the assumption $\ran T\subset\dom
S$.

(iv) This statement is obtained by applying part (ii) to the inverse
$(ST)^{-1}=T^{-1}S^{-1}$ and by taking into account
Lemma~\ref{UNITlem} and the equivalence stated in (i) of
Proposition~\ref{UNIT}.

(v) This is again an immediate consequence of Lemma~\ref{prodlemma};
it can be obtained also from (iii) by means of inverses.

(vi) If $\ran T=\dom S$ and $\ran S=\sK_3$, then $\dom S$ and $\dom
T$ are closed by (i) of Proposition~\ref{UNIT}. Therefore, by part
(v), the relation $ST:\sK_1\to\sK_3$ is unitary and $\ran ST =\ran
S=\sK_3$. Furthermore, $ST$ bounded and single-valued by (iii) of
Lemma~\ref{UNITcor}.

(vii) The relations $S$ and $T$ are assumed to be unitary. If in
addition $S\in[\sK_2,\sK_3]$, then by definition $\dom S=\sK_2$ and
moreover $\ker S=\{0\}$ by Proposition~\ref{UNIT}. Hence the
relation $ST$ is unitary by part (ii). On the other hand, if $T
\in[\sK_1,\sK_2]$, then $\dom T=\sK_1$, $\ran T=\sK_2$, and now part
(v) shows that $ST$ is unitary.

(viii) This is clear and a well-known fact.
        \end{proof}

Observe that in Theorem~\ref{product} the only standard result in
the literature is the last statement (viii). Notice also that (iii)
is in fact a special case of (ii). Indeed, if $\ran T\subset\dom S$
then $\mul T \supset \ker S$ by Proposition~\ref{UNIT} and hence in
this case $T \hplus (\{0\}\times\ker S)=T$ is closed. Likewise (v)
is a special case of (iv).

\begin{corollary}
\label{prodcor} Let the linear relations $T:\sK_1\to\sK_2$ and
$S:\sK_2\to\sK_3$ be unitary. Then:
\begin{enumerate}
\def\labelenumi{\rm (\roman{enumi})}

\item
if $\dom S$ is closed and $\dim\ker S<\infty$ then $S T:
\sK_1\to\sK_3$ is unitary;

\item
if $\dom T$ is closed and $\dim\mul T<\infty$ then $S T:
\sK_1\to\sK_3$ is unitary;

\item
if $\dim \sK_2<\infty$ then $S T: \sK_1\to\sK_3$ is unitary.

\end{enumerate}
\end{corollary}
\begin{proof}
The statements (i) and (ii) are immediate consequences of parts (ii)
and (iv) in Theorem~\ref{product}, respectively.

As to (iii) observe that if $\sK_2$ is finite-dimensional then
automatically the assumptions in (i) and (ii) are satisfied.
\end{proof}

The following examples show that in the case of infinite dimensional
spaces unitary operators may be unbounded and their set does not
form a semigroup, that is, the product of two unitary operators need
not be a unitary operator.

\begin{example}\label{Unbound}
Let $K$ be a densely defined operator on a Hilbert space $\sH$ and
define the block operator matrix $T$ by
\begin{equation}\label{WUnit}
 T=\begin{pmatrix} I_{\sH} & K \\ 0 & I_{\sH}
 \end{pmatrix}.
\end{equation}
Then $T$ is an injective operator, i.e., $\ker T=\{0\}$, $\mul
T=\{0\}$. It is easy to see that $T$ is closed if and only if $K$ is
closed. The inverse of $T$ is given by
\begin{equation}\label{WUnit2}
 T^{-1}=\begin{pmatrix} I_{\sH} & -K \\ 0 & I_{\sH}
 \end{pmatrix}
\end{equation}
and hence $T$ is densely defined with dense range; in fact $\dom
T=\ran T=\sH\oplus\dom K$. Now consider $\sH\oplus\sH$ as the
Kre\u{\i}n space $(\sH^2,J_\sH)$ with the fundamental symmetry
\begin{equation}\label{jh}
 J_\sH:=\begin{pmatrix} 0 & -iI_\sH \\ i I_\sH & 0\end{pmatrix}.
\end{equation}
Then
\begin{equation}\label{WUnit3}
 T^{[*]}=\begin{pmatrix} I_{\sH} & -K^{*} \\ 0 & I_{\sH}
 \end{pmatrix}.
\end{equation}
The identities \eqref{WUnit2} and \eqref{WUnit3} show that $T$ is
isometric (unitary) if and only if $K$ is symmetric (resp.
selfadjoint). Therefore, if $K_1$, $K_2$ are two unbounded
selfadjoint operators in $\sH$ such that $K_1+K_2$ is not
selfadjoint, the product $T_1T_2$ of the unitary operators $T_1$ and
$T_2$,
\[
 T_1T_2=\begin{pmatrix} I_{\sH} & K_1 \\ 0 & I_{\sH}\end{pmatrix}
 \begin{pmatrix} I_{\sH} & K_2 \\ 0 & I_{\sH} \end{pmatrix}
 =\begin{pmatrix} I_{\sH} & K_1+K_2 \\ 0 & I_{\sH}\end{pmatrix},
\]
is not a unitary operator in $(\sH^2,J_\sH)$. Here both assumptions
in \eqref{UnRel} can fail to hold. This is the case if, for
instance, $K_1$ and $K_2$ are selfadjoint operators in $\sH$ such
that $\dom K_1\cap \dom K_2=\{0\}$.

Note also that if $K_1$ is an unbounded selfadjoint operator in
$\sH$ and $K_2=-K_1$ then $\ran T_2=\dom T_1$, cf. \eqref{WUnit2},
so that $\dom T_1T_2=\dom T_2$. Now the product $T_1T_2$ is not
closed and hence it cannot be unitary. In this case the first
assumption in \eqref{UnRel} and \eqref{UnRel2} is satisfied, while
the second assumption in \eqref{UnRel}, \eqref{UnRel2} fails to
hold. The second assumption in (ii) and (iv) of
Theorem~\ref{product} is also satisfied, since $\ker T_1=\{0\}$ and
$\mul T_2=\{0\}$.
\end{example}

Obviously, $ST$ can be unitary even if the assumptions \eqref{UnRel}
and \eqref{UnRel2} are not satisfied. Also only one of the two
conditions in \eqref{UnRel}, \eqref{UnRel2} or in (ii), (iv) of
Theorem~\ref{product} is not sufficient for the product $ST$ to be
unitary.

\begin{example}
Let  $K_2$ be a selfadjoint operator in $\sH$. Then the linear
relation $T_2$ given by
\[
T_2=\left\{\left\{\left(%
\begin{array}{c}
  K_2h \\
  h \\
\end{array}\right),%
\left(%
\begin{array}{c}
  0 \\
  g \\
\end{array}%
\right)\right\}:\, h\in\dom K_2,\,\,g \in\sH\right\}
\]
is unitary in $(\sH^2,J_\sH)$ with $\dom T_2=\ker
T_2=(\mbox{gr\,}K_2)^{-1}$, $\ran T_2=\mul T_2=\{0\}\times
\sH(\subset\sH\times \sH)$ closed. If $T_1$ is as in \eqref{WUnit}
with $K_1$ a selfadjoint operator then the product $T_1T_2$ is
unitary. Here $\dom T_1$ is closed if and only if $K_1$ is bounded,
in which case the assumptions in \eqref{UnRel} are satisfied.
However, if $K_1$ is unbounded then both of the assumptions in
\eqref{UnRel} fail to hold and also the first assumption in
\eqref{UnRel2} is not satisfied. It is not difficult to check that
both assumptions in (iv) of Theorem~\ref{product} are satisfied.

The product $T_2T_1$ is given by
\[
 T_2T_1=\left\{\left\{\left(%
\begin{array}{c}
  (K_2-K_1)h \\
  h \\
\end{array}\right),%
\left(%
\begin{array}{c}
  0 \\
  g \\
\end{array}%
\right)\right\}:\, h\in\dom K_2\cap\dom K_1,\,\,g \in\sH\right\}.
\]
This relation is unitary if and only if $K_2-K_1$ is selfadjoint.
Now the second assumption in \eqref{UnRel} is satisfied, while the
first assumption in \eqref{UnRel} does not hold. If, for instance,
$\dom K_1\cap\dom K_2=\{0\}$, then $T_2T_1$ is not unitary. In this
case both assumptions in \eqref{UnRel2} fail to hold. On the other
hand, if $K_1$ is bounded then the assumptions in \eqref{UnRel2} are
satisfied and $T_2T_1$ is unitary. The product $T_2T_1$ is also
unitary if $K_2$ is bounded, while both of the assumptions in
\eqref{UnRel2} fail to hold if $K_1$ is unbounded.

The first assumption in (ii) of Theorem~\ref{product} holds. The
second assumption in (ii) of Theorem~\ref{product} is equivalent for
the row operator $(K_1\, K_2)$ to be closed, which therefore by part
(ii) implies that $K_2-K_1$ is selfadjoint. Obviously, $K_2-K_1$ can
be selfadjoint even if the row operator $(K_1\, K_2)$ is not closed:
consider e.g. $-K_1=K_2=:K$: here $K_2-K_1=2K$ is selfadjoint, but
the row operator $(K_1\, K_2)=(-K\, K)$ is not closed if $K$ is
unbounded: let $h_n\to h \not\in \dom K$. Then $\col (h_n,h_n)\to
\col (h,h)$ and $(-K\, K)\col (h_n,h_n)\equiv 0$, and closedness
would imply that $h\in\dom K$.

Also observe that the linear relation $T_2^2=T_2$ is unitary in
$(\sH^2,J_\sH)$. However, $\dom T_2\cap\ran T_2=\{0\}$, if
$0\not\in\sigma_p(K_2)$ and in this case $\dom T_2$ and $\ran T_2$
are closed.
\end{example}


\subsection{The main transform}

It is convenient to interpret the Hilbert space $\sH^2=\sH\oplus\sH$
as a Kre\u{\i}n space $(\sH^2,J_\sH)$ whose inner product is
determined by the fundamental symmetry $J_\sH$ of the
form~\eqref{jh};
notice the connection to the definition of the adjoint of linear
relations in \eqref{adjo}. There is a useful and important transform
which gives a connection between the subspaces of a Hilbert space
$\sH\oplus\cH$ and linear relations from the Kre\u{\i}n space
$(\sH^2,J_\sH)$ to the Kre\u{\i}n space $(\cH^2,J_\cH)$, which will
be now recalled from \cite{Trans}. Let $\sH$ and $\cH$ be Hilbert
spaces and let their Cartesian product be denoted by $\wt \sH=\sH
\oplus \cH$. Define the linear mapping $\cJ$ from $\sH^2\times\cH^2$
to $(\sH \oplus \cH)^2$ by
\[
 \cJ:\left\{ \begin{pmatrix} f \\ f' \end{pmatrix},
         \begin{pmatrix} h \\ h' \end{pmatrix}
     \right\}
 \mapsto
     \left\{ \begin{pmatrix} f \\ h \end{pmatrix},
         \begin{pmatrix} f' \\ -h' \end{pmatrix}
     \right\},
\quad f,f'\in\sH,\,\, h,h'\in\cH.
\]
This mapping establishes a one-to-one correspondence between the (closed) linear
relations $\Gamma:\sH^2 \to \cH^2$ and the (closed) linear
relations $\wt A$ in $\wt \sH=\sH\oplus\cH$ via
\begin{equation}
\label{awig}
 \Gamma\mapsto\wt A :=\cJ(\Gamma)=
 \left\{\,
 \left\{ \begin{pmatrix} f \\ h \end{pmatrix},
         \begin{pmatrix} f' \\ -h' \end{pmatrix}
 \right\} :\,
 \left\{ \begin{pmatrix} f \\ f' \end{pmatrix},
         \begin{pmatrix} h \\ h' \end{pmatrix}
 \right\}
 \in \Gamma
 \,\right\}.
\end{equation}
The mapping $\cJ$ plays a principal role  and it is refered to as
the \textit{main transform}. Some basic properties of this transform
are stated in the following proposition.

\begin{proposition}
\label{uksi} Let the linear relation $\Gamma$ from $(\sH^2,J_\sH)$
to $(\cH^2,J_\cH)$ and the linear relation $\wt A$ in $\sH\oplus\cH$
be connected by $\wt A=\cJ(\Gamma)$. The main transform $\cJ$
establishes a one-to-one correspondence between the contractive,
isometric, and unitary  relations $\Gamma$ from $(\sH^2,J_\sH)$ to
$(\cH^2,J_\cH)$ and the dissipative, symmetric, and  selfadjoint
relations $\wt A$ in $\sH \oplus \cH$, respectively.
\end{proposition}


\subsection{Nevanlinna families}

A family of linear relations $M(\lambda)$, $\lambda \in \cmr$, in
a Hilbert space $\cH$ is called a \textit{Nevanlinna family} if:
\begin{enumerate}
\def\labelenumi{\rm (\roman{enumi})}
\item for every $\lambda \in \dC_+ (\dC_-)$
      the relation $M(\lambda)$ is maximal dissipative (resp. accumulative);
\item $M(\lambda)^*=M(\bar \lambda)$, $\lambda \in \cmr$;
\item for some, and hence for all, $\mu \in \dC_+ (\dC_-)$ the
      operator family $(M(\lambda)+\mu)^{-1} (\in [\cH])$ is
      holomorphic for all $\lambda \in \dC_+ (\dC_-)$.
\end{enumerate}
By the maximality condition, each relation $M(\lambda)$, $\lambda
\in \cmr$, is necessarily closed. The \textit{class of all
Nevanlinna families} in a Hilbert space is denoted by $\wt
R(\cH)$.
If  the multi-valued part $\mul
M(\lambda)$ of $M(\cdot) \in \wt R(\cH)$ is nontrivial, then it is independent of $\lambda \in \cmr$, so that
\begin{equation}
\label{ml}
 M(\lambda)=M_s(\lambda) \oplus M_\infty, \quad M_\infty=\{0\}
 \times \mul M(\lambda), \quad \lambda\in\dC\setminus\dR,
\end{equation}
where $M_s(\lambda)$ is a Nevanlinna family of densely defined
operators in $\cH \ominus \mul M(\lambda)$, \cite{KL1}. 

Clearly, if $M(\cdot) \in \wt R(\cH)$, then
$M_\infty \subset M(\lambda) \cap M(\lambda)^*$ for all $\lambda\in\dC\setminus\dR$.
The following \textit{subclasses} of the class $\wt R(\cH)$ will be useful:
\begin{enumerate}
\item[] $R(\cH)$ is the set of all $M(\cdot) \in \wt R(\cH)$
         for which $\mul M(\lambda)=\{0\}$;
\item[] $R^s(\cH)$ is the set of all $M(\cdot) \in \wt R(\cH)$
         for which $M(\lambda) \cap M(\lambda)^*=\{0\}$ for all $\lambda\in\cmr$;
\item[] $R^u(\cH)$ is the set of all $M(\cdot) \in \wt R(\cH)$
         for which $M(\lambda) \hplus M(\lambda)^*=\cH^2$ for all $\lambda\in\cmr$;
\item[] $R[\cH]$ is the set of all $M(\cdot) \in \wt R[\cH]$
         for which $\dom M(\lambda)=\cH$ for all $\lambda\in\cmr$;
\item[] $R^s[\cH]$ is the set of all $M(\cdot) \in R[\cH]$
         for which $\ker \IM M(\lambda)=\{0\}$ for all $\lambda\in\cmr$;
\item[] $R^u[\cH]$ is the set of all $M(\cdot) \in R^s[\cH]$
         for which $0 \in \rho(\IM M(\lambda))$ for all $\lambda\in\cmr$;
\item[] $\wt R^c(\cH)$ is the set of all constant Nevanlinna families.

\end{enumerate}
The  subclasses of $\wt R(\cH)$ can be equivalently defined by
assuming the corresponding property of $M(\lambda)$ only at a single
point $\lambda\in\cmr$,~\cite{Trans}. Moreover,  it is easy to show
that $R^u[\cH]=R^u(\cH)$, see~\cite{Trans}.   The Nevanlinna
functions in $R^s[\cH]$ and $R^u[\cH]$ will be called
\textit{strict} and \textit{uniformly strict}, respectively.

If $M(\cdot) \in R[\cH]$,
then it admits the following integral representation
\begin{equation}
\label{INTrep}
 M(\lambda)
  =A+B\lambda+\int_{\dR}\left(\frac{1}{t-\lambda}-\frac{t}{t^2+1}\right)\,
               d\Sigma(t),
 \quad
 \int_{\dR}\,\frac{d\Sigma(t)}{t^2+1}\in[\cH],
\end{equation}
where $A=A^*\in[\cH]$, $0\le B=B^*\in[\cH]$, the $[\cH]$-valued
family $\Sigma(\cdot)$ is nondecreasing, and the integral is
uniformly convergent in the strong topology, cf. \cite{Br},
\cite{KacK}.

A pair $\{\Phi,\Psi\}$ of holomorphic $[{\cH}]$-valued functions on
$\dC_+\cup\dC_-$ is said to be a {\it Nevanlinna pair} if:
\begin{enumerate}
\def\labelenumi{\rm (\roman{enumi})}
\item[(N1)]
$\IM \Phi(\lambda)^*\Psi(\lambda)/\IM \lambda \geq 0$, $\lambda\in
\dC_+\cup\dC_-$;
\item[(N2)]
$\Psi(\bar \lambda)^*\Phi(\lambda)-\Phi(\bar
\lambda)^*\Psi(\lambda)=0$, $\lambda \in {\dC_+\cup\dC_-}$;
\item[(N3)]
$0\in\rho(\Psi(\lambda)\pm i\Phi(\lambda))$, $\lambda \in\dC_\pm$.
\end{enumerate}
Two Nevanlinna pairs $\{\Phi_1,\Psi_1\}$ and $\{\Phi_2,\Psi_2\}$ are
said to be equivalent, if $\Phi_2
(\lambda)=\Phi_1(\lambda)\chi(\lambda)$ and
$\Psi_2(\lambda)=\Psi_1(\lambda)\chi(\lambda)$ for some operator
function $\chi(\lambda)\in[\cH]$, which is holomorphic and
invertible on $\dC_+\cup\dC_-$. If $\{\Phi,\Psi\}$ is a Nevanlinna
pair, then the following kernel is nonnegative on $\dC_+\cup\dC_-$:
\begin{equation}
\label{Nkern}
 {\sf N}_{\Phi \Psi} (\lambda,\mu)
  =\frac{\Phi(\mu)^* \Psi(\lambda)-
   \Psi(\mu)^* \Phi(\lambda)}{\lambda-\bar \mu},
\quad
 \lambda,\mu \in \dC_+\cup\dC_-.
\end{equation}

The set of Nevanlinna families $\tau(\lambda)$ and the set of
equivalence classes of Nevanlinna pairs $\{\Phi,\Psi\}$ are in a
 one-to-one correspondence via the formula
\begin{equation}
\label{tau1}
 \tau(\lambda)=\{\Phi(\lambda),\Psi(\lambda)\}:=
 \{\,\{\Phi(\lambda)h,\Psi(\lambda)h\}:\, h\in\cH\,\}.
\end{equation}
Moreover, {strict} and {uniformly strict} Nevanlinna families are
charachterized by the conditions $0\not\in\sigma_p({\sf N}_{\Phi
\Psi} (\lambda,\lambda))$ and $0\in\rho({\sf N}_{\Phi \Psi}
(\lambda,\lambda))$ for some $\lambda\in \dC\setminus\dR$,
respectively.

\subsection{Shmul'yan transform of linear relations}

Let $\cH$ and $\cK$ be Hilbert spaces and let $W$ be a linear
relation from the Hilbert space $\sH^2=\cH \oplus \cH$ to the
Hilbert space $\cK^2=\cK \oplus \cK$. For any linear relation
$\Theta$ in $\cH$,
\begin{equation}
\label{Lft0}
 W[\Theta]=\{\, \wh{k} \in \cK^2 :\, \{\wh h, \wh k\} \in W, \,\, \wh h \in
 \Theta\,\},
\end{equation}
defines a linear relation $W[\Theta]$ in $\cK$.

\begin{definition}
The linear relation $W[\Theta]$ in $\cK$, defined by \eqref{Lft0},
is said to be the {\it Shmul'yan transform} of the linear relation
$\Theta$ in $\cH$, induced by the linear relation $W: \cH^2 \to
\cK^2$.
\end{definition}

For any pair of relations $\Theta_1$ and $\Theta_2$ in the Hilbert
space $\cH$, 
there are the
inclusions
\[
W[\Theta_1 \cap \Theta_2] \subset W[\Theta_1] \cap W[\Theta_2],
\]
and
\begin{equation}
\label{Lft2} W[\Theta_1] \hplus W[\Theta_2] \subset W[\Theta_1
\hplus \Theta_2].
\end{equation}
Furthermore, if $\ker W=\{0\}$, then
\begin{equation}
\label{lft3} W[\Theta_1 \cap \Theta_2]=W[\Theta_1] \cap W[\Theta_2],
\end{equation}
and, if $\mul W=\{0\}$, then
\begin{equation}W[\Theta_1] \hplus W[\Theta_2] = W[\Theta_1
\hplus \Theta_2].
\end{equation}

Now interpret $W$ as a linear relation from the Kre\u{\i}n space
$(\cH^2,J_\cH)$ to the Kre\u{\i}n space $(\cK^2,J_\cK)$, where the
inner products are defined as in \eqref{jh}.

\begin{lemma}\label{smul}
Let $W$ be a linear relation from the Kre\u{\i}n space
$(\cH^2,J_\cH)$ to the Kre\u{\i}n space $(\cK^2,J_\cK)$. Then $\{\wh
\alpha, \wh \beta\} \in W^{[*]}$ if and only if
\begin{equation}\label{ADJ}
 [\wh \beta, \wh h]_{\cH^2}=[\wh \alpha, \wh k]_{\cK^2}
\quad \mbox{for all} \quad \{\wh h, \wh k\} \in W.
\end{equation}
Let $\Theta$ be a linear relation in $\cH$  and let  $\{\wh \alpha,
\wh \beta\} \in W^{[*]}$.  Then
\begin{equation}\label{ADJ+}
 \wh \beta \in \Theta^* \Leftrightarrow \wh \alpha \in W[\Theta]^*.
\end{equation}
\end{lemma}

\begin{proof}
Since $\{\wh \alpha, \wh \beta\} \in W^{[*]}$, the identity
\eqref{ADJ} is satisfied for all $\{\wh h, \wh k\} \in W$. The
equivalence in \eqref{ADJ+} is a straightforward consequence of
\eqref{ADJ} and the connection of the inner products in the
Kre\u{\i}n spaces $(\cH^2,J_\cH)$ and $(\cK^2,J_\cK)$ to
the definition of adjoint in \eqref{adjo}.
\end{proof}

\begin{corollary}\label{smul1}
Let $W$ be an isometric linear relation from the Kre\u{\i}n space
$(\cH^2,J_\cH)$ to the Kre\u{\i}n space $(\cK^2,J_\cK)$ and let
$\Theta$ be a linear relation in $\cH$. Then
\[
W[\Theta^*] \subset W[\Theta]^*.
\]
Moreover, $\Theta $ is dissipative (symmetric) if and only if
$W[\Theta]$ is  dissipative (symmetric).
\end{corollary}

\begin{proof}
Let $\wh k \in W[\Theta]$ and $\wh \alpha \in W[\Theta^*]$. Then
there exist elements $\wh h \in \Theta$ and $\wh \beta \in
\Theta^*$ such that $\{\wh h, \wh k\} \in W$ and
$\{\wh \beta, \wh \alpha\} \in W$.
Since $\{\wh \alpha,\wh \beta \} \in
W^{-1}=W^{[*]}$, it follows from Lemma~\ref{smul} that
\[
 [\wh k, \wh \alpha]_{\cK^2}=[\wh h, \wh \beta]_{\cH^2}=0.
\]
This shows that $W[\Theta]$ and $W[\Theta^*]$ are orthogonal in the
Kre\u{\i}n space $(\cK^2,J_\cK)$.

For every $\{\wh h, \wh k\} \in W$ one has $\{\wh k, \wh h\} \in
W^{-1}\subset W^{[*]}$, since $W$ is isometric. Therefore,
\begin{equation}\label{equaa}
 0=[\wh h, \wh h]_{\cH^2}-[\wh k, \wh k]_{\cK^2}=2i\left[ \IM
 (h',h)-\IM (k',k) \right], \quad \{\wh h, \wh k\} \in W,
\end{equation}
where the identity on the left is due to \eqref{ADJ}, and the
identity on the right is due to the definition of the inner
products. Note that if $\wh h \in \Theta$, then there exists $\wh k
\in W[\Theta]$ with $\{\wh h, \wh k\} \in W$. Hence, if $W[\Theta]$
is dissipative or symmetric, then \eqref{equaa} shows that $\Theta$
is dissipative or symmetric, respectively. Conversely, if $\wh k \in
W[\Theta]$, then there exists $\wh h \in \Theta$ with $\{\wh h, \wh
k\} \in W$. Hence, if $ \Theta $ is dissipative or symmetric, then
so is $W[\Theta]$.
\end{proof}

In the general context of Corollary \ref{smul1} it seems difficult
to conclude anything about the maximality of the dissipative
(symmetric) relations. Of course, when $W$ is a standard unitary
operator, then some known properties can be easily recovered,
cf.~\cite{KSh},\cite{Sh}.

\begin{corollary}
\label{LFTprop} Let $W$ be a standard unitary operator from the
Kre\u{\i}n space $(\cH^2,J_\cH)$ onto the Kre\u{\i}n space
$(\cK^2,J_\cK)$. Let $\Theta$ be a linear relation in $\cH$.  Then:
\[
W[\Theta^*] = W[\Theta]^*.
\]
Furthermore,
\begin{enumerate}
\def\labelenumi{\rm (\roman{enumi})}

\item  $\Theta $ is   maximal  dissipative
$\Leftrightarrow$ $W[\Theta]$  is  maximal   dissipative;

\item $\Theta $ is  maximal  symmetric
$\Leftrightarrow$  $W[\Theta ]$  is  maximal  symmetric;

\item $\Theta $  is selfadjoint
$\Leftrightarrow$  $W[\Theta ]$  is selfadjoint.
\end{enumerate}
\end{corollary}


In case $W$ is a standard unitary operator from $\cH^2$ onto
$\cK^2$, the Shmul'yan transform is usually written out in
components. Then $W$ is bounded and with bounded inverse and it can
be represented in the block form
\begin{equation}
\label{W00} W=\begin{pmatrix} W_{00} & W_{01} \\ W_{10} & W_{11}
  \end{pmatrix},\quad W_{ij}\in[\cH,\cK], \quad i,j=0,1.
\end{equation}
 If $\Theta$ is a linear relation in $\cH$,   then $W[\Theta]$ in
\eqref{Lft0} takes the form
\begin{equation}
\label{lft0}
 W[\Theta]=\{\,\{W_{00}h+W_{01}h',W_{10}h+W_{11}h' \}:\, \{h,h'\} \in
\Theta\,\}.
\end{equation}
Clearly, $W[\Theta]$ is contained in the linear relation
\begin{equation}
\label{lft} (W_{10}+W_{11}\Theta)(W_{00}+W_{01}\Theta)^{-1}
=\{W_{00}h+W_{01}h',W_{10}h+W_{11}h''\}:\, \{h,h'\},\{h,h''\}\in
\Theta\,\}.
\end{equation}
In fact, the following equality holds
\begin{equation}
\label{lft1}
  (W_{10}+W_{11}\Theta)(W_{00}+W_{01}\Theta)^{-1}
   =W[\Theta] \hplus \{0,W_{11}(\mul\Theta)\}.
\end{equation}
Hence, if $\Theta$ is a relation with $W_{11}(\mul\Theta)=\{0\}$,
and in particular if $\Theta$ is an operator, the linear relations
in \eqref{lft0} and in \eqref{lft} coincide.

\section{Boundary relations and Weyl families}
\label{BR}

\subsection{Definitions and basic properties}

Let $S$ be a closed symmetric linear relation in the Hilbert space
$\sH$. It is not assumed that the defect numbers of $S$ are equal or
finite. A boundary relation for $S^*$ is defined as follows.

\begin{definition}
\label{GBT} Let $S$ be a closed symmetric linear relation in a
Hilbert space $\sH$ and let $\cH$ be an auxiliary Hilbert space. A
linear relation $\Gamma:\sH^2\mapsto\cH^2$ is called a
\textit{boundary relation} for $S^*$, if:
\begin{enumerate}
\def\labelenumi{\rm (\roman{enumi})}
\item[(G1)] $\dom\Gamma$ is dense in $S^*$ and the identity
\begin{equation}
\label{Green1}
 (f',g)_\sH-(f,g')_\sH=({h'},k)_{\cH}-(h,k')_{\cH},
\end{equation}
holds for every $\{\wh f,\wh h\},\, \{\wh g,\wh k\}\in\Gamma$;
\item[(G2)] $\Gamma$ is maximal in the sense that if
$\{\wh g,\wh k\}\in\sH^2\times\cH^2$
satisfies \eqref{Green1} for every $\{\wh f,\wh h\}\in \Gamma$,
then $\{\wh g,\wh k\}\in\Gamma$.
\end{enumerate}
Here $\wh f=\{f,f'\}$, $\wh g=\{g,g'\}\in\dom\Gamma(\subset\sH^2)$,
$\wh h=\{h,h'\}$, $\wh k=\{k,k'\}\in\ran\Gamma(\subset\cH^2)$).
\end{definition}

The condition \eqref{Green1} in (G1) can be interpreted as an
abstract Green's identity. Using the terminology of Kre\u{\i}n
spaces \eqref{Green1} means that $\Gamma$ is an isometric relation
from  the Kre\u{\i}n space $(\sH^2,J_\sH)$ to the Kre\u{\i}n space
$(\cH^2,J_\cH)$, since
\begin{equation}
\label{Green}
 (J_{\sH}\wh f,\wh g)_{\sH^2}=
 (J_{\cH}\wh h,\wh k)_{\cH^2},
 \quad
 \{\wh f,\wh h\},
 \quad
 \{\wh g,\wh k\}\in\Gamma.
\end{equation}
The maximality condition (G2) and Proposition \ref{UNIT} now yield
the following result.

\begin{proposition} (\cite{Trans})
\label{prop1} Let $\sH$ and $\cH$ be Hilbert spaces and let $S$ be a
closed symmetric linear relation in $\sH$. Then a linear relation
$\Gamma:\sH^2\mapsto\cH^2$ is a boundary relation for $S^*$ if and
only if $\Gamma$ is a unitary relation from the Kre\u{\i}n space
$(\sH^2,J_\sH)$ to the Kre\u{\i}n space $(\cH^2,J_\cH)$ with $S=\ker
\Gamma$.
\end{proposition}
In some cases the following criterion for a linear relation to be a
boundary relation is useful; see \cite[Proposition~3.6]{Trans}.

\begin{proposition}
\label{thm1} The linear relation $\Gamma:\sH^2\mapsto\cH^2$ is a
boundary relation for $S^*$ if and only if the following conditions
hold:
\begin{enumerate}
\def\labelenumi{\rm (\roman{enumi})}
\item $\dom \Gamma$ is dense in $S^*$;
\item $\Gamma$ is closed and isometric
      from the Kre\u{\i}n space $(\sH^2,J_\sH)$
      to the Kre\u{\i}n space $(\cH^2,J_\cH)$;
\item $\ran (\Gamma(\wh \sN_\lambda(T))+\lambda)=\cH$
      for some (and, hence, for all) $\lambda\in\dC_+$
      and for some (and, hence, for all) $\lambda \in \dC_-$.
\end{enumerate}
\end{proposition}
Note that a boundary relation $\Gamma$ is automatically closed and
linear, since it is a unitary relation from the Kre\u{\i}n space
$(\sH^2,J_\sH)$ to the Kre\u{\i}n space $(\cH^2,J_\cH)$. Observe
that the inverse $\Gamma^{-1}:(\cH^2,J_\cH)\to (\sH^2,J_\sH)$ is
also unitary; see Lemma~\ref{UNITlem}. Therefore, in this case
$\Gamma^{-1}$ can be interpreted as a boundary relation for $\wt
S^*\subset \cH^2$ , the adjoint of the closed symmetric relation
\begin{equation}
\label{Stilde}
 \wt S:=\ker\Gamma^{-1}=\mul\Gamma \,(\subset \cH^2).
\end{equation}

Let $\Gamma$ be a boundary relation for $S^*$ and let
$T=\dom\Gamma$. According to~\cite[Proposition~2.12]{Trans}
 the linear relation $T$ in $\sH$ satisfies
\begin{equation}
\label{STS}
 S\subset T\subset S^*,
 \quad
 \clos T=S^*.
\end{equation}
The eigenspaces ${\sN}_{\lambda}(T)$ and $\wh{\sN}_{\lambda}(T)$ for
$T$ are defined by
\begin{equation}
\label{defectT}
 \sN_\lambda(T)=\ker (T-\lambda),
\quad
 \wh \sN_\lambda(T)
 =\{\,  \{f,\lambda f\} \in T :\, f\in\sN_\lambda(T)\,\}.
\end{equation}
For notational convenience the usual defect spaces of $S$ are
denoted here by $\sN_\lambda(S^*)$ and $\wh\sN_\lambda(S^*)$.

For all elements $\{\wh f_\lambda,\wh h\}, \{\wh g_\mu,\wh
k\}\in\Gamma$ with $\wh f_\lambda\in\wh\sN_\lambda(T)$ and $\wh
g_\mu\in\wh\sN_\mu(T)$ one has
\begin{equation}
\label{gr00}
 (\lambda-\bar{\mu})(f_\lambda,g_\mu)_\sH=(h',k)_\cH-(h,k')_\cH,
 \quad
 \lambda,\mu \in \cmr,
\end{equation}
which follows from the identity~\eqref{Green1}.
Hence, the subspace $\wh{\sN}_{\lambda}(T)$
is positive in the Kre\u{\i}n space
$(\sH^2,J_\sH)$ for $\lambda\in\dC_+$ and negative for $\lambda\in\dC_-$.

\begin{corollary}
\label{npmA} Let $\Gamma:\sH^2\mapsto\cH^2$  be a boundary relation
for a symmetric operator $A$. Then:
\begin{enumerate}
\def\labelenumi{\rm (\roman{enumi})}
\item $n_\pm(A)\le\dim \cH$;
\item
if $n_\pm(A)<\infty$, then $\dim \cH-n_\pm(A)=\dim\mul\Gamma$;
\item if $\dim\cH<\infty$, then $n_+(A) = n_-(A)$.
\end{enumerate}
\end{corollary}

\begin{proof}
(i) Let $\wt A=\cJ(\Gamma)$ be the main transform of $\Gamma$. It
follows from~\cite[Lemma 2.14]{Trans} that
\[
 n_\pm(A)=n_\pm(\wt S),
\]
where $\wt S=\mul \Gamma \subset \cH^2$; cf. \eqref{Stilde}. This
implies the statement (i).

(ii) If $n_\pm(A)<\infty$, or equivalently, $n_\pm(\wt S)<\infty$,
one obtains
\[
\dim\mul\Gamma=\dim \wt S=\dim\cH-\dim n_\pm(\wt S),
\]
where the last identity follows from the fact that $\dim \wt S$ in
$\cH^2$ is equal to $\dim \ran(\wt S-\lambda)$ in $\cH$ for all
$\lambda\in\cmr$.

(iii) If $\dim\cH<\infty$, then clearly $n_+(\wt S)=n_-(\wt S)$ and
therefore also $n_+(A)=n_-(A)$.
\end{proof}

Recall from \cite{Trans} that the boundary relation
$\Gamma:\sH^2\mapsto\cH^2$ is said to be \textit{minimal}, if
\[
 \sH=\sH_{min}:=\cspan\{\,\sN_\lambda(T):\,\lambda\in\dC_+\cup\dC_-\,\}.
\]
Since $\sN_\lambda(T)$ is dense in $\sN_\lambda(S^*)$,
the boundary relation $\Gamma:\sH^2\mapsto\cH^2$ is minimal if and
only if $S$ is simple.

\begin{definition}
\label{Weylfam} The \textit{Weyl family} $M(\cdot)$ of $S$
corresponding to the boundary relation $\Gamma:\sH^2\mapsto\cH^2$
is defined by $M(\lambda):=\Gamma(\wh{\sN}_\lambda(T))$, i.e.,
\begin{equation}
\label{Weylf}
 M(\lambda) :=
  \left\{\,\wh h\in\cH^2:\,\{\wh f_\lambda,\wh h\}\in\Gamma \mbox{
for some } \wh f_\lambda=\{f,\lambda f\}\in\sH^2\,\right\},
\end{equation}
where $\lambda \in \cmr$. In the case where $M(\cdot)$ is
operator-valued it is called the \textit{Weyl function} of $S$
corresponding to the boundary relation $\Gamma$.
\end{definition}

\begin{definition}
\label{gammafield} The \textit{$\gamma$-field} $\gamma(\cdot)$ of
$S$ corresponding to the boundary relation $\Gamma : \sH^2 \to
\cH^2$ is defined by
\begin{equation}
\label{gfield}
 \gamma(\lambda) :=
  \left\{\, \{h,f_\lambda\}\in \cH\times\sH :\,\{\wh f_\lambda,\wh h\}\in\Gamma \mbox{
for some } \wh f_\lambda=\{f,\lambda f\}\in\sH^2\,\right\},
\end{equation}
where $\lambda \in \cmr$. Moreover, $\wh\gamma(\lambda)$ stands for
\begin{equation}
\label{ghatfield}
 \wh\gamma(\lambda):
  =\left\{\, \{h,\wh f_\lambda\}\in \cH\times\sH^2:\,
  \{h,f_\lambda\}\in\gamma(\lambda),\quad \wh f_\lambda=\{f,\lambda f\}\in\sH^2\,\right\},
 \quad \lambda \in \cmr.
\end{equation}
\end{definition}

Associate with $\Gamma$ the following linear relations which are not
necessarily closed:
\begin{equation}
\label{g01}
\begin{split}
& \Gamma_0 =\left\{ \, \{\wh f,h\} :\, \{\wh f, \wh h\} \in \Gamma,
\,
           \wh h=\{h,h'\}\,\right\}, \\
& \Gamma_1 =\left\{ \, \{\wh f,h'\} :\, \{\wh f, \wh h\} \in \Gamma,
\,
             \wh h=\{h,h'\}\,\right\}.
\end{split}
\end{equation}
It is clear that
\[
 \dom M(\lambda)=\Gamma_0(\wh\sN_\lambda(T))\subset\ran\Gamma_0,
 \quad
 \ran M(\lambda)=\Gamma_1(\wh\sN_\lambda(T))\subset\ran\Gamma_1.
\]
If the boundary relation $\Gamma$ is single-valued the triplet
$\{\cH,\Gamma_0,\Gamma_1\}$ will be called a \textit{boundary
triplet} associated with the boundary relation
$\Gamma:\sH^2\mapsto\cH^2$. In this case the Weyl family
corresponding to the boundary triplet $\{\cH,\Gamma_0,\Gamma_1\}$
can be also defined via the equality
\begin{equation}\label{0.2++}
 \Gamma_1(\{f_{\lambda},\lambda f_{\lambda}\})
 =M(\lambda)\Gamma_0 (\{f_{\lambda},\lambda f_{\lambda}\}),
 \quad
 \{f_{\lambda},\lambda f_{\lambda}\} \in T.
\end{equation}

The $\gamma$-field $\gamma(\cdot)$ associated with the boundary
relation $\Gamma : \sH^2 \to \cH^2$ is the first component of the
mapping $\wh\gamma(\lambda)$ in \eqref{ghatfield}. Observe that
\[
 \wh\gamma(\lambda):=(\Gamma_0\uphar\wh\sN_\lambda(T))^{-1},
 \quad \lambda\in\cmr,
\]
is a linear mapping from $\Gamma_0(\wh\sN_\lambda(T))=\dom
M(\lambda)$ onto $\wh\sN_\lambda(T)$; it is single-valued in view of
\eqref{gr00}. Consequently, the $\gamma$-field is a single-valued
mapping from $\dom M(\lambda)$ onto $\sN_\lambda(T)$ and satisfies
$\gamma(\lambda)\Gamma_0\wh f_\lambda=f_\lambda$ for all $\wh
f_\lambda\in \wh\sN_\lambda(T)$.

\subsection{Realization theorem}

It follows from the identity~\eqref{gr00}
 that each Weyl family is a Nevanlinna family.
In~\cite{Trans} the converse statement was also proven:
each Nevanlinna family can be realized as the Weyl family of a minimal boundary relation.


Two boundary relations $\Gamma^{(j)}:(\sH^{(j)})^2 \to \cH^2$,
$j=1,2$, are said to be \textit{unitarily equivalent} if there is a
unitary operator $U:\sH^{(1)}\to \sH^{(2)}$ such that
\begin{equation}
\label{UNIT1} \Gamma^{(2)}=\left\{ \left\{\begin{pmatrix}Uf \\
Uf'\end{pmatrix},
\begin{pmatrix}h \\ h'\end{pmatrix}\right\}:
\left\{\begin{pmatrix}f \\ f'\end{pmatrix},
\begin{pmatrix}h\\h'\end{pmatrix}\right\}
\in\Gamma^{(1)} \right\}.
\end{equation}
If the boundary relations $\Gamma^{(1)}$ and $\Gamma^{(2)}$
satisfy~\eqref{UNIT1} and $S_j=\ker\Gamma^{(j)}$,
$T_j=\dom\Gamma^{(j)}$, $j=1,2$, then $S_2=US_1U^{-1}$ and
$T_2=UT_1U^{-1}$.

\begin{theorem}
\label{GBTNP}~\cite{Trans} Let $\Gamma:\sH^2 \to \cH^2$ be a boundary relation
for $S^*$. Then the corresponding Weyl family $M(\cdot)$ belongs to
the class $\wt R(\cH)$.

Conversely, if $M(\cdot)$ belongs to the class $\wt R(\cH)$ then
there exists  a unique (up to unitary equivalence) minimal boundary
relation whose Weyl function coincides with $M(\cdot)$.
\end{theorem}

In the following proposition geometric characterization of boundary
relations, whose Weyl functions belong to certain subclasses of $\wt
R(\cH)$ is given.

\begin{proposition} ~\cite{Trans}
\label{thmA} Let $\Gamma : \sH^2 \to \cH^2$ be a boundary relation
for $S^*$ with the Weyl family $M(\lambda)=\Gamma(\wh\sN_\lambda(T))$. Then:
\begin{enumerate}
\def\labelenumi{\rm (\roman{enumi})}
\item $M(\cdot)\in R(\cH)$ if and only if
      $\mul\Gamma\cap(\{0\}\times\cH)=\{0\}$;

\item $M(\cdot)\in R^s(\cH)$ if and only if
      $\ran\Gamma$ is dense in $\cH^2$;
\item $M(\cdot)\in R[\cH]$ if and only if
      $\Gamma_0(\wh\sN_\lambda(T))=\cH$, $\lambda\in\cmr$;
\item $M(\cdot)\in R^s[\cH]$ if and only if $\mul\Gamma_0=\{0\}$ and
      $\Gamma_0(\wh\sN_\lambda(T))=\cH$, $\lambda\in\cmr$.

\item $M(\cdot)\in R^u[\cH]$ if and only if
      $\ran\Gamma=\cH^2$.
\end{enumerate}
\end{proposition}
The case (ii) is specified in more detail in the following
Proposition.
      \begin{proposition}\label{cormul}
Let $\Gamma:\ \sH^2\to\cH^2$ be a boundary relation for $S^*$ and
let $M(\cdot)=\{\Phi(\cdot), \Psi(\cdot)\}$ be the corresponding
Weyl family. Then
     \begin{equation}\label{mulGamma}
\dim\mul\Gamma=\dim\ker {\sf N}_{\Phi,\Psi}.
 \end{equation}
In particular,
    \begin{equation*}
 \mul\Gamma=\{0\}\,\, \Leftrightarrow\,\, \ker {\sf N}_{\Phi,\Psi}=\{0\}.
     \end{equation*}
                  \end{proposition}
\begin{proof}
Let $T(\lambda_0)$ be a mapping from $\cH$ to $\cH^2$ given by
$T(\lambda)=\begin{pmatrix} \Phi(\lambda) \\
\Psi(\lambda)\end{pmatrix}$. Then $M(\lambda_0)=T(\lambda_0)\cH$. If
$\cH_0:=\ker {\sf N}_{\Phi,\Psi}(\lambda,\lambda)\ne 0$, then
$T(\lambda_0)\cH_0$ is the isotropic subspace of the space
$T(\lambda_0)\cH$ considered as a subspace of the Kre\u{\i}n space
$(\cH^2,J_{\cH})$. Therefore,
\[
T(\lambda)\ker {\sf N}_{\Phi,\Psi}(\lambda,\lambda)=M(\lambda)\cap
M(\lambda)^*.
\]
In view of \cite[Lemma 4.1]{Trans} this yields the
equality~\eqref{cormul}.
\end{proof}

\subsection{Linear   transformations of boundary relations}

Let $\Gamma:\sH^2\to\cH^2$ be a boundary relation for $S^*$ and let
$W$ be a linear relation from the Kre\u{\i}n space $(\cH^2,J_\cH)$
to the Kre\u{\i}n space $(\cK^2,J_\cK)$. When the product $W\Gamma$
is a boundary relation for $S^*$, then the corresponding
$\gamma$-field can be expressed in terms of the $\gamma$-field of
the original boundary relation. The Weyl family for $W\Gamma$ can be
expressed as  a Shmul'yan transform of the original Weyl family.


\begin{proposition}\label{Wlemma}
Let $\Gamma:\sH^2\to\cH^2$ be a boundary relation for $S^*$ with the
$\gamma$-field $\gamma(\lambda)$ and the Weyl family $M(\lambda)$.
Let $W$ be a linear relation from the Kre\u{\i}n space
$(\cH^2,J_\cH)$ to the Kre\u{\i}n space $(\cK^2,J_\cK)$, such that
\begin{equation}\label{gw}
 W\Gamma   \mbox{  is unitary,} \quad \ker W =\{0\}.
\end{equation}
Then:
\begin{enumerate}
\def\labelenumi{\rm (\roman{enumi})}

\item the relation $W\Gamma :\sH^2\to\cK^2$
is a boundary  relation for $S^*$;

\item the $\gamma$-field $\gamma_W(\lambda)$ associated with $W \Gamma$ is given by
\begin{equation}
\label{dfHb}
 \gamma_W(\lambda)=\{\, \{k,f_\lambda\}\in \cK\times\sH :\,
    \wh k=W\wh h,\,\, \{\wh f_\lambda,\wh h\}\in\Gamma \},
    \quad \lambda\in\cmr;
\end{equation}
\item
the corresponding Weyl family $M_W(\lambda)$ is given by Shmul'yan
transform
\begin{equation}
\label{mX}
 M_W(\lambda)=W[M(\lambda)],
  \quad \lambda\in\cmr.
\end{equation}
\end{enumerate}
\end{proposition}

\begin{proof}
(i) This statement is immediate from Proposition~\ref{prop1}, since
$\ker W\Gamma=\ker \Gamma=S$.

(ii) According to \eqref{gfield} the $\gamma$-field associated to
$W\Gamma$ is given by
\[
 \gamma_W(\lambda)=
 \{\, \{k,f_\lambda\}\in \cK\times\sH :\,
   \{\wh f_\lambda,\wh k\}\in W \Gamma \}
   =\{\, \{k,f_\lambda\}\in \cK\times\sH :\,
    \wh k=W\wh h,\,\, \{\wh f_\lambda,\wh h\}\in\Gamma \},
\]
$\lambda\in\cmr$, which gives \eqref{dfHb}.

(iii) By Definition~\ref{Weylfam} the Weyl family $M_W(\cdot)$ of
$S$ corresponding to the boundary relation
$\Gamma_W:\sH^2\mapsto\cH^2$ is given by
\[
 M_W(\lambda) =
  \left\{\,\wh k\in\cK^2:\,\{\wh f_\lambda,\wh k\}\in W\Gamma \,\right\}
  =\left\{\,\wh k\in\cK^2:\, \wh k=W\wh h,\,\,
  \{\wh f_\lambda,\wh h\}\in\Gamma \,\right\},
\]
$\lambda \in \cmr$, and this leads to \eqref{mX}.
\end{proof}

In order to guarantee that the product $W\Gamma$ in \eqref{gw} is
unitary, some sufficient conditions on $W$ are required when
$\Gamma$ is not a standard unitary operator; cf.
Section~\ref{preli}.

\begin{remark}
If $W$ is a standard unitary operator from $\cH^2$ to $\cK^2$, then
the conditions in \eqref{gw} are automatically satisfied, and the
$\gamma$-field and the Weyl function of the boundary relation
$W\Gamma$ can be written as
\begin{equation}
\label{dfHb0}
 \gamma_W(\lambda)=
 \left\{\, \{W_{00}h+W_{01}h',\gamma(\lambda)h\}:\,
 \{h,h'\}\in M(\lambda) \,\right\},
 \quad \lambda\in\cmr,
\end{equation}
and
\begin{equation}
\label{mX0}
 M_W(\lambda)
 =\{\,
 \{W_{00}h+W_{01}h',W_{10}h+W_{11}h'\}:\,
 \{h,h'\}\in M(\lambda)\,\},
 \quad \lambda\in\cmr,
\end{equation}
where $W$ is decomposed as in \eqref{W00}. When $W=J_\cH$ the
boundary relation $W\Gamma$  takes the form
\begin{equation}
\label{gaminv}
 \Gamma^{\top}:=\Gamma_{J_\cH }=
\left\{\, \left\{\wh f,J_\cH \wh h\right\}:\, \{\wh f,\wh h\}
\in\Gamma \,\right\},
\end{equation}
and is called the {\it transposed boundary relation}. As follows
from~\eqref{mX0} the corresponding Weyl family $M^\top(\cdot)$ for
$\Gamma^\top$ coincides with $-M(\cdot)^{-1}$.
\end{remark}

\begin{theorem}
\label{INVcl} Let $W$ be a standard unitary operator in the
Kre\u{\i}n space $(\cH^2,J_\cH)$. The following classes of
Nevanlinna families are invariant under the Shmul'yan transform
induced by $W$:
 \begin{enumerate}
\def\labelenumi{\rm (\roman{enumi})}
\item
the class $\wt R(\cH)$ of all Nevanlinna families;
\item
the class $\wt R^c(\cH)$ of constant Nevanlinna families;
\item
the class  $R^s(\cH)$  of strict Nevanlinna functions;
\item
the class  $R^u(\cH)=R^u[\cH]$  of   uniformly strict Nevanlinna
functions.
\end{enumerate}
\end{theorem}

\begin{proof}
(i) According to Theorem~\ref{GBTNP} every Nevanlinna family
$M\in\wt R(\cH)$ admits a realization as a Weyl family of a boundary
relation $\Gamma$. By Lemma~\ref{Wlemma} the linear fractional
transform $\wt M=WM$ is the Weyl family of the boundary relation
$W\Gamma$ and, therefore by Theorem~\ref{GBTNP}, $\wt M$ belongs to
$\wt R(\cH)$.

(ii) Clearly, if $M(\lambda)$ does not depend on $\lambda\in\cmr$,
the same is true for $\wt M=WM$.

(iii)\&(iv) Since $W$ is a unitary operator in the Kre\u{\i}n space
$(\cH^2,J_\cH)$, $W$ maps $\cH^2$ onto $\cH^2$ and, furthermore, as
a bounded everywhere defined unitary operator it also maps a dense
subspace of $\cH$ onto a dense subspace of $\cH$. Therefore, (iii)
and (iv) follow from parts (v) and (ii) of Proposition~\ref{thmA},
respectively.
\end{proof}

The invariance property for the class $R^u(\cH)=R^u[\cH]$ in (iv) of
Theorem~\ref{INVcl} was proved in a completely different manner by
M.G.~Kre\u{\i}n and Yu.~L.~Shmul'yan~\cite{KSh}; the present proof
reflects the power of the realization in Theorem~\ref{GBTNP}.


\section{Special boundary relations and their Weyl families}

In this section special attention is paid to the  boundary relations
whose Weyl families belong to the class $R[\cH]$. In particular, an
orthogonal decomposition of the auxiliary space $\cH$ leads to Weyl
functions of intermediate extensions. Furthermore, attention is paid
to the subclasses $R^s[\cH]$ and $R^u[\cH]$ of strict and uniformly
strict Nevanlinna functions.

\subsection{Ordinary boundary triplets}

{\rm(\cite{GG})} Let $S$ be a closed symmetric operator in a Hilbert
space $\sH$ with equal defect numbers. A triplet
${\Pi}=\{{\cH},{\Gamma}_0,{\Gamma}_1\}$, where $\cH$ is a Hilbert
space with $\dim {\cH}=n_{\pm}(S)$ and $\Gamma_i\in[S^*,\cH]$,
$i=0,1$, is said to be an \textit{ordinary boundary triplet} (or a
\textit{boundary value space}) for $S^*$  if:
\begin{enumerate}
\def\labelenumi{\rm (\roman{enumi})}
\item[(A1)]
the abstract Green's identity  \begin{equation}\label{Greendef1}
  (f', g)-(f, g')
  = (\Gamma_1\wh f, \Gamma_0\wh g)_{\cH}- (\Gamma_0\wh f, \Gamma_1\wh g)_{\cH}
  \end{equation}
holds for all $\wh f=\{f,f'\},\wh g=\{g,g'\} \in  S^*$;
\item[(A2)]
the mapping $\Gamma:=\{\Gamma_0,\Gamma_1\}:\,S^*\to \cH^2$ is
surjective.
\end{enumerate}

For a densely defined symmetric operator this notion was introduced by A.N.
Kochubei~\cite{Koc}  (see also~\cite{GG}), a close definition has been used for other
purpose
by A.V. \v{S}traus~\cite{Str3}. 
For a nondensely defined symmetric
operator it was introduced in~\cite{MM2}.  In this case the adjoint $S^*$ of a symmetric
operator $S$ in $\sH$ is a closed linear relation in $\sH$; it can be considered as a
Hilbert space with the graph norm.

Simple observations (see~\cite[Proposition~5.3]{Trans}) show that
the following statement holds.

\begin{proposition}\label{OBTpr}
The following statements are equivalent:
\begin{enumerate}
\def\labelenumi{\rm (\roman{enumi})}

\item a triplet $\{\cH,\Gamma_0,\Gamma_1\}$ is an ordinary boundary
triplet for $S^*$;

\item $\Gamma=\{\Gamma_0,\Gamma_1\}:\  \sH^2\mapsto\cH^2$ is a
boundary relation for $S^*$ such that  $\ran \Gamma=\cH^2$;

\item the corresponding Weyl family $M(\cdot)$ belongs to
$R^u[\cH]$.

\end{enumerate}
\end{proposition}

A linear extension $\wt A$ of the operator $S$ is said to be
\textit{intermediate} if $S \subset \wt A \subset S^*$. Ordinary
boundary triplets provide a means to describe all intermediate
extensions of $S$. It is well-known (see \cite{DM1, MM2}) that the
set of all intermediate extensions of $A$ in $\sH$ admits the
parametrization
\begin{equation}
\label{+new} \wt A_\Theta:=\{\wh f\in A^*:\,\Gamma\wh
f\in\Theta\}=\ker (\Gamma_1-\Theta \Gamma_0)
\end{equation}
where  $\Theta$ ranges over the set of all linear relations in
$\cH$. Moreover, in this case the linear relation $\wt A_\Theta$ is
closed (symmetric, selfadjoint) if and only if the linear relation
$\Theta$ is closed (symmetric, selfadjoint, respectively).


The definitions of the Weyl function $M(\cdot)$ and the
$\gamma$-field $\gamma(\cdot)$   corresponding to the ordinary
boundary triplet $\Pi=\{\cH,\Gamma_0,\Gamma_1\}$ can be rewritten in
a simpler form
\begin{equation}
\label{GWeyl}
 \wh{\gamma}({\lambda}):= (\Gamma_0\uphar \wh
 {\sN}_{\lambda})^{-1}, \quad {\gamma}({\lambda}):=
 \pi_1(\Gamma_0\uphar \wh{\sN}_{\lambda})^{-1}, \quad
 M({\lambda})=\Gamma_1{\wh{\gamma}}({\lambda}),
\end{equation}
with $\lambda\in\rho(A_0)$. Here $\wh
{\sN}_{\lambda}:=\wh{\sN}_{\lambda}(S^*)$ and $\pi_1$ stands for the
projection onto the first component of $\cH\oplus\cH$. The Weyl
function $M(\cdot)$ and the $\gamma$-field $\gamma(\cdot)$  satisfy
the following identities:
\begin{equation}
\label{ggam3}
 \gamma(\lambda)=[I+(\lambda-\mu)(A_0-\lambda)^{-1}]\gamma(\mu),
 \quad
 \lambda,\mu\in\dC_+ (\dC_-).
\end{equation}
\begin{equation}
\label{orep0}
  M(\lambda)=M(\mu)^*+(\lambda-\bar{\mu})
        \gamma(\mu)^*[I+(\lambda-\mu)(A_0-\lambda)^{-1}]\gamma(\mu),
  \quad
 \lambda, \mu \in \cmr,
\end{equation}

In the case of an ordinary boundary triplet the resolvent of an
intermediate extension $\wt A$ of $A$ can be calculated in terms of
the corresponding Weyl function.

\begin{proposition}\label{PropRez}
{\rm \cite{DM2}} Let  $\{\cH,\Gamma_0,\Gamma_1\}$ be an ordinary
boundary triplet for $S^*$,  let $M(\cdot)$ be the corresponding
Weyl function, let $\Theta$ be a linear relation in $\cH$, and let
$\lambda\in\rho(A_0)$. Then $\lambda\in\rho(\wt A_{\Theta})$ if and
only if $0\in\rho({\Theta}-M(\lambda))$ and the resolvent of $\wt
A_{\Theta}$ is given by
\begin{equation}
\label{krein00}
 (\wt A_{\Theta}-\lambda)^{-1}=(A_0-\lambda)^{-1}
   +\gamma(\lambda)(\Theta-M(\lambda))^{-1}\gamma(\bar{\lambda})^*.
\end{equation}
\end{proposition}

The following result is well known. However, the very simple proof
is here derived from the definition of boundary triplet.

\begin{proposition}\label{uneq}
Two ordinary boundary triplets
$\Pi^{(j)}=\{\cH,\Gamma_0^{(j)},\Gamma_1^{(j)}\}$ ($j=1,2$) for
$A^*$ are connected via the formula
\begin{equation}
\label{Lftr} \Gamma^{(2)}=W\Gamma^{(1)}, \qquad
W=\begin{pmatrix}W_{00} & W_{01}\\W_{10} & W_{11}\end{pmatrix},
\end{equation}
where $W$ is a $J_\cH$-unitary operator in $[\cH^2]$.
\end{proposition}

\begin{proof}
By Proposition~\ref{OBTpr} ordinary boundary triplets determine
single-valued unitary relations $\Gamma^{(j)}:\sH^2\to\cH^2$
satisfying $\dom \Gamma^{(j)}=A^*$ and $\ran \Gamma^{(j)}=\cH^2$.
The composition $W:=\Gamma^{(2)}\circ {\Gamma^{(1)}}^{-1}$ is a
bounded unitary mapping from $\cH^2$ onto $\cH^2$ such that
$\Gamma^{(2)}=W\Gamma^{(1)}$ (see Theorem~\ref{product}). The
statement is now immediate from Lemma~\ref{Wlemma}.
\end{proof}


\subsection{Ordinary boundary triplets and boundary relations.}
\label{OBTandLFT}

Let $S$ be a closed symmetric relation in $\sH$ with equal defect
numbers. It turns out that all the boundary relations of $S$ can
be obtained by extending Proposition~\ref{uneq} in an appropriate
manner. Namely, it is shown that they naturally arise when the
group of $J$-unitary operators in $[\cH^2]$ is augmented by the
class of all $J$-unitary relations $W$ in $\cH^2$ for which $\ker
W=\{0\}$.

 \begin{proposition}\label{WGamma}
Let $S$ be a closed symmetric relation in $\sH$ with equal defect numbers,
let $\Gamma^{(1)}=\{\cH,\Gamma_0^{(1)},\Gamma_1^{(1)}\}$ be an
ordinary boundary triplet for $S^*$ and let $W$ be a unitary
relation from $(\cH^2,J_\cH)$ to $(\wt\cH^2,J_{\wt\cH})$ such that
$\ker W=\{0\}$. Then the linear relation
\begin{equation}\label{Wconnect}
     \Gamma^{(2)}=W\Gamma^{(1)}
\end{equation}
is a boundary relation for $S^*$.

Conversely, for every boundary relation $\Gamma^{(2)}$ of $S^*$
there exists a unitary relation $W$ with $\ker W=\{0\}$ and such
that $\Gamma^{(2)}$ can be recovered from $\Gamma^{(1)}$
via~\eqref{Wconnect}.

In particular, the formula \eqref{Wconnect} establishes a
bijective correspondence between the set of all boundary relations
for $S^*$ with the fixed parameter space $\cH$ and the set of all
unitary relations $W$ in $(\cH^2,J_\cH)$ with $\ker W=\{0\}$.
         \end{proposition}
             \begin{proof}
Let $W$ be a unitary relation from $(\cH^2,J_\cH)$ to
$(\wt\cH^2,J_{\wt\cH})$ such that $\ker W=\{0\}$ and let
$\Gamma^{(2)}$ be given by \eqref{Wconnect}. Then
$\ker\Gamma^{(2)}=\ker\Gamma^{(1)}=S$. Moreover, since $\dom
\Gamma^{(1)}=S^*$ is closed and $\dom W\subset
\ran\Gamma^{(1)}=\cH^2$, part (v) of Theorem~\ref{product} shows
that $\Gamma^{(2)}$ is a unitary relation from $(\sH^2,J_\cH)$ to
$(\wt\cH^2,J_{\wt\cH})$. Since $\ker\Gamma^{(2)}=S$,
$\Gamma^{(2)}$ is a boundary relation for $S^*$.

Conversely, let $\Gamma^{(2)}:(\sH^2,J_\cH) \to
(\wt\cH^2,J_{\wt\cH})$ be a boundary relation for $S^*$. Then by
(iv) of Theorem~\ref{product} the linear relation
$W^{-1}:=\Gamma^{(1)}\circ (\Gamma^{(2)})^{-1}$ is a unitary
relation from $(\wt\cH^2,J_{\wt\cH})$ to $(\cH^2,J_\cH)$, since
$\dom\Gamma^{(1)}=S^*$ is closed and $\Gamma^{(2)}$ as a boundary
relation for $S^*$ satisfies
\[
 \ran(\Gamma^{(2)})^{-1}=\dom\Gamma^{(2)}\subset S^*=\dom
 \Gamma^{(1)}.
\]
Assume that $h\in\ker W$. Then $\{h,0\}\in W=\Gamma^{2}\circ
(\Gamma^{(1)})^{-1}$ and, hence, there is a vector $g\in\sH^2$
such that
\[
\{h,g\}\in(\Gamma^{(1)})^{-1},\quad\{g,0\}\in\Gamma^{(2)}.
\]
Since $\ker\Gamma^{(2)}=S$, this implies that $g\in S$ and
\[
 \{g,h\}\in\Gamma^{(1)}.
\]
Since $\Gamma^{(1)}=\{\cH,\Gamma_0^{(1)},\Gamma_1^{(1)}\}$ is an
ordinary boundary triplet for $S^*$ this implies $h=0$. This shows
that $\ker W=\{0\}$, and completes the proof of the converse
statement.
\end{proof}

\begin{remark}\label{ADD}
  (i)\ A relation $\Gamma^{(2)}$ from $(\sH^2,J_\sH)$ to $(\cH^2,J_\cH)$
   is a single-valued boundary relation for $S^*$ if and only if $W$ in
   \eqref{Lftr} is a unitary operator in
   the parameter space $(\cH^2,J_\cH)$ with $\ker W=\{0\}$.

 (ii)\ If $\dim\cH<\infty$ then $\dom W$ and,
  therefore, also $\ran W$ is closed, and one has
\[
  \dim (\dom W)+\dim (\ker W)=\dim \cH^2, \quad \dim (\ran
  W)+\dim(\mul W)=\dim \wt\cH^2.
\]
If $W$ is single-valued then the assumption $\ker W=\{0\}$ is
equivalent to the fact that $W$ is a standard unitary operator in
$[\cH^2,\wt\cH^2]$, in which case $\dim \wt\cH=\dim \cH$; cf.
Corollary~\ref{cor1}. If instead $W$ is a unitary relation from
$\cH^2$ to $\wt\cH^2$ with $\ker W=\{0\}$, then $\dim\cH\le\dim\wt\cH$
and
    $$
    \dim\wt\cH - \dim\cH = \dim\mul W.
    $$

(iii)\ If $\cH_1 = \cH_2=: \cH$, then $\dim\ker W=\dim \mul W$.
Therefore if $\ker W\neq \{0\}$, then $\mul W \neq \{0\}$ and
$\Gamma^{(2)}$ is a multi-valued mapping. In this case
$\Gamma^{(2)}$ is a boundary relation for $S^*_1$, where
$S_1:=\ker\Gamma^{(2)}\supset A,$ and $S_1\not =A$.
        \end{remark}

  \begin{corollary}\label{WeakWGamma}
Let $S$ be a closed symmetric relation in $\sH$ with equal defect
numbers $(n,n)$, $n\leq\infty$, and let $\Gamma^{(1)}:\sH^2\to\cH^2$
be an ordinary boundary triplet for $S^*$. Then the class of all
single-valued boundary relations $\Gamma^{(2)}:\sH^2\to\cH^2$ for
$S^*$ satisfying $\dom \Gamma^{(2)}=S^*$ coincides with the class of
ordinary boundary triplets for $S^*$; they are parametrized by the
class of all standard unitary operators $W\in[\cH^2]$ via
\eqref{Wconnect}.

Furthermore, if $n=\dim\cH<\infty$ then the class of all boundary
relations for $S^*$ with the fixed parameter space $\cH$ coincides with
the class of ordinary boundary triplets for $S^*$.
  \end{corollary}
 \begin{proof}
Let $\Gamma^{(2)}:\sH^2\to\cH^2$ be a single-valued boundary relation
for $S^*$ with $\dom\Gamma^{(2)}=S^*$. Then
$\ran \Gamma^{(2)}$ is closed and dense in $\cH^2$,
so that in \eqref{Wconnect}
$\cH^2=\ran \Gamma^{(2)}\subset\ran W$, and consequently
$W$ is a standard unitary operator in $\cH^2$. Therefore,
$\Gamma^{(2)}$ is an ordinary boundary triplet for $S^*$.

If the defect numbers of $S$ are finite, then in \eqref{Wconnect}
$\dim \mul\Gamma^{(2)}=\dim\mul W=0$ by Remark~\ref{ADD}.
Therefore, $\Gamma^{(2)}$ is single-valued and since
$n<\infty$, $\dom\Gamma^{(2)}$ is closed, which means that
$\dom\Gamma^{(2)}=S^*$. Therefore, by the first part of the proof
$\Gamma^{(2)}$ is an ordinary boundary triplet for $S^*$.
 \end{proof}

The next result gives a complete description of all the Weyl families
of a symmetric relation $S$ with equal finite defect numbers
in an arbitrary parameter space $\cK$.


\begin{proposition} \label{prGGamma}
Let $S$ be a closed symmetric relation in $\sH$ with equal defect numbers
$(n,n)$, $n\leq\infty$, let $\Gamma^{(1)}:\sH^2\to \cH^2$ be an ordinary
boundary triplet for $S^*$ with the Weyl family $M(\lambda)$, and let
$\Gamma^{(2)}=W\Gamma^{(1)}:\sH^2\to \cH^2$ be an arbitrary boundary
relation for $S^*$. Then the Weyl family associated with
$\Gamma^{(2)}$ is the Shmul' yan transform of $M(\lambda)$ under $W$:
$M^{(2)}(\lambda)=W[M(\lambda)]=:M_W(\lambda)$.

Furthermore, the class of all Weyl families of boundary relations
$\wt\Gamma^{(2)}:\sH^2\to\cK^2$ of $S^*$ in the parameter space
$\cK$ with $\dim\cK=\dim\cH$ are unitarily equivalent to the
class of all Weyl families $M_W(\lambda)$ of $S$ acting on $\cH$,
and they are connected to each others by unitarily equivalent
Shmul' yan transforms.

If $\dim\cK> \dim\cH$ (so that $\dim \cH<\infty$), then the strict
part $M_{r}(\lambda)$ of $M(\lambda)$ in $\cK$ is unitarily
equivalent to a Weyl function of $S$ in the parameter space $\cH$.
\end{proposition}

\begin{proof}
The first assertion $M^{(2)}(\lambda)=W[M(\lambda)]$ is clear from
Proposition~\ref{WGamma} and the definition of the Weyl family.

Next consider boundary relations $\wt\Gamma^{(2)}:\sH^2\to\cK^2$ of
$S^*$ in the parameter space $\cK$ with $\dim\cK=\dim\cH$. Let $U$
be a unitary mapping from the Hilbert space $\cH$ onto the Hilbert
space $\cK$ and define the unitary mapping $\wt U$ from
$\cH\oplus\cH$ onto $\cK\oplus\cK$ by $\wt U=U\oplus U$. Then $\wt
U$ is also a standard unitary mapping from $(\cH^2,J_\cH)$ onto
$(\cK^2,J_\cK)$ in the Kre\u{\i}n space sense. Hence
$\wt\Gamma^{(2)}$ is represented in the form $\wt\Gamma^{(2)}=\wt
U\Gamma^{(2)}$, where $\Gamma^{(2)}:\sH^2 \to\cH^2$ is a boundary
relation for $S^*$ and clearly all boundary relations for $S^*$ with
the parameter space $\cK$ are obtained in this way. The
corresponding Weyl families are connected by $\wt
M^{(2)}(\lambda)=\wt U M^{(2)}(\lambda)$, and hence they are
unitarily equivalent. Moreover, if $\wt\Gamma^{(1)}=\wt U
\Gamma^{(1)}$ and one defines $\wt W=\wt U W \wt U^{-1}$, where $W$
is as in Proposition~\ref{WGamma}, then $\wt W$ is a unitary
relation in $(\cK^2,J_\cK)$ with $\ker\wt W=\{0\}$. Moreover, the
corresponding Shmul'yan transforms are unitarily equivalent: $\wt
W[\wt M(\lambda)]=\wt U W[M(\lambda)]$.

If $\dim \cK>\dim\cH$ with $n=\dim\cH<\infty$ then $S$ has finite
defect numbers. Consequently, $\dom\wt\Gamma^{(2)}=S^*$,
$\ran\wt\Gamma^{(2)}$ is closed, and for $\wt W$ as in
Proposition~\ref{WGamma} one has
$\dom \wt W=\cH^2=M(\lambda)\hplus M(\lambda)^*$, $\lambda\in\cmr$,
since $\Gamma^{(1)}$ is an ordinary boundary triplet for $S^*$. The
codimension of $\mul\wt\Gamma^{(2)}$ is $2n$ in
$\ran\wt\Gamma^{(2)}$ and $n$ in $\wt\Gamma^{(2)}(\wh
\sN_\lambda(S^*))$ for all $\lambda\in\cmr$. Since
\[
 \wt M^{(2)}(\lambda)\cap (\{0\}\oplus \cK)
 =\{0\}\oplus \mul \wt M^{(2)}(\lambda)=\mul\wt\Gamma^{(2)},
\]
one concludes that the codimension of
$\mul\wt\Gamma_0^{(2)}$ in $\dom\wt M^{(2)}(\lambda)$
is also $n$; see \cite[Lemma~4.1]{Trans}. Moreover, here
$\mul\wt\Gamma_0^{(2)}=\ker (\wt M^{(2)}(\lambda)-\wt M^{(2)}(\lambda)^*)$
and therefore also the codimension of
$\cmul\wt\Gamma_0^{(2)}$ in $\cdom\wt M^{(2)}(\lambda)$
is $n$ for all $\lambda\in\cmr$. Since
$\cdom\wt M^{(2)}(\lambda)=\cran\wt\Gamma^{(2)}_0$ by
\cite[Corollary~4.3]{Trans}, it follows that
$\ran\wt\Gamma^{(2)}_0=\dom\wt M^{(2)}(\lambda)$ for all
$\lambda\in\cmr$. Thus $\wt M^{(2)}(\lambda)\in\wt R_{inv}(\cK)$
and ``the strict part'' $\wt M^{(2)}_r(\lambda)$ of $\wt M^{(2)}(\lambda)$
acts on an $n$-dimensional subspace of $\cK$. It is obvious that
$\wt M^{(2)}_r(\lambda)$ defines a Weyl function for $S^*$
on an $n$-dimensional parameter space and, hence, it is unitarily
equivalent to a Weyl function of $S$ acting on $\cH$.
\end{proof}




Proposition~\ref{prGGamma} shows that for studying the Weyl families
of arbitrary boundary relations for a given symmetric relation $S$
it is enough to select one parameter space $\cH$ whose dimension is
equal to $n$, the defect numbers of $S$. The strict part of the Weyl
family is the one that determines the symmetric operator and its
selfadjoint extension in the model space up to unitary equivalence.
However, for instance, in the connection of generalized resolvents
nonstrict Weyl families naturally appear.

\subsection{Boundary relations whose Weyl functions belong to the
class $R[\cH]$.}

A purely geometric characterization  of this class of boundary
relations is given in the following proposition.

\begin{proposition}\label{GBTb} Let $S$ be a closed symmetric relation in a Hilbert
space $\sH$. Let $\cH$  be  a Hilbert space and let
$\Gamma:\sH^2\to\cH^2$ be a (possibly multivalued) linear relation
such that:
\begin{enumerate}
\def\labelenumi{\rm (\roman{enumi})}
\item[(B1)]
 Green's identity~\eqref{Green1} holds;
\item[(B2)]
$\ran\Gamma_0=\cH$;
\item[(B3)]
$A_0:=\ker\Gamma_0$ is a selfadjoint linear relation in $\sH$.
\end{enumerate}
Then $\Gamma:\sH^2\to\cH^2$ is a boundary relation for
$S^*:=(\ker\Gamma)^*$
such that
\begin{equation}
\label{RanG01}
 \Gamma_0(\wh{\sN}_\lambda(T))=\cH, \quad \lambda\in\cmr.
\end{equation}
Conversely, every closed isometric linear relation
$\Gamma:\sH^2\to\cH^2$ satisfying~\eqref{RanG01} satisfies also the
conditions (B1)--(B3).

If the conditions (B1)--(B3) are satisfied, the corresponding Weyl
function belongs to the class $R[\cH]$. Moreover, every
$R[\cH]$-function is the Weyl function of some boundary relation
$\Gamma:\sH^2\to\cH^2$ with the properties (B1)--(B3).
\end{proposition}

\begin{proof}
The proof of the direct statement was given in~\cite{Trans}.

Assume now that $\Gamma:\sH^2\to\cH^2$ is a closed isometric linear relation
satisfying~\eqref{RanG01}.
Let $\{\wh f_\lambda,\wh h\}$, $\{\wh g_{\bar\lambda},\wh
k\}\in\Gamma$ with
\[
\wh f_\lambda=\left(\begin{array}{c}
   f_\lambda \\
  \lambda f_\lambda \\
\end{array}\right)\in\wh\sN_\lambda(T),
\quad
\wh g_{\bar\lambda}=\left(\begin{array}{c}
   g_{\bar\lambda} \\
  {\bar\lambda}g_{\bar\lambda} \\
\end{array}\right)\in\wh\sN_{\bar\lambda}(T),
\quad
\wh h=\left(\begin{array}{c}
   h \\
  h' \\
\end{array}\right),
\quad
\wh k=\left(\begin{array}{c}
   k \\
  k' \\
\end{array}\right)\in\cH^2.
\]
Then it follows from~\eqref{gr00} that
\[
0=(\lambda f_\lambda,g_{\bar\lambda})_\sH-( f_\lambda,\bar\lambda
g_{\bar\lambda})_\sH =(h',k)_\cH-(h,k')_\cH,
 \quad
 \lambda \in \cmr.
\]
Since $\wh h\in M(\lambda)$, $\wh k\in M(\bar\lambda)$ this implies
that $M(\lambda)\subset M(\bar\lambda)^*$. Next,   the
assumption~\eqref{RanG01},  implies that
\[
\dom M(\lambda)=\dom  M(\bar\lambda)=\cH.
\]
and, hence, $M(\lambda)$ is bounded for all $\lambda\in \cmr$.
Since the operator $M(\lambda)$ is dissipative
for $\lambda\in\dC_+$ this implies
\[
\ran(\Gamma\wh\sN_\lambda+\lambda)=\ran(M(\lambda)+\lambda)=\cH,\quad\lambda\in\dC_+.
\]
Due to Proposition~\ref{thm1} this proves that $\Gamma : \sH^2 \to
\cH^2$ is a boundary relation for $S=T^*$.

Thus $\cH=\Gamma_0(\wh\sN_\lambda(T))\subset\ran \Gamma_0$,
so that $\ran \Gamma_0=\cH$, i.e., (B2) is satisfied.
Also the property (B3) is obtained from
$\Gamma_0(\sN_\lambda(T))=\cH$ by using~\cite[Proposition~4.15]{Trans}.
The condition (B1) for
the boundary relation $\Gamma$ is clearly satisfied.

The fact that every $R[\cH]$-function is the Weyl function of
some boundary relation $\Gamma:\sH^2\to\cH^2$
satisfying the conditions (B1)-(B3) is implied
 by Theorem~\ref{GBTNP} and Proposition~\ref{thmA}.
\end{proof}

Recall that for a boundary relation $\Gamma:\sH^2\to\cH^2$
satisfying the conditions (B1)-(B3) the operator function
$\gamma(\lambda)=\pi_1(\Gamma_0\uphar
\wh\sN_\lambda(T))^{-1}:\cH\to\sN_\lambda(T)$ is bounded and
single-valued for every $\lambda\in\cmr$, see~\cite{Trans}. Clearly,
the Weyl function $M(\cdot)$ and the $\gamma$-field $\gamma(\cdot)$
satisfy the identities~\eqref{ggam3} and~\eqref{orep0}. Let $E(t)$
be the spectral family of $A_0$ and let $P=E(\infty)$ be the
orthogonal projection onto $\cdom A_0$. Then~\eqref{orep0} leads to
the following integral representation of $M(\lambda)$
\begin{equation}
\label{inrep0}
 (M(\lambda)h,h)=a_h+b_h\lambda
 +\int_{\dR}\left(\frac{1}{t-\lambda}-\frac{t}{t^2+1}\right)\,d\sigma_h(t),
 \quad h\in\cH_0,
\end{equation}
where
\[
 a_h=(\RE M(i)h,h)_\cH,\,\, b_h=((I-P)\gamma(i)h,\gamma(i)h),\,\,
 d\sigma_h(t)=(t^2+1)d(E(t)P\gamma(i)h,P\gamma(i)h)_\sH.
\]
The representation~\eqref{inrep0} leads to the following
characterization.

\begin{proposition}\label{MULA0T}
Let  $S$ be a symmetric operator in $\sH$.  Let
$\Gamma:\sH^2\to\cH^2$ be a boundary relation for $S^*$ satisfying
the conditions (B1)-(B3) and let $M(\lambda)$ be the corresponding
Weyl function. Let  $\cH_0=\pi_1\mul\Gamma$, $A_0=\ker\Gamma_0$, and
$T=\dom\Gamma$.  Then:
\begin{enumerate}
\def\labelenumi{\rm (\roman{enumi})}
\item $\mul A_0=\{0\}$  if and only if
\begin{equation}
\label{mulA0}
 \lim_{y \to \infty} \frac{(M(iy)h,h)_\cH}{iy}=0,
 \quad
 h \in \cH;
\end{equation}

\item $\mul T=\{0\}$ if and only if $M$ satisfies the
condition~\eqref{mulA0} and
\begin{equation}
\label{mulT} \lim_{y\uparrow\infty}y \IM \bigl(M(i y)h,h\bigr) =\infty, \qquad
h\in\cH\ominus\cH_0.
\end{equation}
\end{enumerate}
\end{proposition}

\begin{proof} The first statement is immediate from the equality
\begin{equation}
\label{BM}
 \lim_{y \to \infty} \frac{(M_s(iy)h,h)_\cH}{iy}
 =\|(I-P)\gamma(i)h\|_\sH^2
 =\|(I-P)\gamma(\mu)h\|_\sH^2.
\end{equation}
Under the assumption~\eqref{mulA0} the limit in~\eqref{mulT} takes the form
\[
\lim_{y\uparrow\infty}y \IM \bigl(M(i y)h,h\bigr)= \int_\dR
 (t^2+1)\,d\|E_t\gamma(i)h\|_\sH^2.
\]
Remark that the mapping $\gamma(i)$ restricted to $\cH\ominus\cH_0$
is injective and, hence, this limit is finite for some
$h\in\cH\ominus\cH_0$, $h\ne 0$, if and only if $\sN_i(T)\cap\dom
A_0=(A_0-\lambda)^{-1}(\mul T)$ is nontrivial. For the proof of the
last equality see~\cite{Trans}.
\end{proof}

The boundary relations with the additional properties (B1)--(B3) are
invariant under a special class of transforms, cf. Proposition
\ref{Wlemma}. Let $B\in[\cH]$ and let $G\in[\cH]$ be invertible, and
assume that
\[
BG=(BG)^*.
\]
Define the block operator $\wt W$ by
\begin{equation}
\label{Wlin} \wt W=\begin{pmatrix} G^{-1} & 0 \\ B & G^*
\end{pmatrix},\quad \mbox{ with }BG=(BG)^*.
\end{equation}
It is easy to see that $\wt W$ is a $J_\cH$-unitary operator in
$\cH^2$.

\begin{proposition}
\label{linprop} Let $\Gamma:\sH^2\to\cH^2$ be a boundary relation
for $S^*$ which satisfies the conditions (B1)--(B3), let
$\gamma(\lambda)$ and $M(\lambda)$ be the corresponding
$\gamma$-field and the Weyl function, and moreover let $\wt W\in
[\cH\oplus\cH]$ be given by \eqref{Wlin}. Then:
\begin{enumerate}
\def\labelenumi{\rm (\roman{enumi})}

\item the transform $\wt \Gamma=\wt W\Gamma$ of $\Gamma$ given by
\begin{equation}
\label{lintrans}
 \wt\Gamma=\left\{\,
         \left\{\wh f,\begin{pmatrix} G^{-1}h \\Bh+G^*h' \end{pmatrix}\right\}:\,
         \left\{\wh f,\wh h \right\}\in\Gamma \,\right\},
\end{equation}
is a boundary linear relation for $S^*$ with
$\dom\Gamma_W=\dom\Gamma$ and $\ker\Gamma_W=\ker\Gamma=S$ which also
satisfies the conditions (B1)--(B3);

\item the $\gamma$-field and the Weyl function associated to $\wt\Gamma$ are given by
\begin{equation}
\label{gmlin}
 \wt\gamma(\lambda)=\gamma(\lambda)G,
 \quad \wt M(\lambda)=BG+G^*M(\lambda)G\,(\in [\cH]),
  \quad \lambda\in\cmr.
\end{equation}
\end{enumerate}
\end{proposition}

\begin{proof}
(i) Since $\wt W$ defined by \eqref{Wlin} is a $J_\cH$-unitary
operator in $\cH$, the transform $\wt\Gamma=\wt
W\Gamma:\sH^*\to\cH^2$ is a boundary relation for $S^*$ with
$\dom\Gamma_W=\dom\Gamma$ and $\ker\Gamma_W=\ker\Gamma=S$ and
clearly $\wt\Gamma$ admits the representation \eqref{lintrans}.
Moreover, since $\ran\Gamma_0=\cH$ and $G\in[\cH]$ is invertible,
the equality $\ran\wt\Gamma_0=\cH$ holds and
$\ker\wt\Gamma_0=\ker\Gamma_0$ is selfadjoint. Hence, $\wt\Gamma$
satisfies the conditions (B1)--(B3).

(ii) By Lemma~\ref{Wlemma} the Weyl function $\wt M(\lambda)$
associated to $\wt\Gamma$ is given by
\begin{equation}
\label{mlin}
\begin{split}
 \wt M(\lambda)
 &=\{\, \{G^{-1}h,Bh+G^*h'\}:\, \{h,h'\}\in M(\lambda)\,\} \\
 &=\{\, \{k,BGk+G^*M(\lambda)Gk\}:\, h=Gk\in \dom M(\lambda)=\cH \,\} \\
 &=BG+G^*M(\lambda)G,
\end{split}
\end{equation}
where $BG=(BG)^*$. Similarly, the $\gamma$-field
$\wt\gamma(\lambda)$ corresponding to $\wt\Gamma$ takes the form
\begin{equation}
\label{dfWlin}
 \wt\gamma(\lambda)
 =\left\{\, \{G^{-1}h,\gamma(\lambda)h\}:\, \{h,h'\}\in M(\lambda) \,\right\}
 =\left\{\, \{k,\gamma(\lambda)Gk\}:\, h=Gk\in \dom M(\lambda)=\cH
 \,\right\},
\end{equation}
so that $\wt\gamma(\lambda)=\gamma(\lambda)G$, $\lambda\in\cmr$.
\end{proof}
\begin{remark}
In the case when the transposed boundary relation $\Gamma^\top$
satisfies (B1)--(B3) the corresponding Weyl family
$M^\top(\cdot)=-M(\cdot)^{-1}$ is single-valued and belongs to the
class $R[\cH]$.
\end{remark}
Up to this point boundary relations $\Gamma:(\sH^2,J_\sH)\to(\cH^2,J_\cH)$ satisfying the
conditions (B1)--(B3) are in general multi-valued. Next we briefly discuss the case when
it is single-valued.
      \begin{definition}\cite{DM2}
If  a boundary relation $\Gamma:(\sH^2,J_\sH)\to(\cH^2,J_\cH)$ is
single-valued and satisfies the conditions (B1)--(B3), then  the
triplet $\{\cH,\Gamma_0,\Gamma_1\}$ is said to be a
\textit{generalized boundary triplet}.
        \end{definition}

The following corollary is implied by Proposition~\ref{thmA}. 

         \begin{corollary}\cite{DM2}         
A single-valued boundary relation
$\Gamma=\{\Gamma_0,\Gamma_1\}:\sH^2\mapsto\cH^2$ corresponds to a generalized boundary
triplet $\{\cH,\Gamma_0,\Gamma_1\}$ if and only if the corresponding Weyl function
$M(\cdot)$ belongs to the class $R^s[\cH]$.
        \end{corollary}

In the case of a generalized boundary triplet the last condition in
Proposition~\ref{MULA0T} is simplified in the following way
      \begin{equation}
\label{mulT1} \lim_{y\uparrow\infty}y \IM \bigl(M(i y)h,h\bigr) =\infty, \qquad h\in\cH.
       \end{equation}
The next proposition shows how one can reduce a multi-valued boundary relation $\chi$
with the properties (B1)--(B3) to a single-valued boundary relation with the same
properties (B1)--(B3).
\begin{proposition}
\label{GBTCB} Let $\Gamma:(\sH^2,J_\sH)\to(\cH^2,J_\cH)$ be a
multi-valued boundary relation which satisfies (B1)--(B3). Then:
\begin{enumerate}
\item
$\cH_0=\pi_0\mul\Gamma$ is a closed subspace of $\cH$
($\cH_1=\cH\ominus\cH_0$);
\item
$\mul\Gamma$ is the graph of bounded symmetric operator
$K_0\in[\cH_0,\cH]$;
\item
for every bounded selfadjoint extension $K$ of $K_0$ in $\cH$ the
linear relation
\begin{equation}\label{BR13}
\Gamma':=\left\{\left\{\left(%
\begin{array}{c}
  f \\
  f' \\
\end{array}%
\right),
\left(%
\begin{array}{c}
  P_{\cH_1}h \\
  h'-Kh \\
\end{array}%
\right)\right\}: \,\left\{\left(%
\begin{array}{c}
  f \\
  f' \\
\end{array}%
\right),
\left(%
\begin{array}{c}
h \\
  h' \\
\end{array}%
\right)\right\}\in\Gamma\right\}:\sH^2\to\cH_1^2
\end{equation}
is a single-valued boundary relation satisfying (B1)--(B3). The Weyl
functions $M(\lambda)$ and $M_1(\lambda)$, corresponding to the
boundary relations $\Gamma$ and $\Gamma'$ are connected by
\begin{equation}\label{WeylB13}
M(\lambda)=K+\mbox{diag }(0_{\cH_0},M_1(\lambda)),\quad
(\lambda\in\dC_+).
\end{equation}
\end{enumerate}
 \end{proposition}
\begin{proof}
1) Since $\ran\Gamma_0=\cH$ one obtains from~\cite[Lemma 2.1]{Trans}
that $\ran\Gamma\wh{+}(\{0\}\oplus\cH)$ is closed. By~\cite[Theorem
4.8]{Kato} and Proposition~\ref{UNIT}
\begin{equation}\label{Angle}
\mul\Gamma\wh{+}(\{0\}\oplus\cH)\mbox{ is closed}.
\end{equation}
Using~\cite[Lemma 2.1]{Trans} again one obtains
$\cH_0:=\pi_1\mul\Gamma$ is a closed subspace of $\cH$.

2) It follows from \eqref{Angle} that $\mul\Gamma$ is the graph of a
bounded operator $K_0:\cH_0\to\cH$. Since $\mul\Gamma$ is a neutral
subspace in $(\cH^2,J_{\cH})$ the  operator $K_0$ is symmetric in
$\cH$.

3)  Let $K$ be a bounded selfadjoint operator extension of $K_0$,
$K\in[\cH]$. Since $\mul\Gamma=\ran\Gamma^{[\perp]}$, one obtains
from
\[
0=(h',h_0)-(h,K_0h_0)=(h'-Kh,h_0)\quad
(h_0\in\cH_0),\{h,h'\}\in\ran\Gamma,
\]
that $h'-Kh$ is orthogonal to $\cH_0$. This proves that
$\ran\Gamma'\subset\cH_1^2$. The mapping $\Gamma'$ is single-valued
since for $\{h,h'\}\in\mul\Gamma$ one has
\[
P_{\cH_1}h=0,\quad h'-Kh=K_0h-Kh=0.
\]
Clearly, $\ran\Gamma'_0=\cH_1$ since $\ran\Gamma_0=\cH$. Assume that
$\{f,f'\}\in\ker\Gamma_0'$. It means that there is a vector
$h'\in\cH$ such that $\{h,h'\}\in\mul\Gamma$
\begin{equation}\label{A0}
    \left\{\left(%
\begin{array}{c}
  f \\
  f' \\
\end{array}%
\right),
\left(%
\begin{array}{c}
h \\
  h' \\
\end{array}%
\right)\right\}\in\Gamma,\quad P_{\cH_1}h=0.
\end{equation}
Then $h\in\cH_0$ and, hence, there is a vector $h''\in\cH$ such that
$\{h,h''\}\in\mul\Gamma$. Therefore,
\[
    \left\{\left(%
\begin{array}{c}
  f \\
  f' \\
\end{array}%
\right),
\left(%
\begin{array}{c}
  0\\
  h'-h'' \\
\end{array}%
\right)\right\}\in\Gamma,
\]
and, hence, $\{f,f'\}\in\ker\Gamma_0$. This proves that $\Gamma'$
satisfies (B3). The equality~\eqref{WeylB13} is implied
by~\eqref{BR13}.
\end{proof}
\begin{corollary}
\label{CorIndDef} If a boundary relation
$\Gamma:(\sH^2,J_\sH)\to(\cH^2,J_\cH)$ for $S^*$ satisfies
(B1)--(B3), then $n_+(S)=n_-(S)$.
\end{corollary}

\section{Weyl functions for intermediate extensions} \label{Inter}

Let $S$ be a closed symmetric operator in a separable Hilbert space
$\sH$  and let $\Gamma:\sH^2\to\cH^2$ be a boundary relation for
$S^*$ which satisfies the conditions (B1)--(B3), so that the
corresponding Weyl family $M(\lambda)$ belongs to the class
$R[\cH]$. The purpose of this section is to associate  intermediate
symmetric extensions $H$ of $S$ to different types of Nevanlinna
functions (say, linear combinations of $M_{ij}$, Schur complements
and compressions of linear fractional transformations of
$M(\lambda)$), which are obtained as block transforms of the
operator matrix representation of $M(\lambda)$ in
\begin{equation}
\label{decom2}
 \cH=\cH_1 \oplus \cH_2,
\quad
  M(\lambda)=(M_{ij}(\lambda))_{i,j=1}^2.
\end{equation}
Consider the linear relations
\begin{equation}\label{P12}
    \cP^{(j)}=\left\{\left\{\left(%
\begin{array}{c}
  h \\
  h' \\
\end{array}\right),
\left(\begin{array}{c}
  h \\
  P_{\cH_j}h' \\
\end{array}
\right)\right\}:\,h\in\cH_j,\,h'\in\cH\,\right\},\quad j=1,2,
\end{equation}
which, clearly, are unitary from $(\cH^2,J_{\cH})$ to
$(\cH_j^2,J_{\cH_j})$. In general, it is not clear whether
$\cP^{(j)}\circ\Gamma$ is a unitary relation if
$\ran\Gamma\nsubseteq \dom\cP^{(j)}=\cH_j\times \cH$ (cf.
Theorem~\ref{product}). However, in the case when
$\Gamma:\sH^2\to\cH^2$ satisfies the conditions (B1)--(B3) it
turns out that $\cP^{(j)}\circ\Gamma$ is a unitary relation from
$(\cH^2,J_{\cH_j})$ to $(\cH_j^2,J_{\cH_j})$, $(j=1,2)$.

%
%
%

\begin{proposition}
\label{coner} Let $\Gamma:\sH^2\to\cH^2$ be a boundary relation for
$S^*$ which satisfies the conditions (B1)--(B3), let
$\gamma(\lambda)$ be the corresponding $\gamma$-field  and decompose
the corresponding Weyl function $M(\lambda)$ as in \eqref{decom2}.
Then:
\begin{enumerate}
\def\labelenumi{\rm (\roman{enumi})}
\item the linear relation $H_1$ given by
\begin{equation}
\label{h1}
 H^{(1)}=\left\{\, \wh f \in S^*:\, \left\{\wh f,
           \begin{pmatrix} 0 \\ h' \end{pmatrix}\right\} \in\Gamma \,\,
              \mbox{ for some } h'\in \cH_2 \,\right\},
\end{equation}
is closed and symmetric in $\sH$ and has equal defect numbers;

\item the linear relation $\Gamma^{(1)}:\sH^2\to\cH_1^2$ given by
\begin{equation}\label{G1}
 \Gamma^{(1)}:=\cP^{(1)}\circ\Gamma=\left\{\,
         \left\{\wh f,\begin{pmatrix} h \\ P_{\cH_j} h' \end{pmatrix}\right\}:\,
         \left\{\wh f,\begin{pmatrix} h \\ h' \end{pmatrix}\right\}\in\Gamma \,\,
              \mbox{ for some } h \in \cH_1, h'\in \cH \,\right\},
\end{equation}
is a boundary relation for $(H^{(1)})^*$ which satisfies the
conditions (B1)--(B3);

\item the domain $T_1:=\dom\Gamma^{(1)}$ is dense in $H_1^*$ and it can be rewritten
as
\begin{equation}
\label{h1*}
 T^{(1)}=\left\{\, \wh f \in S^*:\, \{\wh f,\wh h\}\in \Gamma, \, \pi_2 h=0 \,\right\};
\end{equation}
\item
the corresponding $\gamma$-field $\gamma_1(\lambda):\cH_1\to\sH$ and
the Weyl function $M_1(\lambda)\in[\cH_1]$ are given by
\begin{equation}
\label{mone}
 \gamma_1(\lambda)=\gamma(\lambda)\uphar \cH_1,
\quad
  M_1(\lambda)=M_{11}(\lambda),
  \quad \lambda\in\cmr.
\end{equation}
\end{enumerate}
\end{proposition}

\begin{proof}
(i) By definition $\Gamma^{(1)}$ is a multi-valued mapping from
$\sH^2$ into $\cH_1^2$. It satisfies the Green's
identity~\eqref{Green1}, since for all $\{\wh f,\wh h\},\, \{\wh
g,\wh k\}\in\Gamma$ with $h,k\in\cH_1$ one has
\[
 (f',g)_\sH-(f,g')_\sH=({h'},k)_{\cH}-(h,k')_{\cH}
  =(\pi_1 h',k)_{\cH}-(h,\pi_1 k')_{\cH}.
\]
The property (B2) of $\Gamma$ implies that $\ran
\Gamma_0^{(1)}=\cH_1$. Moreover, from the property (B3) of $\Gamma$
one concludes that $\ker\Gamma_0^{(1)}=\ker\Gamma_0=A_0$ is
selfadjoint. Hence, $\Gamma^{(1)}$ is a boundary relation for
$(H^{(1)})^*$ which admits the properties (B1)--(B3).

(ii) Since $\Gamma^{(1)}$ is unitary, $H^{(1)}=\ker\Gamma^{(1)}$ is
closed and symmetric. The description of $H^{(1)}$ in \eqref{h1} is
immediate from the definition of $\Gamma^{(1)}$ in (ii). Since
$\Gamma^{(1)}$ satisfies the conditions (B1)--(B3) the defect
numbers of $H^{(1)}$ are equal to $(n_1,n_1)$, $n_1=\dim \cH_1-\dim
\mul\Gamma^{(1)}$.

(iii) The description of $T_1=\dom \Gamma^{(1)}$ in \eqref{h1*} is
clear from the definition of $\Gamma^{(1)}$ in (i) and the denseness
of $T^{(1)}$ in $(H^{(1)})^*$, or equivalently, the identity
$(T^{(1)})^*=\ker\Gamma^{(1)}=H^{(1)}$ holds by the definition of
boundary relations.

(iv) According to \cite[Proposition~5.9]{Trans} the conditions
(B1)--(B3) imply that
\begin{equation}
\label{ggdom}
 \Gamma_0(\wh\sN_\lambda(T)=\cH,\quad
 \Gamma^{(1)}_0(\wh\sN_\lambda(T^{(1)})=\cH_1, \mbox{ for all } \lambda\in\cmr.
\end{equation}
Hence, $\dom \wh\gamma(\lambda)=\cH$ and $\dom
\wh\gamma_1(\lambda)=\cH_1$, and the formulas
\[
\gamma(\lambda)=\{\,\{h,\wh f_\lambda\}:\, \{\wh f_\lambda,\wh
h\}\in \Gamma \,\}, \quad \gamma_1(\lambda)=\{\,\{h,\wh
f_\lambda\}:\, \{\wh f_\lambda,\wh h\}\in \Gamma,\, h\in\cH_1 \,\}
\]
show that these single-valued mappings are connected by
$\gamma_1(\lambda)=\gamma(\lambda)\uphar\cH_1$. Moreover,
\eqref{ggdom} implies that $M_1(\lambda)\in[\cH_1]$,
$M(\lambda)\in[\cH]$, and thus
\[
\begin{split}
 M_1(\lambda)&=\{\,\wh h\in\cH_1^2:\, \{\wh f_\lambda,\wh h\}\in
                \Gamma^{(1)} \,\} \\
  &=\{\,\{h,\pi_1 h'\}\in\cH^2:\, \{\wh f_\lambda,\wh h\}\in
   \Gamma,\, h\in\cH_1 \,\} \\
 &=P_{\cH_1} M(\lambda)\uphar\cH_1.
\end{split}
\]
This completes the proof.
\end{proof}

Replacing $P_{\cH_1}$ by $P_{\cH_2}$ one obtains

\begin{corollary}
\label{cblock} Let $\Gamma:\sH^2\to\cH^2$, $\gamma(\lambda)$, and
$M(\lambda)$ be as in Proposition~\ref{coner}. Then:
\begin{enumerate}
\def\labelenumi{\rm (\roman{enumi})}
\item the linear relation $H_2$ given by
\begin{equation}
\label{h2}
 H^{(2)}=\left\{\, \wh f \in S^*:\, \left\{\wh f,
           \begin{pmatrix} 0 \\ h' \end{pmatrix}\right\} \in\Gamma \,\,
              \mbox{ for some } h'\in \cH_1 \,\right\},
\end{equation}
is closed and symmetric in $\sH$ and has equal defect numbers;

\item the linear relation $\Gamma^{(2)}:\sH^2\to\cH_2^2$ given by
\[
 \Gamma^{(2)}:=\cP^{(2)}\circ\Gamma=\left\{\,
         \left\{\wh f,\begin{pmatrix} h \\ \pi_2 h' \end{pmatrix}\right\}:\,
         \left\{\wh f,\begin{pmatrix} h \\ h' \end{pmatrix}\right\}\in\Gamma \,\,
              \mbox{ for some } h \in \cH_2, h'\in \cH \,\right\},
\]
is a boundary relation for $H_2^*$ which satisfies the conditions
(B1)--(B3);

\item the domain $T^{(2)}:=\dom\Gamma^{(2)}$ is dense in $H_2^*$ and it can be rewritten
as
\begin{equation}
\label{h2*}
 T^{(2)}=\left\{\, \wh f \in S^*:\, \{\wh f,\wh h\}\in \Gamma, \, \pi_1 h=0 \,\right\};
\end{equation}
\item
the corresponding $\gamma$-field $\gamma_2(\lambda):\cH_2\to\sH$ and
the Weyl function $M_2(\lambda)\in[\cH_2]$ are given by
\begin{equation}
\label{mone1}
 \gamma_2(\lambda)=\gamma(\lambda)\uphar \cH_2,
\quad
  M_2(\lambda)=M_{22}(\lambda),
  \quad \lambda\in\cmr.
\end{equation}
\end{enumerate}
\end{corollary}

\begin{corollary}
\label{SchurC} Let $\Gamma:\sH^2\to\cH^2$ be a boundary relation for
$S^*$, such that $\Gamma$, $\Gamma^\top$,
$\left(\Gamma^{(2)}\right)^\top$ satisfy the conditions (B1)--(B3),
and decompose the corresponding Weyl function $M(\lambda)$ as in
\eqref{decom2}.
\begin{enumerate}
\def\labelenumi{\rm (\roman{enumi})}
\item the linear relation $S^{(1)}$ given by
\begin{equation}
\label{s1}
 S^{(1)}=\left\{\, \wh f \in S^*:\, \left\{\wh f,
           \begin{pmatrix} h \\ 0 \end{pmatrix}\right\} \in\Gamma \,\,
              \mbox{ for some } h\in \cH_2 \,\right\},
\end{equation}
is closed and symmetric in $\sH$ and has equal defect numbers;

\item the linear relation 
\[
 \Gamma':=\left(\cP^{(1)}\circ\Gamma^{\top}\right)^\top=\left\{\,
         \left\{\wh f,\begin{pmatrix} P_{\cH_1}h \\ h' \end{pmatrix}\right\}:\,
         \left\{\wh f,\begin{pmatrix} h \\ h' \end{pmatrix}\right\}\in\Gamma \,\,
              \mbox{ for some } h \in \cH, h'\in \cH_1 \,\right\},
\]
is a boundary relation for $(S^{(1)})^*$ which satisfies the
conditions (B1)--(B3);

\item
the corresponding 
Weyl function $M^{(1)}(\lambda)\in[\cH_1]$ is given by
\begin{equation}
\label{moneSchur}
  M^{(1)}(\lambda)=M_{11}(\lambda)-M_{12}(\lambda)M_{22}(\lambda)^{-1}M_{21}(\lambda),
  \quad \lambda\in\cmr.
\end{equation}
\end{enumerate}
\end{corollary}
\begin{proof}
It follows from the assumptions and Proposition~\ref{GBTb} that both
$M(\cdot)$ and $M(\cdot)^{-1}$ belong to the class $R[\cH]$ and
$M(\cdot)^{-1}$ admits a block representation
$M(\cdot)^{-1}=\left((M(\cdot)^{-1})_{ij}\right)_{i,j=1}^2$. Since
$\left(\Gamma^{(2)}\right)^\top$ satisfy the conditions (B1)--(B3),
one obtains from Proposition~\ref{GBTb} that  $M_{22}(\cdot)^{-1}$
belongs to the class $R[\cH_2]$. Then it follows from the Frobenius
formula that
\begin{equation}\label{Frobenius}
\left(M(\lambda)^{-1}\right)_{11}=\left(M_{11}(\lambda)-M_{12}(\lambda)M_{22}(\lambda)^{-1}M_{21}(\lambda)\right)^{-1}.
\end{equation}

Let us apply Proposition~\ref{coner} to the linear relation
$\Gamma^\top$. Then the linear relation
$\cP^{(1)}\circ\Gamma^{\top}$ satisfies the assumptions (B1)--(B3)
and the corresponding Weyl function coincides with
$\left(M(\cdot)^{-1}\right)_{11}$ in~\eqref{Frobenius}. To complete
the proof it remains to show that the linear relation
$\left(\cP^{(1)}\circ\Gamma^{\top}\right)^\top$ satisfies the
assumptions (B1)--(B3). Since the Weyl function
\[
M'(\lambda)=\left(M(\cdot)^{-1}\right)_{11}=M_{11}(\lambda)-M_{12}(\lambda)M_{22}(\lambda)^{-1}M_{21}(\lambda)
\]
belongs to the class $R[\cH_1]$ this fact is implied by
Proposition~\ref{GBTb}. However, we will present also a direct
proof.

Let $h_1\in\cH_1$. Since $\Gamma$ satisfies (B2) there exists $\wh
f\in S^*$, and $h'\in \cH$ such that
\begin{equation}\label{eq:4.1}
    \left\{\wh f,\left(%
\begin{array}{c}
  h_1 \\
  h_1' \\
\end{array}%
\right)\right\}\in\Gamma.
\end{equation}
Next, using the fact that $\left(\Gamma^{(2)}\right)^\top$ satisfies
(B2), we find $\wh g\in S^*$, and $h_2\in \cH_2$, $h_2'\in\cH$, such
that
\begin{equation}\label{eq:4.2}
    \left\{\wh g,\left(%
\begin{array}{c}
  h_2 \\
  h_2' \\
\end{array}%
\right)\right\}\in\Gamma,\quad P_{\cH_2}h_2'=P_{\cH_2}h_1'.
\end{equation}
Now it follows from~\eqref{eq:4.1},~\eqref{eq:4.2} that
\[
    \left\{\wh f-\wh g,\left(%
\begin{array}{c}
  h_1-h_2 \\
  h_1'-h_2' \\
\end{array}%
\right)\right\}\in\Gamma,\quad P_{\cH_2}(h_1'-h_2')=0.
\]
This implies that $\ran\Gamma_0'=\cH_1$ and, hence, $\Gamma'$
satisfies (B2). The assumption (B3) for $\Gamma'$ is implied by the
equality
\[
\ker\Gamma_0'=\left\{\wh f\in S^*:\,\left\{\wh f,\left(%
\begin{array}{c}
  h_1 \\
  h_1' \\
\end{array}%
\right)\right\}\in\Gamma \right\}.
\]

\end{proof}
\begin{proposition}
\label{block2} Let $\Gamma:\sH^2\to\cH^2$ be a boundary relation for
$S^*$ which satisfies the conditions (B1)--(B3), let
$\gamma(\lambda)$ be the corresponding $\gamma$-field, let
$\cH=\cH_1\oplus\cH_2$, decompose the corresponding Weyl function
$M(\lambda)$ as in \eqref{decom2}, and let $T\in [\cH_2,\cH_1]$.
Then:
\begin{enumerate}
\def\labelenumi{\rm (\roman{enumi})}

\item
the linear relation $H_T$ defined by
\begin{equation}
\label{sT}
 H_T=\left\{\, \wh f \in S^*:\, \{\wh f,\wh h\} \in\Gamma,\,\,
 h=0,\,  h_2'=-T^*h_1' \,\right\},
\end{equation}
is closed and symmetric in $\sH$ and has equal defect numbers;

\item the linear relation $\Gamma_T:\sH^2\to\cH_2^2$ given by
\begin{equation}
\label{BRT}
 \Gamma_T=\left\{\,
         \left\{\wh f,\begin{pmatrix} h_2 \\T^*h_1'+h_2' \end{pmatrix}\right\}:\,
         \left\{\wh f,\wh h \right\}\in\Gamma,\, h_1=Th_2 \,\right\},
\end{equation}
is a boundary relation for $H_T^*$ which satisfies the conditions
(B1)--(B3);

\item
the domain of $\Gamma_T$ is given by
\begin{equation}
\label{sT*}
 \dom \Gamma_T=\left\{\,\wh f \in S^*:\, \{\wh f,\wh h\} \in\Gamma, \,\,
               h_1=Th_2 \,\right\};
\end{equation}
\item
the $\gamma$-field $\gamma_T(\lambda):\cH_2\to \sH$ corresponding to
the boundary relation $\Gamma_T$ is given by
\begin{equation}
\label{dfT}
 \gamma_T(\lambda)=\gamma_1(\lambda)T+\gamma_2(\lambda),
\end{equation}
where
$\gamma(\lambda)=\begin{pmatrix}\gamma_1(\lambda)&\gamma_2(\lambda)\end{pmatrix}:\,
\cH_1\oplus\cH_2\to\sH$ is decomposed according to
$\cH=\cH_1\oplus\cH_2$;
\item
the Weyl function $M_T(\lambda)$ associated to $\Gamma_T$ is of the
form
\begin{equation}
\label{mT0}
 M_T(\lambda)
 =T^*M_{11}(\lambda)T+T^*M_{12}(\lambda)+M_{21}(\lambda)T+M_{22}(\lambda).
\end{equation}
\end{enumerate}
\end{proposition}
\begin{proof}
Define the operator $G\in[\cH]$, where $\cH=\cH_1\oplus\cH_2$, and
the operator $W\in[\cH\oplus\cH]$ via the block formulas
\[
G=\begin{pmatrix} I & T \\ 0 & I \end{pmatrix}, \quad
W=\begin{pmatrix} G^{-1} & 0 \\ 0 & G^* \end{pmatrix},
\]
respectively. Then $G$ is invertible, $G^{-1}\in[\cH]$, and $W$ is
$J_\cH$-unitary in $\cH^2=\cH\oplus\cH$. According to
Proposition~\ref{linprop} the product $\wt\Gamma=W\Gamma:\,
\sH^2\to\cH^2$ given by \eqref{lintrans} is a $J$-unitary relation
which satisfies the properties (B1)--(B3). Moreover, according to
\eqref{gmlin} the $\gamma$-field and the Weyl function associated to
$\wt\Gamma$ are given by
\begin{equation}
\label{gCT}
 \wt \gamma(\lambda)h=\gamma(\lambda)Gh
 =\gamma_1(\lambda)(h_1+Th_2)+\gamma_2(\lambda)h_2,
\end{equation}
and
\begin{equation}
\label{mCT}
 \wt M(\lambda) 
 =\begin{pmatrix}
   M_{11}(\lambda) & M_{11}(\lambda)T+M_{12}(\lambda) \\
   T^*M_{11}(\lambda)+M_{21}(\lambda)
   & T^*M_{11}(\lambda)T+T^*M_{12}(\lambda)+M_{21}(\lambda)T+M_{22}(\lambda)
  \end{pmatrix},
\end{equation}
respectively. Since
\[
 G^{-1}h=\begin{pmatrix}h_1-Th_2 \\ h_2\end{pmatrix},
 \quad
 G^*h'=\begin{pmatrix}h_1'\\T^*h_1'+h_2' \end{pmatrix},
\]
it follows from Corollary~\ref{cblock} that $H_T$ in \eqref{sT} is a
closed symmetric relation in $\sH$ and that
$\Gamma_T:\sH^2\to\cH_2^2$ defined by \eqref{BRT} is a boundary
relation for $H_T^*$ which satisfies the conditions (B1)--(B3).
Moreover, the formulas for the $\gamma$-field and the Weyl function
in \eqref{dfT} and \eqref{mT0} are obtained by applying
Corollary~\ref{cblock} to the formulas \eqref{gCT} and \eqref{mCT}.

The formula \eqref{sT*} is immediate from the description of
$\Gamma_T$ in \eqref{BRT}.
\end{proof}
\begin{corollary}
\label{summ} Let $S_j$ be symmetric operators in Hilbert spaces
$\sH_j$ and let $\Gamma^{(j)}:\sH_j^2\to\cH^2$ be  boundary
relations for $S_j^*$ which satisfy the conditions (B1)--(B3), and
let 
$M_j(\lambda)$ be the corresponding Weyl functions of $S_j$, j=1,2.
Then:
\begin{enumerate}
\def\labelenumi{\rm (\roman{enumi})}
\item the linear relation
\begin{equation}
\label{h3}
 H^{(3)}=\left\{\, \wh f=\wh f_1\oplus\wh f_2 
 :\, \left\{\wh f_1,
           \begin{pmatrix} 0 \\ h_1 \end{pmatrix}\right\} \in\Gamma^{(1)}, \,\left\{\wh f_2,
           \begin{pmatrix} 0 \\ -h_1 \end{pmatrix}\right\} \in\Gamma^{(2)} \,
              \mbox{ for some } h_1\in \cH \,\right\},
\end{equation}
is closed and symmetric in $\sH$ and has equal defect numbers;

\item the linear relation $\Gamma^{(3)}:\  \sH^2\to\cH^2$ given by
\begin{equation}\label{G3}
 \Gamma^{(3)}:=  
 \left\{\left\{\wh f_1\oplus\wh f_2 ,\begin{pmatrix} h \\ h_1+h_2
 \end{pmatrix}\right\}:\,\left\{\wh f_1,
           \begin{pmatrix} h \\ h_1 \end{pmatrix}\right\} \in\Gamma^{(1)}, \,\left\{\wh f_2,
           \begin{pmatrix} h \\ h_2 \end{pmatrix}\right\} \in\Gamma^{(2)} \,\,
               h,h_1,h_2 \in \cH \,\right\},
\end{equation}
is a boundary relation for $H^*$ which satisfies the conditions
(B1)--(B3);

\item the corresponding  Weyl function $M(\lambda)$ associated to $\Gamma^{(3)}$ is
\begin{equation}
\label{Msumm} M(\lambda)=
  M_1(\lambda)+M_{2}(\lambda),
  \quad \lambda\in\cmr.
\end{equation}
\end{enumerate}
\end{corollary}
\begin{proof}
To prove the statements (i)-(iii) it is enough to apply
Proposition~\ref{block2} to the boundary relation
\[
 \Gamma^{(1)}\oplus\Gamma^{(2)}:=
 \left\{\left\{\wh f_1\oplus\wh f_2 ,\wh h_1\oplus\wh h_2 \right\}:
           \left\{\wh f_1,\wh h_1 \right\} \in\Gamma^{(1)}, \,\left\{\wh f_2,\wh h_2 \right\} \in\Gamma^{(2)} \,\,
              \mbox{ for some } h_1,h_2 \in \cH \,\right\},
\]
for $S_1^*\oplus S_2^*$ with the corresponding Weyl function
$M(\lambda)=\mbox{diag }(M_1(\lambda),M_{2}(\lambda))$, setting
there $T=I_\cH$.
\end{proof}

\section{Orthogonal couplings}

\subsection{Orthogonal coupling and boundary relations}

Let $\sH_1$ and $\sH_2$ be arbitrary Hilbert spaces and let $\wt A$
be a selfadjoint linear relation in the  orthogonal sum $\wt
\sH=\sH_1\oplus\sH_2$. Then the formula
\begin{equation}\label{ST}
 S_j =\wt A \cap \sH_j^2,
 \quad
 T_j =\left\{\, \begin{pmatrix}P_{j} \varphi \\ P_{j} \varphi'\end{pmatrix}:\,
         \begin{pmatrix}\varphi\\ \varphi'\end{pmatrix}\in \wt A
         \,\right\},
\end{equation}
defines closed symmetric linear relations $S_1$ and $S_2$, and not
necessarily closed linear relations $T_1$ and $T_2$, in $\sH_1$ and
$\sH_2$, respectively. The relation $\wt A$ can be interpreted as a
selfadjoint extension of the orthogonal sum $S_1\oplus S_2$. It is
called the {\it orthogonal coupling} of $S_1$ and $T_2$ (or of $T_1$
and $S_2$), see~\cite{Str3}. The selfadjoint relation $\wt A$ is
said to be \textit{minimal} with respect to the Hilbert space
$\sH_j$ ($j$ is fixed, j=1,2) if
\begin{equation}
\label{minimal}
 \sH_1\oplus\sH_2=
 \cspan \left\{\,\sH_j+(\wt A-\lambda)^{-1}\sH_j:\,
 \lambda\in \rho(\wt A)\,\right\}.
\end{equation}
Associate with $T_j$ the eigenspaces as in \eqref{defect0},
\eqref{defectT},
\begin{equation}
\label{defect}
 \sN_\lambda(T_j)=\ker (T_j-\lambda),
\quad
 \wh \sN_\lambda(T_j)
 =\left\{\,  \begin{pmatrix}f\\ \lambda f\end{pmatrix} \in T_j
 :\, f\in\sN_\lambda(T_j)\,\right\}.
\end{equation}
Observe that $S_2$ is connected to $\wt S=\mul\Gamma$ in
\eqref{Stilde} via $S_2=-\wt S$, cf. \eqref{awig}. Moreover,
according to \cite[Lemma~2.14]{Trans} $\sN_\lambda(T_j)$ is dense in
$\sN_{\lambda}(S_j^*)$ for all $\lambda\in\cmr$, $j=1,2$.


\begin{lemma} \cite{Trans} \label{ran12} Let $\wt A$ be a selfadjoint linear
relation in $\wt \sH=\sH_1 \oplus \sH_2$, and let the linear
relations $S_j$ and $T_j$, $j=1,2$, be defined by~\eqref{ST}. Then:
\begin{enumerate}
\def\labelenumi{\rm (\roman{enumi})}
\item
$ \sN_\lambda(T_1)=P_1(\wt A-\lambda)^{-1}\, \sH_2$,
$\sN_\lambda(T_2) =P_2(\wt A-\lambda)^{-1}\, \sH_1$;
\item
$\sN_\lambda(T_j)$ is dense in $\sN_{\lambda}(S_j^*)$ for all
$\lambda\in\cmr$, $j=1,2$;
\item
The defect numbers of $S_1$ and $-S_2$ coincide:
$n_\pm(S_1)=n_\mp(S_2)$;
\item
$\wt A$ is minimal with respect to $\sH_1$ (resp. $\sH_2$) if and
only if $S_2$ (resp. $S_1$) is simple.
\end{enumerate}
\end{lemma}

The main transform $\wt A=\cJ(\Gamma)$ of a boundary relation
$\Gamma$ defined by \eqref{awig} can be treated as an orthogonal
coupling of symmetric operators $A$ and $\wt S=-\mul\Gamma$. Then
the first statement of the following proposition is just a
reformulation of Proposition~\ref{uksi}.
\begin{proposition}
\label{thm0} Let $\Gamma$ be a subspace in $\sH^2\times\cH^2$ and
let $S=\ker\Gamma$. Then $\Gamma$ is a  boundary relation for $S^*$
if and only if $\wt A=\cJ(\Gamma)$ is a selfadjoint relation in
$\sH\oplus\cH$. In this case the boundary relation $\Gamma$ is
minimal if and only if $\wt A=\cJ(\Gamma)$ is a minimal selfadjoint
extension of $-\wt S=-\mul\Gamma$.
\end{proposition}
\begin{proof}
To  prove the second statement let us mention first that the
boundary relation $\Gamma:\sH^2\to\cH^2$ for $S_2^*$ is minimal if
and only if the symmetric linear relation $S$ is simple, since
$\sN_\lambda(T)$ are dense in $\sN_\lambda(S^*)$ (see
Lemma~\ref{ran12}, (ii)). Combining this with the statement (iv) of
Lemma~\ref{ran12} one proves that $\Gamma$ is minimal if and only if
$\wt A=\cJ(\Gamma)$ is a minimal selfadjoint extension of $-\wt
S=-\mul\Gamma$.
\end{proof}

\subsection{Induced boundary relation}
Let $A$ be a symmetric operator in the Hilbert space $\sH_1$  and
let $\Pi=\{\cH,\Gamma_0,\Gamma_1\}$ be a boundary triplet for $A^*$.
Let $\wt A$ be a selfadjoint extension of $A$ in the Hilbert space
$\sH_1 \oplus \sH_2$ and define the linear relations  $S_2$ and
$T_2$. There is a natural way to define a boundary relation for
$S_2$ in the Hilbert space $\sH_2$ with corresponding Weyl family.

\begin{theorem}
\label{tmtau} Let $A$ be a symmetric operator in $\sH_1$ with equal
defect numbers and let $\Pi=\{\cH,\Gamma_0,\Gamma_1\}$ be an
ordinary boundary triplet for $A^*$. Then:
\begin{enumerate}
\def\labelenumi{\rm (\roman{enumi})}
\item
If $\wt A=\wt A^*$ is a minimal selfadjoint exit space extension of
$A$ in $\wt\sH=\sH_1\oplus\sH_2$ and $S_2$ is defined by \eqref{ST},
then the linear relation $\chi:\sH_2^2\to\cH^2$ defined
by
\begin{equation}
\label{Chi} \chi=\left\{\left\{\wh f_2,\begin{pmatrix}\Gamma_0\wh
f_1\\-\Gamma_1\wh f_1\end{pmatrix}\right\}: \,\wh f_1\oplus \wh
f_2\in \wt A,\,\wh f_1\in A^*,\,\wh f_2\in T_2 \right\}
\end{equation}
is a minimal boundary relation for $S_2^*$.
\item
If $S_2$ is a simple symmetric operator in $\sH_2$ and
$\chi:\sH_2^2\to\cH^2$ is a minimal boundary relation for $S_2^*$,
then the linear relation $\wt A$ defined by
\begin{equation}
\label{coupl}
   \wt A=\left\{\,\wh f_1\oplus\wh f_2\in A^*\oplus S_2^*:\,
\left\{ \wh f_2, \begin{pmatrix}\Gamma_0\wh f_1\\-\Gamma_1\wh
f_1\end{pmatrix}\right\} \in\chi\,\right\}
\end{equation}
is a minimal selfadjoint extension of $A$ which satisfies $\wt A\cap
\sH_2^2=S_2$.
\end{enumerate}
\end{theorem}

\begin{proof}
(i) Let $\wt A$ be a selfadjoint extension of $A$ in the Hilbert
space $\sH_1 \oplus \sH_2$ and define the linear relation
\begin{equation}
\label{Delta} \Delta:=\cJ^{-1}\circ \wt A
=\left\{\left\{\begin{pmatrix}f_1\\f_1'\end{pmatrix},\begin{pmatrix}f_2\\-f_2'\end{pmatrix}\right\}:
\,\{f_1\oplus f_2,f_1'\oplus f_2'\}\in\wt
A,\,f_j,f_j'\in\sH_j,\,j=1,2\right\}.
\end{equation}
It follows from Proposition~\ref{thm0} that $\Delta$ is a unitary
relation from $(\sH_2^2,J_{\sH_2})$ to $(\sH_1^2,J_{\sH_1})$ with
\[
\dom\Delta=T_1,\quad \ran \Delta=-T_2.
\]
Let $\Pi=\{\cH,\Gamma_0,\Gamma_1\}$ be an ordinary boundary triplet
for $A^*$. Since
\[
 \ran \Delta^{-1}=T_1 \subset \dom \Gamma, \quad \ran \Gamma=\cH^2,
\]
the composition $\chi_-$ of the unitary relations $\Delta^{-1}$ and
$\Gamma$
\begin{equation}\label{chi-}
\chi_-=\Gamma \circ\Delta^{-1} =\left\{\left\{
\begin{pmatrix} f_2\\-f_2'\end{pmatrix}, \Gamma\wh
f_1\right\}: \,\wh f_1\oplus \wh f_2\in \wt A,\,\wh f_1\in A^*,\,\wh
f_2\in T_2 \right\}.
\end{equation}
 is a unitary relation from $(\sH_2^2,J_{\sH_2})$ to
$(\cH^2,J_{\cH})$ with $\dom \chi_-=-T_2$, cf.
Theorem~\ref{product}. Changing signs in the second components of
$\chi_-$ gives the linear relation $\chi $  of the form~\eqref{Chi}
 with $\dom \chi=T_2$.
Since $\clos T_2=S_2^*$ it follows that $\chi:\sH_2^2\to\cH^2$ is a
boundary relation for $S_2^*$.

 To complete the proof of (i) it remains to prove that the boundary
relation $\chi:\sH_2^2\to\cH^2$ for $S_2^*$ is minimal. Since
$\sN_\lambda(T_2)$ are dense in $\sN_\lambda(S_2^*)$ (see
Lemma~\ref{ran12}, (ii)) the latter is equivalent to simplicity of
symmetric linear relation $S_2$. But as was shown in
Lemma~\ref{ran12} (iv), $S_2$ is simple  if and only if  $\wt A$ is
a minimal selfadjoint extension of $A$.

(ii) Let $\chi:\sH_2^2\to\cH^2$ be a boundary relation for $S_2^*$.
Then
\begin{equation}\label{chi-2}
\chi_- =\left\{\left\{
\begin{pmatrix} f_2\\-f_2'\end{pmatrix},
\begin{pmatrix} h\\-h'\end{pmatrix}\right\}:
\,\left\{
\begin{pmatrix} f_2\\f_2'\end{pmatrix},
\begin{pmatrix} h\\h'\end{pmatrix}\right\}\in\chi\right\}.
\end{equation}
is a boundary relation for $-S_2^*$. Since
\[
\ran\chi_-\subset\cH^2=\dom\Gamma^{-1},\quad
\ran\Gamma^{-1}=A^*=A^{[\perp]}=\left(\mul\Gamma^{-1}\right)^{[\perp]}
\]
it follows from Theorem~\ref{product} that the linear relation
\begin{equation}
\label{Delta-1}
 \Delta^{-1}:=\Gamma^{-1}\circ \chi_-=\left\{\left\{
 \begin{pmatrix}f_2\\-f_2'\end{pmatrix},\wh f_1\right\}:
\,
\left\{\begin{pmatrix}f_2\\f_2'\end{pmatrix},\begin{pmatrix}\Gamma_0\wh
f_1\\-\Gamma_1\wh f_1\end{pmatrix}\right\}\in\chi,\,\wh f_1\in
A^*\,\right\}
\end{equation}
is a unitary relation from $(\sH_2^2,J_{\sH_2})$ to
$(\sH_1^2,J_{\sH_1})$. Therefore, the linear relation
\begin{equation}
\label{Delta1}
 \Delta=\left\{\left\{\wh f_1,
 \begin{pmatrix}f_2\\-f_2'\end{pmatrix}\right\}:
\,
\left\{\begin{pmatrix}f_2\\f_2'\end{pmatrix},\begin{pmatrix}\Gamma_0\wh
f_1\\-\Gamma_1\wh f_1\end{pmatrix}\right\}\in\chi,\,\wh f_1\in
A^*\,\right\}
\end{equation}
is a unitary relation from $(\sH_1^2,J_{\sH_1})$ to
$(\sH_2^2,J_{\sH_2})$. Due to Theorem~\ref{product} the
composition $\Delta=\Gamma_-^{-1}\circ \chi$ of linear relations
$\chi$ and $\Gamma_-^{-1}$ is a unitary relation from
$(\sH_2^2,J_{\sH_2})$ to $(\sH_1^2,J_{\sH_1})$. Applying the
transform $\cJ$ to the linear relation $\Delta$ one obtains by
Proposition~\ref{thm0} a selfadjoint extension $\wt A$ of $A$
given by~\eqref{coupl}. Now ~\eqref{ST} is implied
by~\eqref{coupl}.

If $\chi:\sH_2^2\to\cH^2$ is a minimal boundary relation for
$S_2^*$, then the minimality of $\wt A$ with respect to $\sH_1$ is
implied by the same reasons as in (i).
           \end{proof}

\begin{proposition}\label{tauW}
Let the assumptions of Theorem~\ref{tmtau} be satisfied. Then the
family $\tau(\lambda)$ defined by
\begin{equation}\label{tau}
\tau(\lambda)=\left\{\left\{\Gamma_0\wh f_1,-\Gamma_1\wh
f_1\right\}\,:\ \  \wh f_1=P_{\sH_1}\wh f,\ \wh f\in\wt A,\
f'-\lambda f\in \sH_1\right\},
\end{equation}
is the Weyl family of $S_2$ corresponding to the boundary relation
$\chi: \ \ \sH_2^2\to\cH^2$.
\end{proposition}

\begin{proof}
Let $\wh f_2=\{f_2,f_2'\}\in\wh\sN_\lambda(T_2)$. Then it follows
from~\eqref{ST} that there are vectors $f_1,f_1'\in\sH_1$ such that
$\wh f=\{f,f'\}=\wh f_1\oplus\wh f_2\in\wt A$ and $f'-\lambda
f\in\sH_1$. Hence by~\eqref{coupl} one obtains
\begin{equation}\label{coupl2}
\left\{\wh f_2,\begin{pmatrix}\Gamma_0\wh f_1\\-\Gamma_1\wh
f_1\end{pmatrix}\right\} \in\chi.
\end{equation}
Since $\wh f_2\in\wh\sN_\lambda(T_2)$ this shows that
$\{\Gamma_0\wh
f_1,-\Gamma_1\wh f_1\}$ belongs to the Weyl family $M_\chi(\lambda)$
of $S_2$ corresponding to the boundary triplet $\chi$. This proves
the inclusion $M_\chi(\lambda)\subset\tau(\lambda)$,
$\lambda\in\dC_+\cup\dC_-$.

Conversely, if $\wh f_1$, $\wh f$ satisfy the
conditions~\eqref{tau}, then $\wh f_2=P_{\sH_2}\wh f$ belongs to
$\wh\sN_\lambda(T_2)$. Due to~\eqref{coupl2} one obtains
$\{\Gamma_0\wh f_1,-\Gamma_1\wh f_1\}\in M_\chi(\lambda)$, which
proves the inclusion $\tau(\lambda)\subset M_\chi(\lambda)$.
\end{proof}

Consider some examples of couplings of differential operators both
single-valued and multi-valued.

\begin{example}
Let $A$ be the symmetric differential operator in $L_2[0,1]$
associated with the differential expression $-D^2$, whose domain of
definition is given by
\[
\dom A=\left\{\,f\in C^1[0,1]:\, f'\in AC[0,1], f''\in L_2[0,1],
f(0+)=f'(0+)=f(1)=0\,\right\}.
\]
Let $S_2$ be a symmetric   differential operator  $-D^2$ on the
interval $[-1,0]$, whose domain of definition is given by
\[
 \dom S_2=\left\{\, f\in C^1[-1,0]:\, f'\in AC[-1,0], f''\in
L_2[-1,0], f(0-)=f'(0-)=f(-1)=0 \,\right\}.
\]
Then the boundary conditions~\eqref{coupl} take the form
\[
f(0+)=f(0-),\quad f'(0+)=f'(0-),
\]
and determine a selfadjoint operator  $\wt A$ in $L_2[-1,1]$
associated with the differential expression $-D^2$ and the boundary
conditions
\[
f(1)=0,\quad f(-1)=0.
\]
\end{example}

\begin{example}
Let $A$ be a minimal differential operator in $L_2[0,1]$ associated
with the differential expression $-D^2$. The domain of $A$ is
characterized by the following conditions
\[
\dom A=\{f\in C^1[0,1]:\,f'\in AC[0,1], f''\in L_2[0,1], f(0+)=f'(0+)=f(1)=f'(1)=0\}.
\]
Let the boundary triplet $\{\dC,\Gamma_0,\Gamma_1\}$ is given by
\[
\Gamma_0 f=\begin{pmatrix}f(0+)\\f(1)\end{pmatrix},\quad
\Gamma_1 f=\begin{pmatrix}f'(0+)\\-f'(1)\end{pmatrix}.
\]
Consider a minimal differential operator $S_2$ generated by the differential expression
$-D^2$ on the interval $[-1,0]$ and let the boundary triplet $\{\dC,\chi_0,\chi_1\}$ for $S_2^*$
is given by
\[
\chi_0 f=\begin{pmatrix}f(0-)\\f(-1)\end{pmatrix},\quad
\chi_1 f=\begin{pmatrix}-f'(0-)\\f'(-1)\end{pmatrix}.
\]
Then the boundary conditions~\eqref{coupl} take the form
\[
f(0+)=f(0-),\quad f'(0+)=f'(0-),\quad f(1)=f(-1),\quad f'(1)=f'(-1)
\]
and determine a selfadjoint operator  $\wt A$ in $L_2[-1,1]$ associated with the differential
expression $-D^2$ and the periodic boundary conditions
\[
f(1)=f(-1),\quad f'(1)=f'(-1).
\]
\end{example}

Theorem~\ref{tmtau} establishes a one-to-one correspondence between
all minimal with respect to $\sH_1$ exit space selfadjoint
extensions of $A$ and all minimal boundary relations
$\chi:\sH_2^2\to\cH^2$ with a fixed space $\cH$. Since minimal
boundary relations are uniquely determined by their Weyl families,
one can consider the correspondence established in
Theorem~\ref{tmtau} as a one-to-one correspondence between all
minimal  exit space selfadjoint extensions of $A$ and all Nevanlinna
families $\tau(\cdot)\in\wt R(\cH)$. This correspondence can be
written  explicitly in terms of generalized resolvents (see Section
4).
\begin{proposition}
Let under the assumptions of Proposition~\ref{tmtau}
$\tau=\{\Phi,\Psi\}$ be the Weyl family of the operator $S_2$
corresponding to the GBT~\eqref{Chi}. Then:
\begin{equation}\label{factor}
\dim S_1/A = \dim\ker{\mathsf N}_{\Phi,\Psi}.
\end{equation}
If, additionally, $S_1=A$ then $T_1\ne T_1^*(=S_1^*)$ if and only if
$\tau\in R_\cH^s\setminus R_\cH^u$, that is $0\in\sigma_c({\mathsf
N}_{\Phi,\Psi}(\lambda,\lambda))$ for each $\lambda\in\dC_+$.
\end{proposition}
\begin{proof}It follows from~\eqref{Chi} that
\[
\left(%
\begin{array}{c}
  h \\
  h' \\
\end{array}%
\right)\in\mul\chi \Longleftrightarrow
\left(%
\begin{array}{c}
  h \\
  -h' \\
\end{array}%
\right)=\Gamma \wh f, \mbox{ where }\wh f\in\wt A\cap\sH_1^2=S_1.
\]
Since $\Gamma$ is an isomorphism between linear spaces $\cH^2$ and
$A^*/A$ this implies that $\mul \chi$ and $S_1/A$ are isomorphic.
Therefore
 \[
\dim S_1/A = \dim\mul \chi.
\]
Making use of the equality~\eqref{mulGamma} one
obtains~\eqref{factor}.
\end{proof}
\begin{example}
Let $A$ be the same as in the previous example and let $S_2$  be a minimal differential
operator generated in $L_2(-\infty,1)$ by the differential expression $-D^2$.
Define a boundary relation $\chi: S_2^*\to\cH^2$ ($\cH=\dC$) for $S_2^*$
by the equality
\[
\chi=\left\{\left\{ \wh f,\mbox{col }(f(0-),c,-f'(0-),hc)\right\}:\,
\wh f\in S_2^*,\,c\in\dC\right\},
\]
where $h\in \dR$ is fixed. The equality~\eqref{coupl} take the form
\[
f(0+)=f(0-),\quad f'(0+)=f'(0-),\quad f(1)=c,\quad f'(1)=ch, \quad c\in\dC,
\]
and determine a selfadjoint operator  $\wt A$  generated in
$L_2(-\infty,1]$ by the differential expression $-D^2$ and the
boundary condition
\[
f'(1)=hf(1).
\]
\end{example}
The operator $S_1$ here is a restriction of $-D^2$ to the domain
\[
\dom S_1=\{f\in C^1[0,1]:\,f'\in AC[0,1], f''\in L_2[0,1],
f(0+)=f'(0+)=f'(1)-hf(1)=0\},
\]
and $\dim S_1/A =1$.
%

The boundary relations $\chi^j:\sH_2^2\to\cH^2$ (see
Theorem~\ref{tmtau}) are induced by the ordinary boundary triplets
$\Pi_j$ $(j=1,2)$ (see~\eqref{Chi}). Hence, due to Theorem
\ref{uneq} the connection between two Weyl families
$\tau_j(\lambda)$ corresponding to boundary relations
$\chi^j:\sH_2^2\to\cH^2$ can be explicitly expressed by means of the
transform $W$ which connects the Weyl functions $M_1(\lambda)$ and
$M_2(\lambda)$.

\begin{proposition}
Let the ordinary boundary triplets
$\Pi_j=\{\cH,\Gamma_0^j,\Gamma_1^j\}$ ($j=1,2$) for $A^*$ be
connected via the formula \eqref{Lftr} and let
$\chi^j:\sH_2^2\to\cH^2$ be boundary relations induced by the
ordinary boundary triplets $\Pi_j$ $(j=1,2)$ via the
formula~\eqref{Chi}. Then the boundary relations
$\chi^j:\sH_2^2\to\cH^2$ and the corresponding Weyl families
$\tau_j(\lambda)$ are connected by the formulas
\begin{equation}
\label{Lftr2} \chi^{(2)}=\wt W\chi^{(1)}, \quad
\tau_2(\lambda)=\wt W[\tau_1(\lambda)],\quad \wt
W=\begin{pmatrix}W_{00} & -W_{01}\\-W_{10} & W_{11}\end{pmatrix}.
\end{equation}
\end{proposition}
\begin{proof}
Let ${\wh f}={\wh f}_1\oplus {\wh f}_2 \in \wt A$.
Then one obtains
 \begin{equation}
\label{Gammaj} \{\Gamma_0^j\wh f_1,-\Gamma_1^j\wh f_1\}\in
\chi^j(\wh f_2)\quad(j=1,2).
\end{equation}
The formula~\eqref{Lftr2} is implied by~\eqref{Gammaj}
and the following equality
\[
\begin{pmatrix}\Gamma_0^2\wh f_1\\-\Gamma_1^2\wh f_1\end{pmatrix}=
\begin{pmatrix}I & 0\\ 0 & -I \end{pmatrix}
\begin{pmatrix}\Gamma_0^2\wh f_1\\\Gamma_1^2\wh f_1\end{pmatrix}=
\begin{pmatrix}I & 0\\ 0 & -I \end{pmatrix}
W\begin{pmatrix}\Gamma_0^1\wh f_1\\\Gamma_1^1\wh f_1\end{pmatrix}=
\wt W\begin{pmatrix}\Gamma_0^1\wh f_1\\-\Gamma_1^1\wh f_1\end{pmatrix}.
\]
The latter formula from~\eqref{Lftr2} is implied by
Lemma~\ref{Wlemma}.
\end{proof}


\subsection{The double Weyl function}

It is shown that associated with every selfadjoint extension $\wt A$
of $A$ there is a special boundary relation involving the linear
relation $A^*\oplus T_2$ and whose parameter space has double
dimension. The corresponding Weyl function of the operator $A\oplus
S_2$ can be written in the block form and as such is frequently
encountered in boundary-eigenvalue problems with boundary conditions
depending on the eigenvalue parameter (see e.g.~\cite{DLS},
\cite{DLS1}).

%

        \begin{theorem}\label{Couple}
Let $A$ be a symmetric operator in $\sH_1$ and let
$\Pi=\{{\cH}_1,\Gamma_0,\Gamma_1\}$ be a boundary triplet for $A^*$
with the Weyl function $M(\lambda)$. Let $S_2$ be a symmetric
operator in a Hilbert space $\sH_2$, let $\chi:\sH_2^2\mapsto\cH^2$
be a boundary relation for $S_2^*$ with the domain $\dom\chi=T_2$
and the Weyl family $\tau(\lambda)=\{\phi,\psi\}\in \wt R(\cH)$ and
let $\wt\sH=\sH_1\oplus\sH_2$. Then:
\begin{enumerate}
\def\labelenumi{\rm (\roman{enumi})}
\item the linear relation $\Gamma^{coupl}:\  \wt\sH^2\mapsto\cH_\Omega^2$ given by
\begin{equation}
\label{Btrip}
 \Gamma^{coupl}=\left\{\left\{\wh f_1\oplus\wh f_2,
   \begin{pmatrix}h'+\Gamma_1\wh f_1\\ h-\Gamma_0\wh
   f_1\end{pmatrix}
    \oplus
    \begin{pmatrix}-\Gamma_0\wh f_1\\ h'\end{pmatrix}\right\}:\,
\wh f_1\in A^*,\quad \left\{\wh f_2,\begin{pmatrix}h\\ h'\end{pmatrix}\right\}\in\chi\right\},
\end{equation}
is a boundary relation for $A^*\oplus S_2^*$, which satisfies (B1)--(B3) (see Proposition
\ref{GBTb});
\item   the corresponding Weyl function $M_{coupl}(\cdot)$  belongs to the class
$R{[\cH]}$ and is given by
\begin{equation}\label{Comega}
M_{coupl}(\lambda)=
 \begin{pmatrix}
-\Phi (\Psi +M \Phi )^{-1}   &  I- \Phi (\Psi +M \Phi )^{-1}  M  \\
\Psi (\Psi +M \Phi )^{-1}& \Psi (\Psi +M \Phi )^{-1}M
\end{pmatrix}.
\end{equation}
\end{enumerate}
\end{theorem}
\begin{proof}
(i) Clearly, the linear relation
       \begin{equation}\label{6.18A}
 \wt\Gamma=\left\{\left\{\wh f_1\oplus\wh f_2,
\left\{\begin{pmatrix}\Gamma_0\wh f_1 \\ -h'\end{pmatrix} ,
\begin{pmatrix}\Gamma_1\wh f_1 \\ h \end{pmatrix}\right\}\right\}:\, \quad
\wh f_1\in A^*,\quad \left\{\wh f_2,\begin{pmatrix}h\\
h'\end{pmatrix}\in\chi\right\}\right\},
         \end{equation}
forms a boundary relation for $A^*\oplus T_2$ and the corresponding
Weyl family is
\[
\wt\tau(\lambda)=M(\lambda)\oplus(-\tau(\lambda)^{-1})
=\{I\oplus(-\psi(\lambda)),M(\lambda)\oplus\phi(\lambda)\}.
\]
Let $W$ be a $J_{\cH_\Omega}$-unitary operator defined by
         \begin{equation}\label{6.19A}
W= \begin{pmatrix} W_{00} & I_{\cH^2} \\ -I_{\cH^2} & 0
\end{pmatrix}, \quad W_{00}= \begin{pmatrix} 0 & -I_{\cH} \\
-I_{\cH} & 0 \end{pmatrix}.
             \end{equation}
By Lemma~\ref{Wlemma} $\Gamma^{coupl} = W\wt\Gamma$ is a new
 boundary relation $\Pi^{coupl}$ for $A^*\oplus S_2^*$
whose Weyl family takes the form
      \begin{equation} \label{Gomega}
      M_{coupl}(\lambda) = W[\wt\tau(\lambda)]
  =\left\{\Omega_0(\lambda),
  \begin{pmatrix} -I & 0 \\ 0 & \psi(\lambda) \end{pmatrix}\right\},
\quad \text{ where }\quad  \Omega_0=
\begin{pmatrix}M(\lambda)& \psi(\lambda) \\
                                -I  & \phi(\lambda)\end{pmatrix}.
\end{equation}
Since $\Omega_0(\lambda)$ is invertible, see \cite[Proposition
A5]{DHMS04}, this implies that
$\Gamma_0^{\Omega}\uphar\left(\wh\sN_\lambda(A)\oplus\wh\sN_\lambda(T)\right)$
is a surjective mapping and by Proposition~\ref{GBTb}
$\Gamma^{coupl}$ is a boundary relation for $A\oplus S_2$ which
satisfies (B1)--(B3). This proves the statement~(i).

(ii) Setting $\omega:= (\psi +M \phi )^{-1}$ one easily derives from
\eqref{Gomega} the formula for the corresponding Weyl function
$M_{coupl}(\cdot)$:
  \[
                 \begin{split}
 M_{coupl}(\lambda)&=
\begin{pmatrix} -I & 0 \\ 0 & \psi(\lambda) \end{pmatrix}
\Omega_0(\lambda)^{-1}\\
&=
\begin{pmatrix} -I & 0 \\ 0 & \psi(\lambda) \end{pmatrix}
\begin{pmatrix}
 \phi(\lambda)\omega(\lambda) & \phi(\lambda)\omega(\lambda)M(\lambda)-I \\
 \omega(\lambda) & \omega(\lambda)M(\lambda)
  \end{pmatrix}    \\
&= \begin{pmatrix}
-\Phi (\Psi +M \Phi )^{-1}   &  I- \Phi (\Psi +M \Phi )^{-1}  M  \\
\Psi (\Psi +M \Phi )^{-1}& \Psi (\Psi +M \Phi )^{-1}M
\end{pmatrix}.
              \end{split}
              \]
Moreover, by Proposition~\ref{GBTb} $M_{coupl}(\cdot)\in R{[\cH]}$.
This gives (ii).
             \end{proof}
           \begin{remark}
(i)\  If the boundary relation $\chi$ in Theorem~\ref{tmtau} is single-valued then it can
be decomposed into a boundary triplet $\Pi''=\{\,\cH,\chi_0,\chi_1\,\}$,  where the
boundary operators $\chi_j$ are given by
\[
\chi_j=\pi_j\chi: \  T_2\to\cH, \quad j=0,1.
\]
In this case the  boundary relation  $\wt\Gamma$  of the form \eqref{6.18A}  becomes a
boundary triplet $\wt\Pi=\{\,{\cH}^2, \wt\Gamma_0, \wt\Gamma_1\,\}$ where
$$
 \wt\Gamma=  \begin{pmatrix}\wt\Gamma_0 \\  \wt\Gamma_1  \end{pmatrix},
 \qquad \wt\Gamma_0 = \begin{pmatrix}\Gamma_0 \\  -\chi_1
 \end{pmatrix},
 \quad \text{and} \quad \wt\Gamma_1
 =\begin{pmatrix}\Gamma_1 \\  \chi_0  \end{pmatrix},
   $$
and the equality~\eqref{coupl} takes the form
         \begin{equation} \label{coupl1A}
\wt A=  \ker(\wt\Gamma_1 - B \wt\Gamma_0) \quad  \text{with} \quad
B = \begin{pmatrix} 0 & I_{\cH} \\
I_{\cH} & 0
\end{pmatrix}.
         \end{equation}
In other words the coupling  $\wt A$ is determined by
           \begin{equation}\label{coupl1}
   \wt A=\left\{\,\wh f_1\oplus\wh f_2\in A^*\oplus T_2:\,
\Gamma_0\wh f_1-\chi_0\wh f_2=\Gamma_1\wh f_1+\chi_1\wh f_2=0\,\right\}.
      \end{equation}
In such a form  a construction of the coupling   $\wt A$ of two boundary triplets has
been introduced in \cite{DHMS1} under an additional assumption that
$\Pi''=\{\,\cH,\chi_0,\chi_1\,\}$ is an ordinary  boundary triplet.

(ii)\  Suppose that in Theorem  \ref{Couple}  the Nevanlinna family $\tau(\cdot)$ belongs
to $R^s(\cH)$. Then due to  \eqref{coupl1A} the  Weyl function 
corresponding  to the triplet $\Pi^{coupl}=\{\,{\cH}^2, \wt\Gamma_1
-B \wt\Gamma_0, -\wt\Gamma_0\,\}$ is $(B-\wt\tau (\cdot))^{-1}.$
Using the Frobenious formula we easily get
      \begin{equation}\label{Comega+}
(B-\wt\tau (\cdot))^{-1} = \left( \begin{pmatrix}
0   &   1  \\
1  &  0
\end{pmatrix} -
 \begin{pmatrix}
M    &  0  \\
0   &  -\tau^{-1}
\end{pmatrix}
\right) ^{-1}=
 \begin{pmatrix}
- (\tau +M )^{-1}   &   (\tau +M )^{-1}\tau  \\
\tau (\tau +M )^{-1}&  (\tau^{-1} +M^{-1} )^{-1}
\end{pmatrix} =  M_{coupl}.
       \end{equation}
Note that the matrix of  linear fractional transformation $\wt\tau (\cdot) \to (B-\wt\tau
(\cdot))^{-1}$ coincides  with  the block matrix $W$ determined by \eqref{6.19A}, that is
$(B-\wt\tau (\cdot))^{-1} = W[\wt\tau(\cdot)].$

Comparing \eqref{Comega} with  \eqref{Comega+} we see  that in this
case $M_{coupl}$ coincides  with the Weyl function corresponding to
the boundary triplet  $\Pi^{coupl}$.
Moreover, these reasonings, borrowed from  \cite{DHMS1}, explain the
appearance of the linear fractional transformation $W$ in formula
\eqref{Gomega}.

Under an additional assumption that $\Pi''=\{\,\cH,\chi_0,\chi_1\,\}$ is an ordinary
boundary triplet Theorem \ref{Couple} has been proved  in our paper \cite{DHMS1}.

(iii)\  In the case of finite defect numbers the function $
M_{coupl}(\cdot)$ appears for instance in the connection of
Sturm-Liouville operators and Hamiltonian systems  with
``$\lambda$-depending'' boundary conditions which are expressed by
means of a Nevanlinna pair $\{\Phi(\cdot),\Psi(\cdot)\}$ (equivalent
to $\tau(\cdot)$). In this case the function $ M_{coupl}(\cdot)$ is
known as the spectral matrix induced by the Nevanlinna pair
$\{\Phi(\cdot),\Psi(\cdot)\}$; cf. 
\cite{DLS}, \cite{DLS1}.
        \end{remark}

Finally, we demonstrate applicability of some results on
intermediate extensions from Section~\ref{Inter} when applied to the
Weyl function $M_{coupl}(\cdot)$ in Theorem \ref{Couple}. Namely, we
can prove that the diagonal elements of the matrix
$M_{coupl}(\lambda)$ are also Weyl families of some intermediate
extensions of the operator $A\oplus S_2$. In particular, this result
gives a geometric interpretation of the Nevanlinna function
$(\tau(\cdot)+M(\cdot))^{-1}$ appearing in the Krein-Naimark formula
for generalized resolvents (see \eqref{gres0}), as a Weyl function
of some intermediate extension. The importance of this result is
demonstrated in Section \ref{sec7}.
            \begin{theorem}
\label{GBTGB} Under the assumptions of Theorem~\ref{Couple}  the following statements hold:
\begin{enumerate}
\def\labelenumi{\rm (\roman{enumi})}
\item the linear relation
\begin{equation}
\label{couplH1}
   H^{(1)}=\left\{\,\wh f_1\oplus\wh f_2\in S_1^*\oplus T_2:\,
\Gamma_0\wh f_1=0, \left\{\wh f_2,\begin{pmatrix}0 \\
                                -\Gamma_1\wh f_1\end{pmatrix}\right\}
\in\chi\,\right\},
\end{equation}
is a closed symmetric linear relation in $\wt \sH=\sH_1\oplus\sH_2$;
\item the linear relation $\wt\Gamma^{(1)}:\wt\sH^2\mapsto\cH^2$ given by
\begin{equation}
\label{wtG}
   \Gamma^{(1)}=\left\{\,\wh f_1\oplus\wh f_2,\left(\begin{array}{c}
     \Gamma_1\wh f_1+h' \\
      -\Gamma_0\wh f_1\\
                                 \end{array}
                                 \right):\,
\left\{\wh f_2,\left(\begin{array}{c}
  \Gamma_0\wh f_1 \\
   h' \\
\end{array}\right)\right\}\in\chi\,\right\}.
\end{equation}
is a boundary relation for ${H^{(1)}}^*$ which satisfies the conditions (B1)--(B3);
\item
the Weyl function $M^{(1)}(\lambda)$ of $H^{(1)}$ corresponding ton
the boundary relation $\wt\Gamma^{(1)}$ is given by
\begin{equation}
\label{WeylH*}
M^{(1)}(\lambda)
=-\Phi(\lambda)(\Psi(\lambda)+M(\lambda)\Phi(\lambda))^{-1}.
\end{equation}
\item the linear relation
\begin{equation}
\label{couplH2}
   H^{(2)}=\left\{\,\wh f_1\oplus\wh f_2\in A^*\oplus T_2:\,
\Gamma_1\wh f_1=0, \left\{\wh f_2,\begin{pmatrix}-\Gamma_0\wh f_1 \\ 0
                                \end{pmatrix}\right\}
\in\chi\,\right\},
\end{equation}
is a closed symmetric linear relation in $\wt \sH=\sH_1\oplus\sH_2$;
\item the linear relation $\wt\Gamma^{(1)}:\wt\sH^2\mapsto\cH^2$ given by
\begin{equation}
\label{wtG2}
   \Gamma^{(2)}=\left\{\,\wh f_1\oplus\wh f_2,\left(\begin{array}{c}
     -\Gamma_0\wh f_1+h \\
      -\Gamma_1\wh f_1\\
                                 \end{array}
                                 \right):\,
\left\{\wh f_2,\left(\begin{array}{c}
        h\\
 -\Gamma_1\wh f_1  \\
\end{array}\right)\right\}\in\chi\,\right\}.
\end{equation}
is a boundary relation for ${H^{(2)}}^*$ which satisfies the conditions (B1)--(B3);
\item
the Weyl function $M^{(2)}(\lambda)$ of $H^{(2)}$ corresponding to
the boundary relation $\wt\Gamma^{(2)}$ is given by
\begin{equation}
\label{WeylH*2}
M^{(2)}(\lambda)
=\Psi(\lambda)(\Psi(\lambda)+M(\lambda)\Phi(\lambda))^{-1}M(\lambda).
\end{equation}
\end{enumerate}
\end{theorem}
\begin{proof}
Let us apply Proposition~\ref{coner} to the boundary relation
$\Gamma^{coupl}$ in Theorem~\ref{Couple}. Then the boundary
conditions in~\eqref{h1} take the form
\[
\Gamma_1\wh f_1+h'=-\Gamma_0\wh f_1+h=\Gamma_0\wh f_1=0,
\]
or, equivalently,
\[
\Gamma_1\wh f_1=-h',\quad\Gamma_0\wh f_1=h=0.
\]
Then it follows from Proposition~\ref{coner} that the linear
relation $H^{(1)}$ is a closed symmetric linear relation in $\wt
\sH$. The equality~\eqref{wtG} is implied by~\eqref{G1}. Due to
Proposition~\ref{coner} the Weyl function corresponding to the
boundary relation $\Gamma^{(1)}$ is the upper left corner of the
block matrix $M_{coupl}(\lambda)$.

Similarly, the statements (iv)-(vi) are implied by Theorem~\ref{Couple} and Corollary~\ref{cblock}.
\end{proof}

\section{Generalized resolvents and admissibility} \label{Sect7}
\subsection{Kre\u{\i}n's formula for generalized resolvents}
Let $A$ be a symmetric operator in a Hilbert space $\sH$ with equal
defect numbers. Let $\wt A$ be a selfadjoint extension of $A$ in a
Hilbert space $\wt\sH$ containing $\sH$ as a closed subspace. The
compression ${\mathbf R}_\lambda=P_\sH(\wt A-\lambda)^{-1}\uphar\sH$
of the resolvent of $\wt A$ to $\sH$ is said to be a generalized
resolvent of $A$.

Using the  coupling method, we easily obtain
the classical Kre\u{\i}n-Naimark formula, parametrizing all
generalized resolvents of $A$ by means of maximal dissipative
relations (Nevanlinna pairs) $\tau(\lambda).$ Namely, combining
Theorem~\ref{tmtau}, Proposition~\ref{tauW} and Theorem \ref{GBTNP}
we arrive at the following formula for generalized resolvents (in
the Straus form).
                    \begin{theorem}\label{tmtau2}
         Let $A$ be a symmetric operator in a Hilbert space $\sH$ and $n_+(A)= n_-(A)$. Let $\wt A$ be a
selfadjoint extension of $A$ in a Hilbert space $\wt \sH\supset\sH$
and let $\Pi=\{{\cH},\Gamma_0,\Gamma_1\}$ be a boundary triplet for
$A^*$. Then there is a unique Nevanlinna family $\tau(\lambda)\in\wt
R_\cH$ such that
\begin{equation}\label{Sres1}
P_\sH(\wt A-\lambda)^{-1}\uphar\sH=(\wt
A_{\tau(\lambda)}-\lambda)^{-1}.
\end{equation}

Moreover, for any $h\in \sH$,  vector $f_1=P_\sH(\wt
A-\lambda)^{-1}h$ is a solution of the ``boundary-value problem''
with spectral parameter $\tau(\lambda)$ in "boundary condition"
\begin{equation}\label{bpr1}
  \left\{
 \begin{array}{ll}
   &  f_1'-\lambda f_1=h,
 \quad \wh f_1=\{f_1,f_1'\}\in A^*,\\
   &  \{\Gamma_0\wh f_1,-\Gamma_1\wh f_1\}\in\tau(\lambda);
 \end{array}
 \right.
\end{equation}

Conversely, given $\tau(\lambda)\in\wt R_\cH$ there is a minimal
selfadjoint extension $\wt A$ of $A$ in a Hilbert space $\wt
\sH\supset\sH$ such that \eqref{Sres1} holds.
           \end{theorem}
\begin{proof}
(i) Let $\lambda\in\rho(\wt A)$ and let $h\in \sH$. Then there is a
vector $\wh f=\begin{pmatrix}f\\f'
\end{pmatrix}\in \wt A$
such that
\begin{equation}
\label{flfh} f'-\lambda f=h.
\end{equation}
Projecting of \eqref{flfh} to $\sH_1$ and $\sH_2$ gives the
following equations
\begin{equation}
\label{flfh12} f'_1-\lambda f_1=h,\quad f_2'-\lambda f_2=0,
\end{equation}
where $f_j=P_{\sH_j}f$, $f'_j=P_{\sH_j}f'$, $j=1,2$. It follows
from~\eqref{coupl} that
$\left\{\wh f_2,\begin{pmatrix}\Gamma_0 \wh f_1\\
-\Gamma_1 \wh f_1\end{pmatrix}\right\}\in\chi$, where $\wh
f_j=\begin{pmatrix}f_j\\f'_j\end{pmatrix}$. Since $\wh
f_2\in\sN_\lambda(T_2)$ this implies
\begin{equation}
\label{proptau}
\{\Gamma_0 \wh f_1, -\Gamma_1 \wh f_1\}
\in\tau(\lambda),
\end{equation}
where $\tau(\cdot)$ is the Weyl family of $S_2$ corresponding to the
boundary relation $\chi$. This proves the statement (i).

(ii)  Conversely, starting with $\tau(\cdot)$ and applying Theorem
\ref{GBTNP} we find a simple symmetric operator $S_2$ in $\sH_2$ and
a minimal boundary relation $\chi:\sH_2^2\to\cH^2$  for $S_2^*$ such
that the corresponding Weyl family is $\tau(\lambda).$  Then by
Theorem \ref{tmtau} the linear relation $\wt A$ in  $\wt
\sH=\sH_1\oplus\sH_2$ (a coupling of $T_1$ and $S_2$) defined by
\eqref{coupl} is an exit space selfadjoint extension of $A$ which
satisfies~\eqref{ST} and \eqref{Sres1} with some $\tau_1(\cdot)\in
\wt R(\cH)$ in place of $\tau(\cdot).$  By  Proposition~\ref{tauW}
$\tau_1(\cdot) = \tau(\cdot).$
                     \end{proof}

Combining Theorem \ref{tmtau2} with
formula \eqref{krein00} for canonical resolvents we arrive at the
following



       \begin{theorem}\label{krein}~{\rm (\cite{KL1})}.
       Let $A$ be a symmetric operator in ${\sH}$ with $n_+(A)= n_-(A)$, let $\Pi=\{{\cH}_1,\Gamma_0,\Gamma_1\}$
be a boundary triplet for $A^*$, and let $M(\cdot)$ and
$\gamma(\cdot)$ be the corresponding Weyl function and the
$\gamma$-field. Then the formula
\begin{equation} \label{gres0}
 {\mathbf R}_\lambda=(A_0-\lambda)^{-1}-\gamma(\lambda)(M(\lambda)
     +\tau(\lambda))^{-1}\gamma(\bar{\lambda})^*,\quad\lambda\in\rho(A_0)\cap \rho(\wt A)
\end{equation}
with $A_0=\ker\Gamma_0$ establishes a bijective correspondence
between the generalized resolvents ${\mathbf R}_\lambda$ of $A$ and
Nevanlinna families $\tau(\cdot)\in \wt R_\cH$.
      \end{theorem}
\begin{proof}
Let $\lambda\in\rho(A_0)$. According to Proposition~\ref{PropRez}
$\lambda\in\rho(A_{-\tau(\lambda)})$ if and only if
$0\in\rho(M(\lambda)+\tau(\lambda))$. In this case
(see~\eqref{krein00})
\begin{equation}
\label{Krein} (\wt
A_{-\tau(\lambda)}-\lambda)^{-1}=(A_0-\lambda)^{-1}-\gamma(\lambda)(M(\lambda)
     +\tau(\lambda))^{-1}\gamma(\bar{\lambda})^*.
\end{equation}
Now the statement follows from Theorem~\ref{tmtau2}.
\end{proof}

\begin{remark}
(i) Note that for "good"  $\tau(\cdot)$ Theorem \ref{krein} can
easily be derived from Theorem~\ref{Couple} and formula
\eqref{krein00} for canonical resolvents with double Weyl function
$M_{coupl}(\cdot)$ (see \eqref{Comega}). We explain the proof
confining ourself to the case $\tau(\cdot)\in R^u[\cH].$
By Theorem~\ref{GBTGB} and Proposition~\ref{OBTpr} there exists an
ordinary boundary triplet  $\{\cH,\chi_0,\chi_1\}$ for $S_2^*$ such
that the corresponding Weyl function is $\tau(\cdot).$ Consider  a
boundary triplet $\{\cH^2,\Gamma^\Omega_0,\Gamma^\Omega_1\}$ for
$A^*\oplus S_2^*$ of the form \eqref{Btrip}. The corresponding Weyl
function $M(\cdot)$ is of the form \eqref{Comega},
$M(\cdot)=\Omega(\cdot).$  Then $\wt A$ and $A_0\oplus A_1^2$
($A_1^{(2)}=\ker\chi_1$) are canonical selfadjoint  extensions of
$A\oplus S_2$ and the formula~\eqref{krein00} implies
\begin{equation}
\label{Krein2} (\wt
A-\lambda)^{-1}\begin{pmatrix}h_1\\h_2\end{pmatrix}
=\begin{pmatrix}(A_0-\lambda)^{-1}h_1\\(A_1^{(2)}-\lambda)^{-1}h_2\end{pmatrix}
-\begin{pmatrix}\gamma(\lambda)& 0\\0 &
\gamma^{(2)}(\lambda)\end{pmatrix} \Omega(\lambda)
\begin{pmatrix}\gamma(\bar\lambda)^*h_1\\ \gamma^{(2)}(\bar\lambda)^*h_2\end{pmatrix},
\end{equation}
where $\gamma^{(2)}(\lambda)$ is the $\gamma$-field corresponding to
the boundary triplet $ \{\cH,-\chi_1,\chi_0\}$ and
$A_1^{(2)}=\ker\chi_1$. Setting $h_2=0$ and applying the projection
$P_1$ onto $\sH_1$ to~\eqref{Krein2} we arrive at \eqref{gres0}.

(ii)  Note, that in fact, both formulas \eqref{Sres1} and
\eqref{gres0} are equivalent to each other and can easily be deduced
one from another (cf. \cite{MM2, DM2}).
            \end{remark}

\begin{remark}
The description of all generalized resolvents was originally given
in different forms by M.G.~Kre\u{\i}n~\cite{Kr0} and
M.A.~Naimark~\cite{Naj1}. It has been extended to  the case of
infinite indices by Saakyan (see \cite{KL1, DM1} and references
therein). Another description in a form close to \eqref{Sres1} was
given by A.V.~\v{S}traus~\cite{Str0}.
A connection of the Kre\u{\i}n-Naimark formula  with boundary
triplets has been discovered in \cite{DM1}, \cite{DM2}, \cite{MM2}.
Moreover, other proofs as well as generalizations of the
Kre\u{\i}n-Naimark formula for nondensely defined symmetric
operators can be found in \cite{DM2}, \cite{MM2}, \cite{DLS90},
\cite{LT}; see also the references therein.
\end{remark}

\subsection{Admissibility} \label{sec7}

In this section some new admissibility criteria will be given, which
guarantee that a generalized resolvent of a symmetric operator
corresponds to a selfadjoint operator extension. Their relation to
some other conditions which have been found earlier in \cite{DM2},
\cite{MM2},  \cite{LT} will be discussed.

Let $A$ be a symmetric operator in $\sH$ with equal defect numbers
$n_+(A)=n_-(A)<\infty$ and let $\Pi=\{\cH,\Gamma_0,\Gamma_1\}$ be a
boundary triplet for $A^*$. According to Theorem \ref{krein} the
generalized resolvents ${\mathbf R}_\lambda$ of $A$ are in
one-to-one correspondence with Nevanlinna families $\tau(\lambda)\in
\wt R_\cH$ via the Kre\u{\i}n-Naimark formula \eqref{gres0}. Let
$\wt A$ be a minimal selfadjoint extension of $A$ whose compressed
resolvent is equal to ${\mathbf R}_\lambda$. Then the family
$\tau(\lambda)$ associated to $\wt A$ via \eqref{gres0} is said to
be $\Pi$-admissible, if $\wt A$ is an operator extension of $A$,
i.e., if $\mul \wt A=\{0\}$.

The next theorem gives a general criterion for the $\Pi$-admissibility
of the family $\tau(\lambda)=\{\phi(\lambda),\psi(\lambda)\}$.

\begin{theorem}\label{admisthm}
Let $A$ be a (nondensely defined) closed symmetric operator in $\sH$
with equal defect numbers $n_+(A)=n_-(A)\le \infty$, let $\Pi=\{\cH,\Gamma_0,\Gamma_1\}$
be a boundary triple for $A^*$ with Weyl function $M(\lambda)$, and let $\tau(\lambda)=
\{\phi(\lambda),\psi(\lambda)\}$ be a Nevanlinna pair in $\cH$. Then:
\begin{enumerate}
\def\labelenumi{\rm (\roman{enumi})}

\item The pair $\{\phi(\lambda),\psi(\lambda)\}$
is $\Pi$-admissible if and only if the following two conditions
are satisfied:
\begin{equation}
\label{Adm1}
  w-\lim_{y \uparrow \infty}
  \frac{\phi(\lambda)(\psi(iy)+M(iy)\phi(\lambda))^{-1}}{y}=0
\end{equation}
and
\begin{equation}
\label{Adm2}
  \lim_{y \uparrow \infty}
  \frac{\psi(\lambda)(\psi(iy)+M(iy)\phi(\lambda))^{-1}M(\lambda)}{y}=0.
\end{equation}

\item If, in addition, $A_0=\ker \Gamma_0$ is an operator,
then the $\Pi$-admissibility of $\{\phi(\lambda),\psi(\lambda)\}$
is equivalent to the single condition \eqref{Adm1}.

\item If $A_1=\ker \Gamma_1$ is an operator,
then the $\Pi$-admissibility of $\{\phi(\lambda),\psi(\lambda)\}$
is equivalent to the single condition \eqref{Adm2}.
\end{enumerate}
\end{theorem}
\begin{proof}
(i) By Theorem \ref{GBTNP} there are a Hilbert space $\sH_2$, a
symmetric operator $S_2$ and a boundary relation
$\chi:\sH^2\to\cH^2$ whose Weyl family is
$\tau(\lambda)=\{\phi(\lambda),\psi(\lambda)\}$. Let the selfadjoint
extension $\wt A$ of $A\oplus S_2$ be as in Lemma~\ref{Couple}.
Moreover, by Lemma~\ref{Couple} the function $\Omega(\lambda)$ given
by \eqref{Comega}, is the Weyl function of $A\oplus S_2$
corresponding to the boundary relation
$\Gamma_\Omega:\wt\sH^2\to\cH^2_\Omega$ of the form~\eqref{Btrip}.
According to Proposition~\ref{MULA0T} the multivalued part of the
linear relation $\wt A$ is trivial if and only if
\begin{equation}
\label{Adm} w-\lim_{y\uparrow \infty}\frac{\Omega(iy)}{y}=0.
\end{equation}
Now it remains to note that~\eqref{Adm} is equivalent to
the conditions ~\eqref{Adm1}, ~\eqref{Adm2}.

(ii)
Assume that $A_0$ is an operator and consider the boundary relation
\[   \Gamma^{(1)}=\left\{\,\wh f_1\oplus\wh f_2,\left(\begin{array}{c}
     \Gamma_1\wh f_1+h' \\
      -\Gamma_0\wh f_1\\
                                 \end{array}
                                 \right):\,
\left\{\wh f_2,\left(\begin{array}{c}
  \Gamma_0\wh f_1 \\
   h' \\
\end{array}\right)\right\}\in\chi\,\right\}.
\]
for $(H^{(1)})^*$, where $H^{(1)}$ is a closed symmetric linear
relation in $\wt \sH=\sH_1\oplus\sH_2$ given by~\eqref{couplH1}. As
was shown in Theorem~\ref{GBTGB} the boundary relation
$\wt\Gamma^{(1)}:\wt\sH^2\mapsto\cH^2$
satisfies the conditions (B1)-(B3) and the corresponding Weyl
function of $H^{(1)}$ is $-(M(\lambda)+\tau(\lambda))^{-1}$. To see
that $H^{(1)}$ is an operator assume that $\wh f_1\oplus \wh
f_2=\{0,f_1'\}\oplus\{0,f_2'\}\in H^{(1)}$. Then $\wh f_1=0$, since
$\Gamma_0\wh f_1=0$ and $\mul A_0=\{0\}$. Using the last condition
in the definition of $H^{(1)}$ in \eqref{couplH1} one obtains $\{\wh
f_2,0\}\in\chi$, which due to Proposition~\ref{GBTb} implies $\wh
f_2\in S_2$. Since $S_2$ is an operator, it follows that $\wh
f_2=0$, and hence $H$ is also an operator. In view of
Proposition~\ref{MULA0T} $\wt A$ is an operator if and only if
\eqref{Adm1} holds.

(iii)
In the case where $A_1=\ker \Gamma_1$ is an operator one
can replace the boundary triplet $\Pi=\{\cH,\Gamma_0,\Gamma_1\}$
by $\wt \Pi=\{\cH,\Gamma_1,-\Gamma_0\}$.
Then the corresponding Weyl families are transformed to
$\wt M(\lambda)=-M(\lambda)^{-1}$ and
$\wt\tau(\lambda)=-\tau(\lambda)^{-1}$,
and the statement in the part (iii) is obtained from
the part (ii).
       \end{proof}
      \begin{remark}
Other approaches to the admissibility  problem have been proposed in
\cite{LT}, \cite{MM2} and  \cite{DM2}. Namely,   a  direct deduction
of Theorem \ref{admisthm} (ii) from Krein-Naimark formula has been
obtained  in \cite{MM2}. This proof is  more complicated than the
one proposed here. Furthermore, under the additional assumption that
$A_1=\ker \Gamma_1$ is an operator, another criterion of
admissibility (with rather complicated proof) has been obtained in
\cite{DM2}. This criterion is equivalent to the statement
Theorem~\ref{admisthm}(iii), while we don't know a direct proof of
their equivalence.

Another criterion of admissibility (without additional assumptions)
has been obtained in \cite{LT}. A connection of
Theorem~\ref{admisthm} with the Langer-Textorious result is
discussed in Section 7.3.
       \end{remark}
In the next proposition another admissibility criterion is obtained, when $\wt A$ is
viewed as an extension of the symmetric intermediate extension $H_T$ defined in
Proposition~\ref{block2}.

\begin{proposition}
\label{sum2} Let $A$ be a simple symmetric operator satisfying the
assumptions of Theorem~{\rm \ref{admisthm}}. Assume that $T \in
[{\cH}]$ and let $M_T(\lambda)$ be defined by
\begin{equation}
\label{MT}
\begin{split}
 M_T(\lambda)= & -T^*(M(\lambda)+\tau(\lambda))^{-1}T
    -T^*(M(\lambda)+\tau(\lambda))^{-1}\tau(\lambda) \\
    &-\tau(\lambda)(M(\lambda)+\tau(\lambda))^{-1}T
     +\tau(\lambda)(M(\lambda)+\tau(\lambda))^{-1}M(\lambda).
\end{split}
\end{equation}
Then for the extension $A^{(\tau)}\,(=\wt A)$
in \eqref{coupl1} to be an operator it is necessary and,
if $\wt A_{T^*}=\ker(\Gamma_1-T^*\Gamma_0)$ is an operator,
it is also sufficient that the following condition holds:
\begin{equation}
\label{mTlim0}
  \mbox{s}-\lim_{y\uparrow\infty}\frac{M_T(iy)}{y}=0;
\end{equation}
\end{proposition}
\begin{proof}
It follows from Proposition~\ref{block2} and Theorem~\ref{Couple}
that $M_T(\lambda)$ is the Weyl function of the linear relation
\begin{equation}
\label{HT}
 H_T=\left\{\, \wh f_1 \oplus \wh f_2
            \in A^*\oplus S_2^*:\,
     \Gamma_1\wh f_1+h'
     =\Gamma_0\wh f_1-h
     =\Gamma_1\wh f_1-T^*\Gamma_0\wh f_1=0,
\left\{\wh f_1,\begin{pmatrix} h\\h'\end{pmatrix}\right\}\in\Delta
      \,\right\}
\end{equation}
corresponding to the boundary relation
\begin{equation}
\label{BROmegaT}
 \Gamma^{T}=\left\{\left\{\wh f_1\oplus\wh f_2,
   \begin{pmatrix}-\Gamma_0\wh f_1+h\\ -T^*\Gamma_0\wh
   f_1+h'\end{pmatrix}\right\}:\,
\wh f_1\in A^*,\quad \left\{\wh f_2,\begin{pmatrix}h\\ h'\end{pmatrix}\right\}\in\Delta\right\},
\end{equation}

The necessity of the condition~\eqref{mTlim0} follows immediately
from (\ref{Adm1}) and (\ref{Adm2}) in Theorem \ref{admisthm}. To
prove the sufficiency let us show that the following implication
holds:
\begin{equation}
\label{MulATHT}
 \mul\wt A_{T^*}=\{0\}\Rightarrow \mul H_{T}=\{0\}.
\end{equation}
Indeed, if $\wh f=(\wh f_1,\wh f_2)^\top \in H_T$ and $\wh
f_i=\{0,f_i'\}$, $i=1,2$, then \eqref{HT} implies that $\Gamma_1\wh
f_1-T^*\Gamma_0\wh f_1=0$ and, hence, $\wh f_1\in \wt A_{T^*}$.
Since $\mul \wt A_{T^*}=\{0\}$, one obtains $\wh f_1=0$. Now it
follows from \eqref{HT} that $\{\wh f_2,0\}\in\Delta$. Thus, $\wh
f_2 \in S_2=\ker \Delta$, and consequently $\wh f_2=0$, since $S_2$
is a simple operator. This proves that, $\mul H_{T}=\{0\}$. Since
$A^{(\tau)}=\ker\Gamma_0^{T}$, it follows from \eqref{HT} that
$M_T(\lambda)$ is a Weyl function of the pair $(H_T,A^{(\tau)})$. By
Proposition~\ref{MULA0T} the condition \eqref{mTlim0} implies that
$A^{(\tau)}$ is an operator. This completes the proof.
\end{proof}

\subsection{The Langer-Textorius criterion}
In this subsection a new proof for the admissibility
criterion in~\cite{LT}  will be given.

Following~\cite{LT} introduce the operator function $Q_{\tau}$
with values in $[{\cH}]$ by
\begin{equation}
\label{Qtau}
 Q^{\tau}_{LT}(\lambda; z_0):=
  M(\lambda)-(M(\lambda)-M(z_0)^*)
    (M(\lambda)+\tau(\lambda))^{-1}(M(\lambda)-M(z_0)).
\end{equation}
The function $Q^{\tau}_{LT}(\lambda;z_0)$ is a $Q$-function
of a pair $(H_{LT},A^{(\tau)})$,
where $A^{(\tau)}$ is a minimal selfadjoint
exit space extension of $A$ corresponding to
$\tau(\lambda)$ in \eqref{gres0} and $H_{LT}$ is a symmetric
restriction of $A^{(\tau)}$, cf. \cite{KL2}, \cite{LT}.
In the following proposition the symmetric linear relation $H_{LT}$ is calculated explicitely. This allows to derive the Langer-Textorius criterion from Proposition~\ref{sum2}.

The next theorem specifies the operator $H_{LT}$ with the help of boundary operators.
\begin{proposition}
\label{LTsum}
Let the assumptions be as in Proposition~\ref{sum2}
and let $z_0\in \dC_+$ be fixed. Then:
\begin{enumerate}
\def\labelenumi{\rm (\roman{enumi})}

\item the linear relation $H_{LT}$ defined by
\begin{equation}
\label{HT2}
 H_{LT}=\left\{\, \wh f= \wh f_1 \oplus \wh f_2
            \in A^*\oplus S_2^*:\,
     \begin{array}{c}
       \Gamma_1\wh f_1+h'
     =\Gamma_0\wh f_1-h=0, \\
       \Gamma_1\wh f_1-M(z_0)^*\Gamma_0\wh f_1=0, \\
     \end{array}
     \left\{\wh f_1,\begin{pmatrix} h\\h'\end{pmatrix}\right\}\in\Delta
      \,\right\}
\end{equation}
is a closed symmetric operator;
\item a linear relation
\begin{equation}
\label{BROmegaLT}
 \Gamma^{LT}=\left\{\left\{\wh f_1\oplus\wh f_2,
   \begin{pmatrix}-\Gamma_0\wh f_1+h\\ -M(z_0)^*\Gamma_0\wh
   f_1+h'\end{pmatrix}\right\}:\,
\wh f_1\in A^*,\quad \left\{\wh f_2,\begin{pmatrix}h\\ h'\end{pmatrix}\right\}\in\Delta\right\},
\end{equation}
is a boundary relation for $H_{LT}^*$;

\item the Weyl function corresponding to $\Gamma^{LT}$ is given by
\begin{equation}
\label{6.36} M_{LT}(\lambda)=  Q^{\tau}_{LT}(\lambda; z_0)-2\RE M(z_0);
\end{equation}
\item
$\tau(\lambda)$ is admissible if and only if
\begin{equation}
\label{6.37}
 \mbox{s}-\lim_{y\uparrow \infty }\frac{M_{LT}(iy)}{y}=0.
\end{equation}
\end{enumerate}
\end{proposition}
\begin{proof}
As was shown in Proposition~\ref{sum2} $H_{LT}$ is a closed symmetric linear relation in $\sH_1\oplus\sH_2$. Moreover, the linear relation $\wt A_{M(z_0)^*}$ takes the form
\[
\wt A_{M(z_0)^*}=A\dotplus\wh\sN_{\bar z_0}.
\]
Since $\sN_{\bar z_0}\cap\dom A=\{0\}$, this implies that $\wt
A_{M(z_0)^*}$ is an operator. Now it follows from \eqref{MulATHT}
that $H_{LT}$   is an operator. The expressions \eqref{HT2},
 and \eqref{BROmegaLT} are obtained from
Proposition~\ref{block2} and the identities in \eqref{Btrip}.
\end{proof}

The functions in \eqref{Qtau} and \eqref{6.36}
are related by
\[
 Q_{LT}^{\tau}(\lambda;z_0)
  =M_{LT}(\lambda)+2\RE M(z_0).
\]
Therefore, Proposition~{\rm \ref{LTsum}~(iv)} yields
the following theorem in \cite{LT}.

\begin{theorem}
\label{Th6.9}
{\rm (\cite{LT})}
Let $z_0\in \dC_+.$
Then the minimal selfadjoint extension
$A^{(\tau)}$ of $A$ in Kre\u{\i}n's formula \eqref{gres0}
is an operator if and only if
\begin{equation}
\label{6.34}
 \lim_{y\uparrow \infty }
  \frac{(Q_{LT}^{\tau }(iy; z_0)h,h)}{y}=0,
\quad  h\in {\cH}.
\end{equation}
\end{theorem}


\end{document}